\input amstex
\magnification=\magstep1
\input epsf
\input amssym.def
\input amssym
\pageno=1
\baselineskip 14 pt
\def \pop#1{\vskip#1 \baselineskip}

\font\gr=cmbx12

\def \rad{\operatorname {rad}}

\def \et{\operatorname {et}}

\def \Fr{\operatorname {Fr}}
\def \Spec{\operatorname {Spec}}
\def \Spf{\operatorname {Spf}}
\def \lim{\operatorname {lim}}

\def \Sp  {\operatorname {Sp  }}

\def \fppf{\operatorname {fppf }}

\def \mod{\operatorname {mod }}

\def \card{\operatorname {card}}

\def \Spf{\operatorname {Spf}}

\def \gcd{\operatorname {gcd}}

\def \Vert{\operatorname {Vert}}
\def \Fr{\operatorname {Fr}}

\def \Deg{\operatorname {Deg}}
\def \SS-Deg{\operatorname {SS-Deg}}
\def \SNS-Deg{\operatorname {SNS-Deg}}
 \def \DS-Deg{\operatorname {DS-Deg}}
\def \DNS-Deg{\operatorname {DNS-Deg}}
\def \F{\operatorname {F}}
\def \Id{\operatorname {Id}}
\def \sup{\operatorname {sup}}
\def \ir{\operatorname {ir}}

\pop {4}
\par
\noindent                                          
\centerline {\bf \gr Cyclic $p$-groups and semi-stable reduction of curves} 

\par
\noindent 

\centerline {\bf \gr in equal characteristic $p>0$. I}

\pop {3}
\noindent                                          
\centerline {\bf \gr Mohamed Sa\"\i di}

\pop {4}
\par
\noindent
\centerline {\bf \gr  Abstract}
\par
In this paper we study the semi-stable reduction of $p$ and $p^2$-cyclic covers of curves in equal characteristic $p>0$.
The main tool we use is the classical Artin-Schreier-Witt theory for $p^n$-cyclic covers in characteristic $p$. Although the results of this paper
concern only the cases of degree $p$ and $p^2$-cyclic cover we develop the techniques and the framework in which the general cyclic case 
can be studied.

\pop {2}
\par
\noindent
{\bf \gr 0. Introduction.}\rm \ This paper is a further attempt to understand the degeneration of wildly ramified Galois covers of algebraic varieties.
Let $R$ be a complete discrete valuation ring with fraction field $K$ and residue field $k$ of characteristic $p>0$. Let $X$ be a smooth 
(or more generally semi-stable) $R$-scheme and let $f:Y\to X$ be a finite Galois cover with group $G$. The basic problem we are interested in 
is to understand the reduction of $Y$ and in particular its semi-stable reduction if it exists (e.g. case of curves). The case where the 
cardinality of $G$ is prime to $p$ is by now classical and well understood (cf. [SGA-1], and [S-2] for the case of semi-stable curves). Although
the case of higher dimensional semi-stable schemes is not available in the literature it can be developed using the ideas in [S-2] and higher 
dimensional formal patching techniques. The case where the cardinality of $G$ is divisible by $p$ is much less understood and presents more 
technical difficulties du to a phenomena of wild ramification that appears. Raynaud first attempted to understand this situation in [R]
where a rather general ``qualitative'' result is proven in the case where $G$ is a $p$-group, $X$ is a proper and smooth curve, and $f$ is 
\'etale above the generic fibre of $X$. 
Beside this the only case which is ``explicitly'' understood is the basic case where $X$ is a curve, $G$ is cyclic of order $p$, and the ring $R$ 
has inequal characteristics (cf. [G-Ma], [H], [L], [Ma], [S], [S-1], [S-3]...). Though one also has some results in inequal characteristics 
in the special case where $p$ strictly divides the cardinality of $G$ (cf. [R-1] and [We]): this case is essencially 
a by-product of the above discussed cases. One of the main objectives of this work is to go beyond 
these basic cases. In this paper we study the case where $G$ is a cyclic $p$-group and $R$ has equal characteristics. 
The above situation in equal characteristics $p$ has been inexplored so far despite the important potential applications of such a theory 
(cf. explanations below). Though this case is in principle easier than the case of inequal 
characteritics, but still the two situations are well related, as it is known to experts (this is part of the philosophy that the 
Artin-Schreier-Witt theory can be viewed as a degeneration of the Kummer theory). Also the results in equal characteristics, beside their own 
interest, can shed some lights on the situation in inequal characteristics. Although the results stated in this paper
concern only the cases of degree $p$ and $p^2$-cyclic cover we develop, pursueing the ideas and the set-up developed in [S], [S-1], [S-3], 
the techniques and the framework in which the general $p^n$-cyclic case can be studied. In what follows we review the content and results of 
this paper.

\pop {.5}
\par 
In section I we recall the classical Artin-Schreier-Witt theory, which is the main tool we use in this paper, and which provides (generically) 
explicit equations for $p^n$-cyclic covers in characteristic $p$. In section II we consider the following situation: let $X$ be a formal 
$R$-scheme of finite type which is normal connected and flat over $R$, with geometrically integral fibres, and let $f:Y\to X$ be a 
$p^n$-cyclic cover, with $Y$ normal, which is an \'etale torsor on the generic fibre of $X$. We assume that the special fibre $Y_k$ of $Y$ is reduced.
We are interested in describing the map $f$ and its special fibre $f_k:Y_k\to X_k$. One of our main results is the following:

\pop {.5}
\par
\noindent
{\bf \gr Theorem 2.2.1.}\rm \ {\sl Assume $\deg (f)=p$. Then the cover $f:Y\to X$ has the structure of a torsor under a finite and flat 
$R$-group scheme of rank $p$ over $X$}.

\pop {.5}
\par
Moreover we give an explicit description of the group schemes which appear as the group of the torsor in this situation (cf. 2.1) 
and in particular we provide ``integral'' equations  for the torsor $f$ which also provide (by reduction) equations for its special fibre. 
Next we study the case of $p^2$-cyclic covers. Our main result is theorem 2.4.3. One of the main consequences of 2.4.3 is 
that the analog of the above theorm doesn't hold: namely
one can find (generic) examples of $f$ as above of degree $p^2$ which doesn't have the structure of a torsor under a finite and flat 
$R$-group scheme of rank $p^2$ over $X$ (cf. 2.4.5). Such examples were not known in the literature before. We are however able in 2.4.3 
to find ``integral'' equations for $f$ which provide by reduction equations for
its special fibre $f_k:Y_k\to X_k$. In other terms we explain how the Artin-Schreier-Witt equations of degree $p^2$ degenerate. 
The proof of 2.4.3 is rather involved and uses the technical lemma 2.4.2. Moreover we exhibit the cases as above where $f$ has 
the structure of a torsor under a finite and flat $R$-group scheme of rank $p^2$ over $X$ (cf. 2.4.3 and 2.4.4), in which case we 
explicit the group schemes which appear as groups of the torsor. These group schemes are basically obtained by ``twisting'' 
the Artin-Schreier-Witt theory (cf. 2.3 for more details). We are also able to associate some degeneration data to the cover $f$
which determine explicitly the cover $f_k$ (cf. 2.4.8). These degeneration data will play an important role in VI were we are able
to reconstruct $p^2$-cyclic covers of formal germs of $R$-curves from the `` degeneration data''. Finally we classify in II $p$ and $p^2$-cyclic
covers above boundaries of formal germs of $R$-curves depending on their reduction type (cf. 2.5.1 and 2.6.1). 

\pop {.5}
\par
In section III we develop in equal characteristic $p$ the technique of computation of vanishing cycles as initiated in [S-1] (cf. also [Ma-Y] in the tame case) 
and which is based on the technique of compactification of covers between formal germs of $R$-curves (cf. 3.2). Our main result is 
theorem 3.2.3 which gives an explicit formula which compares the 
dimensions of the spaces of vanishing cycles in a $p$-cyclic cover between formal germs of $R$-curves and which depends on the degeneration type 
of the cover on the boundaries. This formula can easily be generalised to the case of $p^2$-cyclic covers (cf. proposition 4.1.1). We also study examples
of $p$ and $p^2$-cyclic covers above formal germs of semi-stable curves (i.e. above a smooth or a double point). Theses examples play an 
important role in the next sections V and VI. We are for example able to classify $p$ and $p^2$-cyclic covers $f:\Cal Y\to \Cal X$
between formal germs of double points wich are \'etale above the generic fibre $\Cal X_K$ of $\Cal X$ (cf. 3.3.9 and  4.2.11).

\pop {.5}
\par
In section V and VI we consider the following local situation: let $f:\Cal Y\to \Cal X$ be a cyclic cover of degree $p$ or $p^2$
above the formal germ $\Cal X\simeq \Spf R[[T]]$ of an $R$-curve at a smooth point. We study the semi-stable reduction of $f$. We are able
to exhibit in 5.2.2 and 6.2.1 ``degeneration'' data which completely determine the geometry of a semi-stable model of $\Cal Y$. These data
are of geometric and combinatorial nature and leave over the residue field $k$ of $R$ (cf. 5.2.2 and 6.2.1 for more details).
One of the main results of this paper is that we are able to 
reconstruct the cover $f$ from these data. More precisely let $\bold {Deg}$ be the set of isomorphism classes of the degeneration data
or rank $p$ (resp. rank $p^2$) as defined in 6.2.1 (resp. 5.2.2) and let $L$ be the function field of the geometric fibre of $\Cal X$. 
Then we prove using the results in III and IV the following:

\pop {.5}
\par
\noindent
{\bf \gr Theorem. 5.2.4, 6.2.3}\rm \ {\sl Let $G$ be a cyclic $p$-group of order $p$ (resp. of order $p^2$). Then there exists a canonical
specialisation map $\Sp:H^1_{\et}(\Spec L,G)\to \bold {Deg}$ which is surjective.}

\pop {.5}
\par
Although the above result is local the techniques and results we use in the proof of 5.2.4 and 6.2.3 can be used to prove a global result of lifting for
finite covers of degree $p$ or $p^2$ between semi-stable curves (such a result is proven in [M] for $p$-covers).

\pop {.5}
\par
Finally let's mention some potential applications of this work. The results of this paper can be used in order to determine the semi-stable reduction
of the Drinfeld modular curves and that of Galois covers of the projective line in characteristic $p$ with Galois groups having a 
cyclic $p$-Sylow subgroup of order $p$ or $p^2$. These results also can be used in the study of automorphism groups of curves in positive characteritics.
One of the motivations of this work was also to explore the possibility of defining reasonable Hurwitz spaces for wildly ramified Galois covers of curves.
It seems that the only case where such a definition is possible is the case of (Galois) covers of degree $p$: the right moduli problem to classify being  
torsors of rank $p$. Unfortunately our results show that when we degenerate Galois covers of degree higher than $p$ then we loose in general 
the torsor structure. this makes it really difficult to imagine a reasonable definition of what the Hurwitz spaces should be in general. 
May be the right point of view to adopt here is to restrict simply to covers where one keeps the torsor structure.

\pop {.5}
\par
Although the results of this paper concern only the cyclic groups of order $p$ and $p^2$ they can be generalized, using 
the set-up and techniques we develop, to the general $p^n$-cyclic case. I hope to come back to this question in the future. This work
was done during my visit to the Max-Planck-Institut F\"ur Mathematik in Bonn. I would like very much to thank the directors of the Institut
for their invitation and for the wonderful working atmosphere.

\pop {1}
\par
\noindent
{\bf \gr I. Artin-Schreier-Witt theory of $p^n$-cyclic covers in 
characteristic $p$}.

\rm
\pop {.5}
\par
\noindent
{\bf \gr 1.1.}\rm\ In this section we review the Artin-Schreier-Witt theory (first developed in [W])
which provides explicit equations describing cyclic covers of degeree $p^n$ in 
characteristic $p$. We refer the reader to the modern treatment of the theory as in [D-G].

\pop {.5}
\par
Let $X$ be a scheme of characteristic $p>0$ and denote by $X_{\et}$ the 
\'etale site on $X$. Let $n>0$ be an integer. We denote by $W_{n,X}$ 
(or simply $W_n$ if there is no confusion) the sheaf of Witt vectors of 
length $n$ on $X_{\et}$ (cf. [D-G], chapitre 5, 1). In the sequel any addition
or substraction of Witt vectors will mean the addition and substraction
in the sense of Witt theory. We denote by $\F$ the 
Frobenius endomorphism of $W_{n,X}$ which is locally defined by
$F.(x_1,x_2,...,x_n)=(x_1^p,x_2^p,...,x_n^p)$, for a Witt vector
$(x_1,x_2,...,x_n)$ of length $n$, and by $\Id$ the identity 
automorphism of $W_{n,X}$. The following sequence is exact on $X_{\et}$:

$$(1)\ \ 0\to (\Bbb Z/p^n\Bbb Z)_X@>i_n>> W_{n,X} @> \F-\Id>> W_{n,X} \to 0$$
where $(\Bbb Z/p^n\Bbb Z)_X$ denotes the constant sheaf 
$(\Bbb Z/p^n\Bbb Z)$ on $X_{\et}$ and $i_n$ is the natural monomorphism 
which applies $1\in\Bbb Z/p^n\Bbb Z$ to $1\in W_n$ 
(cf [G-D], chapitre 5, 5.4). From the long cohomology exact sequence 
associated to (1) one deduces the following exact sequence:

$$(2)\ \ W_{n,X}(X)@>\F-\Id>>W_{n,X}(X)\to H^1_{\et}(X,\Bbb Z/p^n\Bbb Z)\to H^1_{\et}
(X,W_n)@>\F-\Id>> H^1_{\et}(X,W_n)$$

Assume that $X=\Spec A$ is affine in which case we have $H^1_{\et}(\Spec A,W_n)=0$ and hence 
an isomorphism:
$H^1_{\et}(\Spec A,\Bbb Z/p^n\Bbb Z)\simeq W_{n,A}(A)/(\F-\Id)(W_{n,A}(A))$. This 
isomorphism has the following interpretation: to an \'etale 
$\Bbb Z/p^n\Bbb Z$-torsor $f:Y\to X=\Spec A$ above $X$ 
corresponds a Witt vector $(a_1,a_2,...,a_n)\in W_{n,A}(A)$ 
of length $n$ which is uniquely determined 
modulo addition of elements of the form $F.(b_1,b_2,...,b_n)-(b_1,b_2,...,b_n)$. 
Moreover the equations 
$F.(x_1,x_2,...,x_n)-(x_1,x_2,...,x_n)=(a_1,a_2,...,a_n)$, where the $x_i$ 
are indeterminates, are equations for the torsor $f$. More precisely there is 
a canonical factorisation of $f$ as $Y=Y_n@>f_{n}>> Y_{n-1}@>f_{n-1}>>...
@>f_2>> Y_1@>f_1>> Y_0:=X$ where each $Y_i=\Spec B_i$ is affine and $f_i:Y_i:=\Spec B_i\to f_{i-1}:
=\Spec B_{i-1}$ is the \'etale $\Bbb Z/p\Bbb Z$-torsor  corresponding to the algebra extension 
$B_{i+1}:=B_i[x_i]$. In the general case (where $H^1_{\et}(X,W_n)\neq 0$) the above equations provide local equations for an 
\'etale $\Bbb Z/p^n\Bbb Z$-torsor in characteristic $p$.

\pop {1}
\par
\noindent
{\bf \gr 1.2. Examples.}\rm\  In what follows $X$ is a scheme of characteristic $p$.

\pop {.5}
\par
\noindent
{\bf \gr 1.2.1. $\Bbb Z/p\Bbb Z$-Torsors.}\rm\ Let $f:Y\to X$ be an \'etale 
$\Bbb Z/p\Bbb Z$-torsor. Then $f$ is locally given by an equation $x^p-x=a$ 
where $a$ is a regular function on $X$ which is uniquely defined up to addition of 
elements of the form $b^p-b$ for some regular function $b$.

\pop {.5}
\par
\noindent
{\bf \gr 1.2.2. $\Bbb Z/p^2\Bbb Z$-Torsors.}\rm\ Let $f:Y\to X$ be an \'etale 
$\Bbb Z/p^2\Bbb Z$-torsor. Then we have a canonical factorisation of $f$ 
as: $Y=:Y_2@>f_2>> Y_1@>f_1>>X$ where $f_2$ and $f_1$ are \'etale 
$\Bbb Z/p\Bbb Z$-torsors. The torsor $f$ is locally given by equations
of the form:
$$F.(x_1,x_2)-(x_1,x_2):=(x_1^p-x_1,x_2^p-x_2-p^{-1}\sum _{k=1}^{p-1}{p \choose k}x_1^{pk}(-x_1)^{p-k}) =(a_1,a_2)$$ 
for some regular functions $a_1$ and $a_2$ on $X$
and the Witt vector $(a_1,a_2)$ is uniquely determined up to addition (in the Witt
theory) of vectors of the form: 
$$(b_1^p,b_2^p)-(b_1,b_2):=(b_1^p-b_1,b_2^p-b_2-p^{-1}\sum _{k=1}^{p-1}{p \choose k}b_1^{pk}(-b_1)^{p-k})$$ 

Thus locally the torsor $f_1$ is defined by the equation: 
$$x_1^p-x_1=a_1$$ 
and $f_2$ by the equation: 
$$x_2^p-x_2=a_2+p^{-1}\sum _{k=1}^{p-1}{p \choose k}x_1^{pk}(-x_1)^{p-k}$$ 

Moreover if we replace the vector $(a_1,a_2)$ by
the vector $(a_1,a_2)+(b_1^p,b_2^p)-(b_1,b_2)$ the above equations are 
replaced by:  
$$x_1^p-x_1=a_1+b_1^p-b_1$$ 
and: 
$$x_2^p-x_2=a_2+b_2^p-b_2+p^{-1}\sum _{k=1}^{p-1}{p \choose k}x_1^{pk}(-x_1)^{p-k}-
p^{-1}\sum _{k=1}^{p-1}{p \choose k}b_1^{pk}(-b_1)^{p-k}-$$
$$p^{-1}\sum _{k=1}^{p-1}{p \choose k}(b_1^p-b_1)^{k}(a_1)^{p-k}$$ 
respectively.

\pop {1}
\par
\noindent
{\bf \gr II. Degeneration of $p$ and $p^2$-cyclic covers in equal characteristic $p>0$.}
\rm
\pop {.5}
\par
In all this paragraph we use the following notations: $R$ is a complete 
discrete valuation ring of equal characteristic $p>0$ with perfect residue field $k$
and fraction field $K:=\Fr R$. We denote by $\pi$ a uniformising parameter 
of $R$.

\pop {.5}
\par
\noindent
{\bf \gr 2.1. The group schemes $\Cal M_{n}$}\rm\ (cf. also [M], 3.2). Let $n\ge 0$ be an integer and let
$\Bbb G_{a,R}=\Spec R[X]$ be the additive group scheme over $R$. The map: 
$$\phi_n:\Bbb G_{a,R}\to \Bbb G_{a,R}$$ 
given by: 
$$X\to X^p-\pi^{(p-1)n}X$$ 
is an isogeny of group schemes. The kernel of $\phi_n$ is denoted by $\Cal M_{n,R}:=\Cal M_{n}$. 
We have $\Cal M_{n}:=\Spec R[X]/(X^p-\pi^{(p-1)n}X)$ and
$\Cal M_{n}$ is a finite and flat $R$-group scheme of rank $p$. Further the following 
sequence is exact in the fppf topology:

$$(3)\ \ 0\to \Cal M_{n}\to \Bbb G_{a,R}@>\phi_n>>\Bbb G_{a,R}\to 0$$

If $n=0$ then the sequence $(3)$ is the Artin-Schreir sequence which is exact in the \'etale topology
and $\Cal M_{0}$ is the \'etale constant group scheme $(\Bbb Z/p\Bbb Z)_R$. If $n>0$ the sequence
$(3)$ has a generic fibre which is isomorphic to the \'etale Artin-Schreier sequence and a special fibre
isomorphic to the radicial exact sequence:

$$(4)\ \ 0\to \alpha_{p}\to \Bbb G_{a,k}@>x^p>>\Bbb G_{a,k}\to 0$$

Thus if $n>0$ the group scheme $\Cal M_{n}$ has a generic fibre which is \'etale isomorphic to
$(\Bbb Z/p\Bbb Z)_K$ and its special fibre is isomorphic to the infinitesimal group scheme 
$\alpha_{p,k}$. Let $X$ be an $R$-scheme. The sequence $(3)$ induces a long cohomology exact sequence:

$$(5)\ \ \Bbb G_{a,R}(X)@>\phi_n>>\Bbb G_{a,R}(X)\to H^1_{\fppf}(X,\Cal M_{n})
\to H^1_{\fppf}(X,\Bbb G_{a,R})@>\phi_n>> H^1_{\fppf}(X,\Bbb G_{a,R})$$

The cohomology group $H^1_{\fppf}(X,\Cal M_{n})$ classifies the isomorphism classes of
$\fppf$-torsors with group $\Cal M_{n}$ above $X$. The above sequence allows the 
following description of $\Cal M_{n}$-torsors:
locally a torsor $f:Y\to X$ under the group scheme $\Cal M_{n}$ is given by an equation 
$T^p-\pi^{(p-1)n}T=a$ where $T$ is an indeterminate and $a$ is a regular function on $X$ which is 
uniquely defined up to addition of elements of the form $b^p-\pi^{(p-1)n}b$ for
some regular function $b$. In particular if $H^1_{\fppf}(X,\Bbb G_{a,R})=0$ (e.g. if $X$ is affine) 
then an $\Cal M_{n}$-torsor above $X$ is globally defined by an equation as above.

\pop {.5}
\par
\noindent
{\bf \gr 2.2. Degeneration of \'etale $\Bbb Z/p\Bbb Z$-torsors.}\rm \ In what follows let $X$ be a formal
$R$-scheme of finite type which is normal connected and flat over $R$. Let $X_K:=X\times _RK$ (resp. 
$X_k:=X\times _Rk$) be the generic (resp. special) fibre of $X$. By generic fibre of $X$ we mean the 
associated $K$-rigid space (cf. [B-L]). We assume further that the special fibre $X_k$ is integral. 
Let $\eta$ be the generic point of the special fibre $X_k$ and let $\Cal O_{\eta}$ be the local ring 
at $\eta$ which is a discrete valuation ring with fraction field
$K(X)$: the function field of $X$. Let $f_K:Y_K\to X_K$ be a non trivial \'etale $\Bbb Z/p\Bbb Z$-torsor
with $Y_K$ geometrically connected and let $K(X)\to L$ be the corresponding extension of function fields. 
The main result of this section is the following:

\pop {.5}
\par
\noindent
{\bf \gr 2.2.1. Theorem.}\rm\ {\sl Assume that the ramification index above $\Cal O_{\eta}$ in the 
extension $K(X)\to L$ equals $1$. Then the torsor $f_K:Y_K\to X_K$ extends to a torsor $f:Y\to X$ 
under a finite and flat $R$-group scheme of rank $p$ with $Y$ normal. Let $\delta$ be the degree of the 
different above $\eta$ in the extension $K(X)\to L$. Then the following cases occur:
\par
{\bf a)}\ $\delta=0$ in which case $f$ is an \'etale torsor under the group scheme $\Cal M_0$ and
$f_k:Y_k\to X_k$ is then an \'etale $\Bbb Z/p\Bbb Z$-torsor.

\par
{\bf b)}\ $\delta>0$ in which case $\delta=n(p-1)$ for a certain integer $n\ge1$ and $f$ is a torsor
under the group scheme $\Cal M_{n}$. Further $f_k:Y_k\to X_k$ is in this case a radicial torsor under the $k$-group
scheme $\alpha_p$.}

\pop {.5}
\par
Note that starting from a torsor $f_K:Y_K\to X_K$ as in 2.2.1 the condition that the ramification 
index above $\Cal O_{\eta}$ equals $1$ is always satisfied after eventually a finite extension of $R$
(cf. e.g. [E]).

\pop {.5}
\par
\noindent
{\bf \gr Proof.}\rm\ We denote by $v$ the discrete valuation of $K(X)$ corresponding to the valuation ring 
$\Cal O_{\eta}$ which is normalised by $v(\pi)=1$ (note that $\pi$ is a uniformiser of $\Cal O_{\eta}$). 
We first start with the special case where $H^1_{\et}(X_K,\Bbb G_{a,K})=0$. 
In this case the torsor $f_K$ is given by an Artin-Schreier equation of the form $T^p-T=a_K$ 
where $a_K$ is a regular function on $X_K$ and we have $a_K=\pi^ma$ where 
$m\in \Bbb Z$ is an integer and $a$ is a regular function on $X$ (i.e. a function with 
$v(a)=0$). First note that
necessarily $m\le 0$ for if $m>0$ then $a_K=b^p-b$ where $b=a\pi^m+(a\pi^m)^p+(a\pi^m)^{p^2}+...
+(a\pi^m)^{p^i}+..$ (since $X_K$ is complete for the $\pi$-adic topology the above 
expression converges) and this contradicts the fact that $f_K$ is a non trivial torsor. If $m=0$ then 
the equation $T^p-T=a$ defines an \'etale $\Bbb Z/p\Bbb Z$-torsor $f:Y\to X$ above $X$ which
coincides with $f_K$ on the generic fibre and we are in the case a). 
In this case the \'etale torsor $f_k:Y_k\to X_k$ is given by the Artin-Schreier equation 
$T^p-T=\bar a$ where 
$\bar a$ is the image of $a$ modulo $\pi$. Next we treat the case where $m<0$. In this case $-m=np$ 
is necessarily
divisible by $p$ for otherwise the extension $K(X)\to L$ is totally ramified above $\Cal O_{\eta}$. 
Assume now first that the image $\bar a$ of $a$ modulo $\pi$ via the canonical map $\Gamma (X,\Cal O(X))\to 
\Gamma (X,\Cal O(X))/\pi \Gamma (X,\Cal O(X))$ is not a $p$-power. Consider the cover $f:Y\to X$ given by the equation 
$\Tilde T^p-\pi^{n(p-1)}\Tilde T=a$. Then $f$ is an $\fppf$-torsor under the group scheme $\Cal M_n$ 
which coicides with $f_K$ on the generic fibre (consider the change of variables $T:=\Tilde T/\pi^n$) 
and its special fibre  $f_k:Y_k\to X_k$
is the $\alpha_p$-torsor given by the equation $t^p=\bar a$. 
\par
In the case where $\bar a$ is a $p$-power then two cases can occur. First: $\bar a$ is a $p^s$-power for every integer $s$
which implies necessarily that $\bar a \in k$. In this case, and after some modifications which do not change the torsor $f_K$, 
we reduce to an equation of the above form and where $\bar a$ doesn't belong to $k$. To explain this assume for simplicity that $n=1$. 
Then $a={a'}^p+\pi^{\alpha}b$ where $b\in A$
and $a'\in R$. Thus the equation defining $f_K$ is $T^p-T={a'}^p/\pi ^p+\pi^{\alpha}b/\pi ^p$ which after some modification can be written as
$T^p-T={a'}/\pi +\pi^{\alpha}b/\pi ^p$ and this equation ramifies above $\pi$ which is not the case by assumption. 
This leads to the second case: there exists a positive integer $r$ such that $\bar a$ is a $p^r$ 
power but not a $p^{r+1}$ power. We assume for simplicity that $r=1$ (the general case $r>1$ is treated 
in a similar way). Let $\bar a=\bar b^p$ so that
$a=b^p+\pi \tilde b$ where $b$ and $\tilde b$ are functions on $X_K$ and $b$ is a function which reduces to 
$\bar b$ modulo $\pi$. Our equation is then of the form $T^p-T=(b/\pi^n)^p+\tilde b/\pi^{(pn-1)}$ 
and after adding
$(b/\pi^n)-(b/\pi^n)^p$ to the right hand side, which doesn't change the torsor $f_K$, we get the equation
$T^p-T=(b/\pi^n)+\tilde b/\pi^{(pn-1)}$ which can also be written in the form $T^p-T=(b/\pi^n)+b'/\pi^{n'}$ 
where $b'$ is a function with $v(b')=0$ and $n'\le pn-1$. If $n>n'$ then $n=ps$ is necessarily 
divisible by $p$ and the equation $\Tilde T^p-\pi^{s(p-1)}\Tilde T=b+\pi^{n-n'}b'$ 
defines a torsor $f:Y\to X$ under the group scheme $\Cal M_s$ which coincides with $f_K$ on 
the generic fibre and its special fibre $f_k:Y_k\to X_k$ is the $\alpha _p$-torsor given by the equation
$\Tilde t^p=\bar b$. In the case where $n'\ge n$ then $n'=s'p$ is necessarily divisible by $p$. In this case
if $\bar b'$ (resp. $\bar b'+\bar b$ in case $n'=n$) is not a p power (where 
$\bar b$ and $\bar b'$ denote the reduction of $b$ resp. $b'$ modulo $\pi$) then the equation $\Tilde T^p-\pi^{s'(p-1)}\Tilde T=
\pi^{n'-n}b+b'$ defines a torsor $f:Y\to X$ under the group scheme $\Cal M_{s'}$, 
which coincides with $f_K$ on the generic fibre, and its special fibre $f_k:Y_k\to X_k$ is the 
$\alpha _p$-torsor given by the equation $\Tilde t^p=\bar b'$
(resp. $\Tilde t^p=\bar b'+\bar b$ in case $n=n'$). Otherwise if $\bar b'$ (or $\bar b'+\bar b$ in case $n=n'$)
is a $p$-power then we repeat the same procedure as above and since 
$n$ and $n'$ decrease at each step this process must stop at some step and we end up
with an equation of the form $\Tilde T^p-\pi^{r(p-1)}\Tilde T=\tilde b$ where $\tilde b$ is a
function whose reduction modulo $\pi$ is not a $p$-power, for some positive integer $r$. Hence the required result.

\par
The argument in the general case is similar to the one used in [S] proof of 2.4. More precisely 
in general there exists an open covering $(U_i)_i$ of $X$ and 
regular functions $\tilde a_i\in \Gamma (U_{i,K},\Cal O_X)$ 
(where $U_{i,K}:=U_i\times _RK$) which are defined up to addition of functions of the form $b_i^p-b_i$ 
and such that the torsor $f_K$ is defined by the equation $T_i^p-T_i=\tilde a_i$ above $U_{i,K}$.
Now the above discussion shows that after some modifications (of the type used above) the torsor $f_K$ can be 
defined above each open $U_{i,K}$ by an equation $\Tilde T_i-\pi^{n_i(p-1)}\Tilde T=a_i$ for some uniquely 
detemined integer $n_i\ge 0$ such that in case $n_i>0$ the image $\bar a_i$ of 
$a_i$ modulo $\pi$ is not a $p$-power and such that the degree of the different $\delta_i$ above 
the generic point $\eta$ of $U_{i,k}:=U_i\times _Rk$ equals $n_i(p-1)$. From this we deduce 
that all $n_i=:n$ are equal. Then the $\Cal M_{n}$-torsor $f:Y\to X$ given locally by the equation 
$\Tilde T_i-\pi^{n(p-1)}\Tilde T_i=a_i$ above the open $U_i$ coincides on the generic fibre with 
the torsor $f_K$.

\pop {.5}
\par
\noindent
{\bf \gr 2.2.2.}\rm\ It follows from 2.2.1 that an \'etale $\Bbb Z/p\Bbb Z$-torsor above the generic fibre $X_K$ of $X$
induces canonically a {\it degeneration data} which is a torsor above the special fibre $X_k$ of $X$ under a finite and flat $k$-group scheme 
which is either \'etale or of type $\alpha_p$. Reciprocally we have the following result of {\it lifting} of such a degeneration data.

\pop {.5}
\par
\noindent
{\bf \gr 2.2.3. Proposition.}\rm\ {\sl Assume that $X$ is {\bf affine}. Let $f_k:Y_k\to X_k$ be a torsor under a finite and flat $k$-group
scheme which is \'etale (resp. of type $\alpha_p$). Then $f_k$ can be lifted to a torsor $f:Y\to X$ under a finite and flat $R$-group 
scheme of rank $p$ which is \'etale (resp. isomorphic to $\Cal M_n$ where $n>0$).}

\pop {.5}
\par
\noindent
{\bf \gr Proof.}\rm\ Since $X$ is affine the torsor $f_k$ is given by an equation $x^p-x=\bar a$ where $\bar a$ is a regular function
on $X_k$ (resp. an equation $x^p=\bar a$  where $\bar a$ is a regular function on $X_k$). Let $a$ be a regular function on $X$ which reduces to
$\bar a$ modulo $\pi$. The equation $X^p-X=a$ (resp. $X^p-\pi^{n(p-1)}X=a$ where $n>0$ is an integer) defines a cover $f:Y\to X$ above $X$ which has the 
structure of a torsor under the \'etale group scheme $(\Bbb Z/p\Bbb Z)_R$ (resp. under the group scheme $\Cal M_{n}$) and which clearly
induces the torsor $f_k$ on the special fibre $X_k$.

\pop {.5}
\par
\noindent
{\bf \gr 2.2.4. Remark.}\rm\ If $X$ is no more affine one can find examples, where $X$ is actually a proper and smooth $R$-curve, of an 
$\alpha_p$-torsor above the special fibre $X_k$ of $X$ which can not be lifted to a torsor above $X$ under a finite and flat $R$-group scheme 
of rank $p$ which is \'etale above the generic fibre of $X$. If $X_k$ is a proper and smooth curve over $k$ such a lifting is however
possible after eventually replacing $X$ by another $R$-curve $X'$ with special fibre $X_k$.

\pop {.5}
\par
\noindent
{\bf \gr 2.3. The group schemes $W_{m_1,m_2}$.}\rm \ Let $m_1$ and $m_2$ be two positive integers such that $m_2-pm_1 \ge 0$. 
We define the {\it twisted $R$-Witt group scheme} $W_{m_1,m_2}$ of length two as follows:
scheme theoretically $W_{m_1,m_2}\simeq \Bbb G_{a,R}^2$ and the group law is defined by: 
$$(x_1,x_2)+(y_1,y_2):=(x_1+y_1,x_2+y_2-p^{-1}\pi^{m_2-pm_1}\sum_{k=1}^{p-1}{p \choose k}x_1^ky_1^{p-k})$$
\par
The generic fibre  $(W_{m_1,m_2})_K$ of $W_{m_1,m_2}$ is isomorphic to the Witt group scheme $W_{2,K}$ via the map: 
$$(W_{m_1,m_2})_K\to W_{2,K}$$
$$(x_1,x_2)@>>>(x_1/\pi ^{m_1},x_2/\pi ^{m_2})$$
and its special fibre $(W_{m_1,m_2})_k$ is isomorphic either to the Witt group scheme $W_{2,k}$ if $m_2-pm_1=0$, or to the group scheme 
$\Bbb G_{a,k}^2$ otherwise.  Note finally that we have an exact sequence:
$$0\to \Bbb G_a@>V>>W_{m_1,m_2}@>R>>\Bbb G_a\to 0$$
where $V:\Bbb G_a\to W_{m_1,m_2}$ is the {\it vershiebung} homomorphism defined by $V(x)=(0,x)$ and $R:W_{m_1,m_2}\to \Bbb G_a$ is the projection
$R(x_1,x_2)=x_1$.

\pop {.5}
\par
\noindent
{\bf \gr 2.3.1. The group schemes $\Cal H_{m_1,m_2}$.}\rm \ We use the same notations as in 2.3. The following maps $I_{m_1,m_2}$ and 
$F$ are group scheme homomorphisms:
$$I_{m_1,m_2}:W_{m_1,m_2}\to W_{pm_1,pm_2}$$
$$(x_1,x_2)\to (\pi ^{m_1(p-1)}x_1,\pi ^{m_2(p-1)}x_2)$$
and: 
$$F:W_{m_1,m_2}\to W_{pm_1,pm_2}$$
$$(x_1,x_2)\to (x_1^p,x_2^p)$$
Consider the following isogeny:
$$\varphi_{m_1,m_2}:=F-I_{m_1,m_2}:W_{m_1,m_2}\to W_{pm_1,pm_2}$$
which is given by: 
$$(x_1,x_2)\to (x_1^p-\pi ^{m_1(p-1)}x_1, x_2^p-\pi ^{m_2(p-1)}x_2-p^{-1}\sum_{k=1}^{p-1}{p \choose k}\pi ^{m_2p-m_1(pk+p-k)}x_1^{pk}
(-x_1)^{p-k})$$
\par
We define the group scheme  $\Cal H_{m_1,m_2}$ to be the kernel of the above isogeny. Thus we have an exact sequence:
$$(6)\ \ 0\to \Cal H_{m_1,m_2}\to W_{m_1,m_2}@>F-I_{m_1,m_2}>>W_{pm_1,pm_2}\to 0$$
and $\Cal H_{m_1,m_2}$ is a finite and flat commutative $R$-group scheme of rank $p^2$. Further we have the following commutative diagram:

$$
\CD
@. 0 @. 0 @. 0 \\
@.   @VVV     @VVV    @VVV \\
0 @>>> \Cal M_{m_2}  @>>>  \Bbb G_a @>\phi_{m_2}>>               \Bbb G_a @>>> 0 \\
@.           @VVV                @VVV                                  @VVV    @.        \\
0 @>>> \Cal H_{m_1,m_2} @>>>  W_{m_1,m_2}   @>\varphi _{m_1,m_2}>>   W_{pm_1,pm_2} @>>> 0 \\
@.           @VVV                @VVRV                                  @VVRV    @.        \\
0 @>>>    \Cal M_{m_1}  @>>>       \Bbb G_a      @>\phi_{m_1}>>                      \Bbb G_a @>>> 0 \\
@.           @VVV               @VVV                                               @VVV              @.\\
@.              0       @.       0            @.                                      0 \\
\endCD
$$

The group scheme $\Cal H_{m_1,m_2}$ is then an extension of the group scheme $\Cal M_{m_1}$ by $\Cal M_{m_2}$.
Its generic fibre $(\Cal H_{m_1,m_2})_K$ is isomorphic to the \'etale constant group $\Bbb Z/p^2\Bbb Z$ and its special fibre
$(\Cal H_{m_1,m_2})_k$ is either the group scheme $H_k$ which is an extension of $\Bbb Z/p\Bbb Z$ by $\alpha_p$ if $m_1=0$ and $m_2>0$ (such an 
extension is necessarily split), or the group scheme $G_k$ which is an extension of $\alpha_p$ by  $\alpha_p$ if $m_1>0$. 
We have the following exact sequences:
$$(7)\ 0\to H_k\to \Bbb G_{a,k}^2 @>(x_1^p-x_1,x_2^p)>> \Bbb G_{a,k}^2\to 0$$
and:
$$(8)\ 0\to G_k\to \Bbb G_{a,k}^2@>(x_1^p,x_2^p)>> \Bbb G_{a,k}^2\to 0$$

Let $X$ be an $R$-scheme. The sequence $(6)$ induces a long cohomology exact sequence (9):
$W_{m_1,m_2} (X)@>\varphi _{m_1,m_2}>> W_{pm_1,pm_2} (X)   \to H^1_{\fppf}(X,\Cal H_{m_1,m_2})
\to H^1_{\fppf}(X,W_{m_1,m_2})@>\phi_n>> H^1_{\fppf}(X,W_{pm_1,pm_2})$. 
\par
The cohomology group $H^1_{\fppf}(X,\Cal H_{m_1,m_2})$ 
classifies the isomorphism classes of
$\fppf$-torsors with group $\Cal H_{m_1,m_2}$ above $X$. The above exact sequence allows the 
following description of $\Cal H_{m_1,m_2}$-torsors:
locally a torsor $f:Y\to X$ under the group scheme $\Cal H_{m_1,m_2}$ is given by the equations: 
$$T_1^p-\pi^{m_1(p-1)}T_1=a_1$$ 
and:  
$$T_2^p-\pi^{m_2(p-1)}T_2=a_2+p^{-1}\sum_{k=1}^{p-1}{p \choose k}\pi ^{pm_2-m_1(pk+p-k)}T_1^{pk}(-T_1)^{p-k}$$ 
where $T_1$ and $T_2$ are indeterminates and $a_1$, $a_2$ are regular functions on $X$. Its special fibre is either the $H_k$-torsor
given by the equations:
$$t_1^p-t_1=\bar a_1$$ 
and:  
$$t_2^p=\bar a_2$$
if $m_1=0$, or the $G_k$-torsor given by the equations: 
$$t_1^p=\bar a_1$$ 
and:  
$$t_2^p=\bar a_2$$
otherwise. In particular if $H^1_{\fppf}(X,\Bbb G_{a,R})=0$ (e.g. if $X$ is affine) then an $\Cal H_{m_1,m_2}$-torsor above $X$ is globally 
defined by an equation as above.

\pop {.5}
\par
\noindent
{\bf \gr 2.4. Degeneration of \'etale $\Bbb Z/p^2\Bbb Z$-torsors.}\rm \ Our aim in this section
is to describe explicitely the degeneration of \'etale $\Bbb Z/p^2\Bbb Z$-torsors. In what follows we use 
the same notations as in  2.2. Let $f_K:Y_K\to X_K$ be a non trivial $\Bbb Z/p^2\Bbb Z$-torsor. Let $K(X)\to L$ be the cyclic 
$p^2$-extension of function fields corresponding to the torsor $f_K$ which canonically factorises as $K(X)\to L_1\to L_2:=L$ where 
$K(X)\to L_1$ is a cyclic $p$-extension. We assume that the ramification index above 
the generic point $\eta$ of $X_k$ in the extension $K(X)\to L$ equals $1$. There 
exists a canonical factorisation $Y_K=Y_{2,K}@>f_{2,K}>>Y_{1,K}@>f_{1,K}>>X_K$ of $f_K$ where $f_{i,K}$ is a  
$\Bbb Z/p\Bbb Z$-torsor for $i\in \{1,2\}$. Moreover by 2.2.1 the torsor $f_{2,K}$ (resp. $f_{1,K}$) extends 
to a torsor $f_2:Y_2\to Y_1$ (resp. $f_1:Y_1\to X$) under a finite and flat $R$-group scheme of rank $p$ and the composite 
$f:=f_1\circ f_2$ is a finite and flat cover which coincides on the generic fibre with $f_K$. We assume that the special fibre $Y_{2,k}$ of 
$Y_2$ is {\bf irreducible}. In particular above the generic point $\eta$ there exists a unique generic point $\eta_1$ in $Y_1$ which lies above $\eta$. 
We denote by $\delta$ (resp. $\delta_1$ and $\delta_2$) the degree of different in the extension $L$ above the point $\eta$ (resp. the degree of 
different in the extension $L_1$ above the point $\eta$ and that of the different in the extension 
$L_2$ above the point $\eta_1$). Note that $\delta=\delta_1+\delta_2$.

\pop {.5}
\par
\noindent
{\bf \gr 2.4.1.}\rm\ We start with the following technical lemma 2.4.2 which will be used later on. In what follows we assume 
that $X=\Spec A$ is affine and that $f_1:Y_1:=\Spec B\to X$ is a torsor under the group scheme $\Cal M_{n}$ for some integer $n>0$. Thus
$f_1$ is given by an equation $T^p-\pi^{n(p-1)}T=v$ where $v\in A$ is such that its 
image $\bar v\in \overline A:=A/\pi A$ is not a $p$-power. In particular the special fibre  
$\bar f_1:Y_{1,k}=\Spec \overline B\to X_k=\Spec \overline A$ of the torsor
$f$ (here $\overline B:=B/\pi B$ and $\overline A:=A/\pi A$) is the $\alpha_p$-torsor given by the equation $t^p=\bar v$ 
and $\overline B$ is a free $\overline A$-algebra with basis $\{1,t,t^2,...,t^{p-1}\}$. We need to characterize 
elements of $A$ which become $p$-powers modulo $\pi$ in $\overline B$ but are not necessarily 
$p$-powers modulo $\pi$ in $\overline A$.

\pop {.5}
\par
\noindent
{\bf \gr 2.4.2. Lemma.}\rm\ {\sl Let $u\in A$. Assume that the image $\bar u$ of $u$
is a $p$-power in $\overline B$. Then $u=f(v)+\pi u'$
where $u'\in A$ and $f(v)$ belongs to the additive subgroup $A_v:=A^p\oplus A^p.v\oplus....\oplus A^p.v^{p-1}$ of $A$. 
Moreover let $f(v):=a_0^p+a_1^pv+...+a_{p-1}^pv^{p-1}\in A_v$ and let $m$ be an integer.
Consider the element $g:=f(v)\pi^{-pm}\in A_K$. Then $g=(a_o^p+a_1^p(T^p-\pi^{n(p-1)}T)+
...+a_{p-1}^p(T^p-\pi^{n(p-1)}T)^{p-1}\pi^{-pm}$ in $B_K$ and after addition of elements of $B_K$ of the form
$b^p-b$ one can transform $g$ in $\tilde g:=\pi ^{-m} (a_0+a_1T+...+a_{p-1}T^{p-1})+\pi^{-(pm-n(p-1))}
(-\sum _{j=1}^{p-1} ja_j^pT^{p(j-1)}T)+\pi^{-(pm-2n(p-1))}h(T)$ where 
$h(T)\in B$. Moreover the image of $-\sum _{j=1}^{p-1} ja_j^pT^{p(j-1)}T=
-a_1^pT-2a_2^pT^{p+1}-...-(p-1)a_{p-1}^pT^{p(p-2)+1}$ in $\overline B$ is not a $p$-power.

\pop {.5}
\par
\noindent
{\bf \gr Proof.}\rm\ We have $\overline B=\overline A\oplus \overline A.t\oplus...\oplus \overline A.t^{p-1}$,
hence $\overline B^p=\overline A^p\oplus \overline A^p.\bar v\oplus...\oplus\overline A.\bar v^{p-1}$
and the first assertion of the lemma follows. Now let $g:=(a_0^p+a_1^pv+...+a_{p-1}^pv^{p-1})\pi^{-pm}\in B_K$
then since $T^p-\pi^{n(p-1)}T=v$ in $B_K$ we can write $g=\pi^{-pm}(a_o^p+a_1^p(T^p-\pi^{n(p-1)}T)+
...+a_{p-1}^p(T^p-\pi^{n(p-1)}T)^{p-1}$ in $B_K$. After developing the terms $(T^p-\pi^{n(p-1)}T)^j$
for $j\in \{1,p-1\}$ according to the binomial expansion and puting together the terms with the same
power of $\pi$ we get that $g=\pi ^{-pm} (a_0^p+a_1^pT^p+...+a_{p-1}^pT^{p(p-1)})
+\pi^{-(pm-n(p-1))}(-a_1^pT-2a_2^pT^{p+1}-...-(p-1)a_{p-1}^pT^{p(p-2)+1})+\pi^{-(pm-2n(p-1))}h(T)$
where $h(T)\in B$. Finally after adding $(a_0+a_1T+...+a_{p-1}T^{p-1})/\pi ^{m}-(a_0^p+a_1^pT^p+...
+a_{p-1}^pT^{p(p-1)})/\pi ^{pm}$ we get the desired expression for $\tilde g$.

\pop {.5}
\par
The next theorem is the main result of this section. It describes locally and explicitely the degeneration of \'etale $\Bbb Z/p^2\Bbb Z$-torsors.
Although we loose the structure of torsors in this case of cyclic covers of degree $p^2$ (cf. Proposition 2.4.5) we are able to find 
``canonical integral equations'' which describe the reduction of $p^2$-cyclic covers in equal characteristic $p$.

\pop {.5}
\par
\noindent
{\bf \gr 2.4.3. Theorem.}\rm\ {\sl We use the same notations as in 2.4. Assume that $X=\Spec A$ is {\bf affine}. Then the torsor 
$f_K$ can be described by an equation of the form: 
$$(T_1^p,T_2^p)-(T_1,T_2)=(\pi ^{m_1}a_1,\pi ^{m_2}a_2)$$
where $a_1$ and $a_2$ are regular functions on $X_K$ with $v(a_1)=v(a_2)=0$, $m_1\le 0$, $m_2$ is an integer, with the following cases which occur:
\pop {.5}
\par
{\bf a)}\ $m_1=0$ and $m_2\ge 0$. In this case $f$ is an \'etale $\Bbb Z/p^2\Bbb Z$-torsor above $X$ given by the equation: 
$$(T_1^p,T_2^p)-(T_1,T_2)=(a_1,\pi ^{m_2}a_2)$$
Its special fibre $f_k:Y_k\to X_k$ is the \'etale $\Bbb Z/p^2\Bbb Z$-torsor given by 
the equation: 
$$(t_1^p,t_2^p)-(t_1,t_2)=(\bar a_1,\overline {\pi ^{m_2}a_2)}$$ 
and $\delta=\delta_1=\delta_2=0$.

\pop {.5}
\par
{\bf b)}\ $m_1=0$, $m_2=-pm_2'<0$, and $a_2$ is not a $p$-power modulo $\pi$. In this case $f$ is a torsor under the $R$-group scheme 
$W_{0,m_2'}$ given by the equations: 
$$T_1^p-T_1=a_1$$ 
and 
$$\Tilde T_2^p-\pi ^{m_2'(p-1)}\Tilde T_2=a_2+p^{-1}\pi ^{-m_2}\sum _{k=1}^{p-1}{p \choose k}T_1^{pk}(-T_1)^{p-k}$$ 
Its special fibre is the torsor under the $k$-group scheme $(W_{0,m_2'})_k\simeq H_k$ given by the equations: $t_1^p-t_1=\bar a_1$,
and $\tilde t_2^p=\bar a_2$, where $\bar a_1$ (resp. $\bar a_2$) is the image of $a_1$ (resp. of $a_2$) modulo $\pi$.
In this case $\delta_1=0$ and $\delta=\delta_2=m_2'(p-1)$.

\pop {.5}
\par
{\bf c)}\ $m_1=-pm_1'<0$ and the image $\bar a_1$ of $a_1$ modulo $\pi$ is not a $p$-power. 
In this case $f_1$ is an $\Cal M_{m_1'}$-torsor given by the equation: 
$$\Tilde T_1^p-\pi ^{m_1'(p-1)}\Tilde T_1=a_1$$ 
and its special fibre $f_{1,k}:Y_{1,k}\to X_{k}$ is the 
$\alpha_p$-torsor given by the equation $\tilde t_1^p=\bar a_1$. 
We have $\delta_1=m_1'(p-1)$. For $f_2$ one has the following cases:

\pop {.5}
\par
{\bf c-1)}\ $m_1'(p(p-1)+1)>-m_2$ (resp. $m_1'(p(p-1)+1)=-m_2$)
In this case $m_1'=pm_1''$ is necessarily divisible by $p$ and if $\tilde m_1:=m_1''(p(p-1)+1)$ then $f_2$ is a torsor under
$\Cal M_{\tilde m_1,R}$ given by the equation:
$$\Tilde T_2^p-\pi ^{\tilde m_1(p-1)}\Tilde T_2=\pi^{\tilde m_1p+m_2}a_2+p^{-1}\sum _{k=1}^{p-1}{p \choose k}
\Tilde T_1^{pk}(-\Tilde T_1)^{p-k}\pi ^{m_1'((p(p-1)+1)-pk-p+k)}$$
Its special fibre is the $\alpha_p$-torsor given by the equation: $\tilde t_2^p
=-\tilde t_1^{p(p-1)+1}$ (resp. the equation: $\tilde t_2^p=-\tilde t_1^{p(p-1)+1}+\bar a_2$).

\pop {.5}
\par
Otherwise $-m_2>m_1'(p(p-1)+1)$ in which case $-m_2=pm_2'$ is necessarily divisible by $p$ and we have the following description 
for $\pi^{m_2}a_2$:
$$\pi^{m_2}a_2=f_1(a_1)/\pi^{pm_2'}+f_2(a_1)/\pi^{pm_2'-t_1}+...+f_r(a_1)
/\pi^{pm_2'-t_1-...-t_{r-1}}+g/\pi^{pm_2'-t_1-...-t_r}$$ where $f_i(a_1)$ belongs to the subgroup $A_{a_1}$ of $A$ (cf. 2.4.2), 
$g\in A$, and the $t_i$ are positive integers (note that $g$ and the $f_i$ can be $0$). Moreover the torsor $f_{2,K}$ 
is given by the equation: $T_2^p-T_2=f_1(a_1)/\pi^{pm_2'}+f_2(a_1)/\pi^{pm_2'-t_1}+...+f_r(a_1)/\pi^{pm_2'-t_1-...-t_{r-1}}$
\newline
$+g/\pi^{pm_2'-t_1-...-t_r}+p^{-1}\sum _{k=1}^{p-1}{p \choose k}\Tilde T_1^{pk}(-\Tilde T_1)^{p-k}\pi^{-m_1'(pk+p-k)}.$ 
\par
And the following distinct cases occur:

\pop {.5}
\par
{\bf c-2)}\ $pm_2'-(p-1)m_1'>\sup (m_1'(p(p-1)+1),pm_2'-t_1-...-t_r)$ (resp. $pm_2'-(p-1)m_1'=m_1'(p(p-1)+1)>pm_2'-t_1-...-t_r$). 
In this case $pm_2'-m_1'(p-1)=m_2''p$ is divisible by $p$ and $\delta_2=m_2''(p-1)$. 
Let $f_1(a_1):=c_0^p+c_1^pa_1+...+c_{p-1}^pa_1^{p-1}$. Then $f_2$ is a torsor under 
$\Cal M_{m_2''}$ and its special fibre is the $\alpha_p$-torsor given by the equation: $\tilde t_2^p=
-\bar c_1^p\tilde t_1-2\bar c_2^p\tilde t_1^{p+1}-...-(p-1)\bar c_{p-1}^p\tilde t_1^{p(p-2)+1}$ (resp. the equation: 
$\tilde t_2^p=-\bar c_1^p\tilde t_1-2\bar c_2^p\tilde t_1^{p+1}-...-(p-1)\bar c_{p-1}^p\tilde t_1^{p(p-2)+1}-\tilde t_1^{p(p-1)+1}$)
where $\bar c_i$ is the image of $c_i$ modulo $\pi$.

\pop {.5}
\par
{\bf c-3})\ $pm_2'-t_1-...-t_r>\sup (pm_2'-(p-1)m_1',m_1'(p(p-1)+1))$ (resp. $pm_2'-(p-1)m_1'=pm_2'-t_1-...-t_r>m_1'(p(p-1)+1)$) 
and the image of $g$ modulo $\pi$ is not a $p$-power in
$\Cal O(Y_{1,k})$. In this case $pm_2'-t_1-...-t_r=pm_2''$ is divisible by $p$, $\delta_2=m_2''(p-1)$, and 
$f_2$ is an $\Cal M_{m_2''}$-torsor. Its special fibre is the $\alpha_p$-torsor given by the equation: $\tilde t_2^p=\bar g$ (resp. the equation:
$\tilde t_2^p=-\bar c_1^p\tilde t_1-2\bar c_2^p\tilde t_1^{p+1}-...-(p-1)\bar c_{p-1}^p\tilde t_1^{p(p-2)+1}+\bar g$).

\pop {.5}
\par
{\bf c-4)}\ $m_1'(p(p-1)+1)>\sup (pm_2'-t_1-...-t_r,pm_2'-(p-1)m_1')$ (resp. $pm_2'-t_1-...-t_r=m_1'(p(p-1)+1)>pm_2'-(p-1)m_1'$)
In this case $m_1'=pm_1''$ and if $\tilde m_1:=m_1''(p(p-1)+1)$ then $f_2$ is a torsor under
$\Cal M_{\tilde m_1,R}$. Its special fibre is the $\alpha_p$-torsor given by the equation: $\tilde t_2^p
=-\tilde t_1^{p(p-1)+1}$ (resp. the equation: $\tilde t_2^p=-\tilde t_1^{p(p-1)+1}+\bar g$). 

\pop {.5}
\par
{\bf c-5)}\ $pm_2'-t_1-...-t_r=m_1'(p(p-1)+1)=pm_2'-(p-1)m_1'$. In this case $m_1'=pm_1''$ and if $\tilde m_1:=m_1''(p(p-1)+1)$ 
then $f_2$ is a torsor under $\Cal M_{\tilde m_1,R}$. Its special fibre is the $\alpha_p$-torsor given by the equation: $\tilde t_2^p
=-\bar c_1^p\tilde t_1-2\bar c_2^p\tilde t_1^{p+1}-...-(p-1)\bar c_{p-1}^p\tilde t_1^{p(p-2)+1}-\tilde t_1^{p(p-1)+1}+\bar g$. 

\pop {.5}
\par
In all the above cases if $f_1$ (resp. $f_2$) is a torsor under the group scheme $\Cal M_{\tilde m_1}$ (resp. $\Cal M_{\tilde m_2}$) 
then necessarily $\tilde m_2\ge \tilde m_1(p(p-1)+1)/p$.
Moreover in all the cases c-2, c-3, c-4, and c-5 above the functions $\bar c_1$, $\bar c_2$,..., $\bar c_{p-1}$ (resp. $\bar g$) are uniquely
determined (resp. is uniquely determined up to addition of $\bar h^p$ where $\bar h$ is a regular function on $X_k$). In the case c-1 the function
$\bar a_2$ is uniquely determined up to addition of $\bar b^p$ where $\bar b$ is a regular function on $X_k$.}

\pop {.5}
\par
\noindent
{\bf \gr Proof.}\rm\ The torsor $f_K$ is given by the Artin-Schreier-Witt theory by an equation
of the form $(T_1^p,T_2^p)-(T_1,T_2)=(\tilde a_1,\tilde a_2)$ where $\tilde a_1$ and $\tilde a_2$ are regular 
functions on $X_K$. We can write $\tilde a_1=\pi^{m_1}a_1$ (resp. $\tilde a_2=\pi^{m_2}a_2$) 
with $v(a_1)=v(a_2)=0$ and it follows from 2.2.1 that we have necessarily $m_1\le 0$. If $m_1=0$ and $m_2\ge 0$ we are
in case a) and the assertion there is clear. Assume that $m_1=0$ and $m_2<0$. Then it follows from 2.2.1 that $m_2=-pm_2'$ is 
necessarily divisible by $p$ and after eventually some modifications (as in the proof of 2.2.1) we may assume 
that $a_2$ is not a $p$-power modulo $\pi$ (here one uses the fact that a regular 
function $u$ on $X$ which is not a $p$-power modulo $\pi$ in $X_k$ can not become a $p$-power in $Y_{1,k}$ since $f_{1,k}$ is an 
\'etale torsor hence not radicial). Then $f$ is defined by the equations:
$$T_1^p-T_1=a_1$$ 
and 
$$\Tilde T_2^p-\pi ^{m_2'(p-1)}\Tilde T_2=a_2+p^{-1}\pi ^{-m_2}\sum _{k=1}^{p-1}{p \choose k}T_1^{pk}(-T_1)^{p-k}$$ 
where $\Tilde T_2:=\pi ^{m_2'}T_2$. The rest of the assertion in case b) follows then easily.
Assume now that $m_1<0$. Then the assertion concerning $f_1$ follows from 2.2.1. 
Assume first that $m_2\ge 0$. The assertion concerning $f_2$ follows then easily after adapting the equation defining the torsor $f_{2,K}$
to the change of variables $T_1=\Tilde T_1/\pi ^{m_1'}$ and we are in the case c-1). 
In this case  $m_1'=pm_1''$ is divisible by $p$ and the cover $f$ is given by the equations:
$$\Tilde T_1^p-\pi ^{m_1'(p-1)}\Tilde T_1=a_1$$
and 
$$\Tilde T_2^p-\pi ^{\tilde m_1(p-1)}\Tilde T_2=\pi^{\tilde m_1p+m_2}a_2+
p^{-1}\sum _{k=1}^{p-1}{p \choose k}\tilde T_1^{pk}(-\tilde T_1)^{p-k}\pi^{\tilde m_1p-m_1'(pk+p-k)}$$
where $\tilde m_1:=m_1''(p(p-1)+1)$. From which we deduce that its special fibre is given by the equations: $\tilde t_1^p=\bar a_1$ and 
$\tilde t_2^p=-\tilde t_1^{p(p-1)+1}$ where $\tilde t_2$ (resp. $\tilde t_1$) is the image of $\tilde T_2$ (resp. image of $\tilde T_1$)
modulo $\pi$.
\par
So lastly we assume that $m_2<0$. Then $f_{2,K}$ is given by the equation:
$$T_2^p-T_2=\pi^{m_2}a_2+p^{-1}\sum _{k=1}^{p-1}{p \choose k}\tilde T_1^{pk}(-\tilde T_1)^{p-k}\pi^{-m_1'(pk+p-k)}$$
The highest power of $\pi$ in the denominators of the summands in: 
$$p^{-1}\sum _{k=1}^{p-1}{p \choose k}\tilde T_1^{pk}(-\tilde T_1)^{p-k}\pi^{-m_1'(pk+p-k)}$$ 
is $m_1'(p(p-1)+1)$ and in order to understand the reduction of the torsor $f_2$ we have to compare this to $m_2$. Assume first that
$m_1'(p(p-1)+1)> -m_2$. Then it follows from 2.2.1 that $m_1'=m_1''p$ must be divisible by $p$.
Let $\tilde m_1:=m_1''(p(p-1)+1)$. Then we are in case c-1) and $f_2$ is a torsor
under the group scheme $\Cal M_{\tilde m_1,R}$ defined by the equation:
$$\Tilde T_2^p-\pi^{\tilde m_1(p-1)}\Tilde T_2=\pi^{\tilde m_1p+m_2}a_2+p^{-1}\sum _{k=1}^{p-1}
{p \choose k}\pi ^{\tilde m_1p-m_1'(pk+p-k)}\tilde T_1^{pk}(-\tilde T_1)^{p-k}$$ 
Its special fibre is the $\alpha_p$-torsor given by the equation:
$\tilde t_2^p=-\tilde t_1^{p(p-1)+1}$. Assume next that $m_1'(p(p-1)+1)<-m_2$ (the case where $m_1'(p(p-1)+1)=-m_2$
is easily treated and is left to the reader). Then it follows from 2.2.1 that $m_2=-pm_2'$ is divisible by $p$
and two cases occur depending on whether or not the image $\bar a_2$ of $a_2$ modulo $\pi$ is or is not a $p$ power in $\Cal O(Y_{1,k})$.
If $\bar a_2$ is not a $p$ power in $\Cal O(Y_{1,k})$ then $f_2$ is a torsor under $\Cal M_{m_2',R}$ given by the equation:
$$\Tilde T_2^p-\pi^{m_2'(p-1)}\Tilde T_2=a_2+p^{-1}\sum _{k=1}^{p-1}{p \choose k}\pi ^{pm_2'-m_1(pk+p-k)}\tilde T_1^{pk}(-\tilde T_1)^{p-k}$$
Its special fibre is the $\alpha_p$-torsor given by the equation $\tilde t_2^p=\bar a_2$ and we are in the case c-3).
Assume that $\bar a_2$ is a $p$ power in $\Cal O(Y_{1,k})$. Then either $a_2$ is already a $p$-power in 
$\Cal O(X_{k})$ in which case we can transform, using the kind of transformations used in the proof of 2.2.1,
the term $\pi ^{m_2}a_2$ into $\pi ^{\tilde m_2}\tilde a_2$ where $0>\tilde m_2>m_2$ and $\tilde a_2\in A$.
Or $a_2$ is not a $p$-power in $\Cal O(X_{1})$ but becomes a $p$-power in $\Cal O(Y_{1,k})$. In
the later case it follows from 2.4.2 that $a_2=f_1(a_1)+\pi^{t_1}g_1$ where $f_1(a_1):=c_0^p+c_1^pa_1+...+ 
c_{p-1}^pa_1^{p-1}$ belongs to the subgroup $A_{a_1}$ of $A$, $t_1>0$, and $g_1\in A$. Moreover the term  
$\pi ^{m_2}a_2=f_1(a_1)/\pi^{pm_2'}+g_1/\pi^{pm_2'-t_1}$ can be transformed to 
$\tilde f_1(T_1)/\pi^{pm_2'-m_1'(p-1)}+g_1/\pi^{pm_2'-t_1}$ where the image $\overline {\tilde f_1(T_1)}:=
-\bar c_1^p\tilde t_1-2\bar c_2^p\tilde t_1^{p+1}-...-(p-1)\bar c_{p-1}^p\tilde t_1^{p(p-2)+1}$ of $\tilde f_1(T_1)$ 
modulo $\pi$ is not a $p$-power (cf. loc. cit.). At this point we can repeat the same argument as above. Namely 
if in the first case the image $\bar {\tilde a}_2$ of $\tilde a$ modulo $\pi$ is not a $p$-power in 
$\Cal O(Y_{1,k})$ then we conclude as above that we are either in case c-3) if $\tilde m_2\ge m_1'(p(p-1)+1)$
in which case $\tilde m_2=pm_2''$ and $f_2$ is a torsor under $\Cal M_{m_2'',R}$ whose special fibre is 
the $\alpha_p$-torosr given by the equation $\tilde t_2^p=\bar {\tilde a}_2$. Otherwise we repeat the same 
process as above. And in the second case if $pm_2'-(p-1)m_1'>\sup (pm_2'-t_1,m_1'(p(p-1)+1))$ then $pm_2'-(p-1)m_1'=:pm_2''$
is divisible by $p$ and $f_2$ is a torsor under the group scheme  $\Cal M_{m_2'',R}$ defined by the equation:
$$\Tilde T_2^p-\pi^{m_2''(p-1)}\Tilde T_2=\tilde f_1(T_1)+\pi^{pm_2''-pm_2'+t_1}g_1
+p^{-1}\sum _{k=1}^{p-1}{p \choose k}\pi ^{pm_2''-(pk+p-k)m_1}\tilde T_1^{pk}(-\tilde T_1)^{p-k}$$ 
Its special fibre is the $\alpha_p$-torsor given by the equation: $\tilde t_2^p=
-\bar c_1^p\tilde t_1-2\bar c_2^p\tilde t_1^{p+1}-...-(p-1)\bar c_{p-1}^p\tilde t_1^{p(p-2)+1}$ and we are in the case c-2).
If $m_1'(p(p-1)+1)>\sup (pm_2'-t_1,pm_2'-(p-1)m_1')$ then $m_1'(p(p-1)+1)$ is divisible by $p$, 
$f_2$ is a torsor under the group scheme  $\Cal M_{(m_1'(p(p-1)+1))/p,R}$, its special fibre is the $\alpha_p$-torsor 
given by the equation: 
\newline
$\tilde t_2^p=-\tilde t_1^{p(p-1)+1}$ and we are in the case c-4). If
$pm_2'-t_1>\sup (m_1'(p(p-1)+1),pm_2'-(p-1)m_1')$ and the image $\bar g_1$ of $g_1$ in $\Cal O(Y_{1,k})$
is not a $p$-power then $pm_2'-t_1=:pm_2''$ is divisible by $p$, $f_2$ is a torsor under the group scheme  
$\Cal M_{m_2'',R}$ defined by the equation: 
$$\Tilde T_2^p-\pi^{m_2''(p-1)}\Tilde T_2=\pi^{\tilde m_2}
\tilde f_1(T_1)+g_1+p^{-1}\sum _{k=1}^{p-1}{p \choose k}\pi ^{pm_2''-(pk+p-k)m_1}\tilde T_1^{pk}(-\tilde T_1)^{p-k}$$ 
where $\tilde m_2:=pm_2''-pm_2'-(p-1)m_1'$, and its special fibre is the $\alpha_p$-torsor given by the equation: 
$\tilde t_2^p=\bar g_1$ and we are in the case c-3). Now in the general case we repeat the same argument as above
if in the first case the image $\bar {\tilde a}_2$ of $\tilde a_2$ modulo $\pi$ is a $p$-power or if in the second 
case $pm_2'-t_1>\sup (m_1'(p(p-1)+1),pm_2'-(p-1)m_1')$ and the image $\bar g_1$ of $g_1$ in $\Cal O(Y_{1,k})$
is a $p$-power. As the exponent of $\pi$ in the denominators in the equation defining $f_{2,K}$ decreases at each step
we conclude that this process must stop after finitely many steps and we end up with an equation as claimed in
the statement c) and the rest of the conclusion follows then easily.

\pop {.5}
\par
\noindent
{\bf \gr 2.4.4. Remark.}\rm\ Assume that we are in the case c-3) of 2.4.3 and that $t_1=...=t_r=0$ and $f_1=...=f_r=0$. 
Then $f$ is a torsor under the $R$-group scheme $\Cal H_{m_1',m_2'}$ given by the equations: 
$$\Tilde T_1^p-\pi ^{m_1'(p-1)}\Tilde T_1=a_1$$ 
and: 
$$\Tilde T_2^p-\pi ^{m_2'(p-1)}\Tilde T_2=g
+p^{-1}\sum _{k=1}^{p-1}{p \choose k}\Tilde T_1^{pk}(-\Tilde T_1)^{p-k}\pi ^{pm_2'-pk-p+k}$$ 
Its special fibre is the $(\Cal H_{m_1',m_2'})_k\simeq H_k$-torsor given by the equations: $\tilde t_1^p=\bar a_1$ and $\tilde t_2^p=\bar g$.

\pop {.5}
\par
The next proposition states that in general the cover $f_2$ in 2.4.3 doesn't have the structure of a torsor under a finite and flat 
$R$-group scheme of rank $p^2$.

\pop {.5}
\par
\noindent
{\bf \gr 2.4.5. Proposition.}\rm\ {\sl We use the same notations as in 2.4.3. Assume that we are in the case c-1) and that $p>2$. Then the finite cover $f$
doesn't have the structure of a torsor under a finite and flat $R$-group scheme of rank $p^2$.}

\pop {.5}
\par
\noindent
{\bf \gr Proof.}\rm\ Recall that in the case c-1 of 2.4.3 the \'etale torsor $f_K$ is given by the equations (1):  
$(T_1^p,T_2^p)-(T_1,T_2)=(a_1/\pi^{pm_1'},a_2':=\pi ^{m_2}a_2)$. Let $b_1\in A$. Then the equations:
$(S_1^p,S_2^p)-(S_1,S_2)=(\pi^{-pm_1'}a_1,a_2')+(\pi^{-pm_1'}b_1^p,0)-(\pi^{-m_1'}b_1,0)$ are also defining equations for the torsor $f_K$. 
The cover $f$ is then given by the equations (1):
$$\Tilde T_1^p-\pi ^{m_1'(p-1)}\Tilde T_1=a_1$$ 
and  
$$\Tilde T_2^p-\pi ^{\tilde m_1(p-1)}\Tilde T_2=
\pi^{\tilde m_1p+m_2}a_2+p^{-1}\sum _{k=1}^{p-1}{p \choose k}\Tilde T_1^{pk}(-\Tilde T_1)^{p-k}\pi ^{\tilde m_1p-m_1'(pk+p-k)}$$ 

And $f$ is also given by the equations (2): 
$$\Tilde S_1^p-\pi ^{m'(p-1)}\Tilde S_1=a_1+b_1^p-\pi^{m(p-1)}b_1$$ 
and:\ \  $S_2^p-S_2=a_2'+p^{-1}(\sum _{k=1}^{p-1}{p \choose k}\Tilde S_1^{pk}(-\Tilde S_1)^{p-k}\pi ^{-m_1'(pk+p-k)}$
\newline
$-\sum _{k=1}^{p-1}{p \choose k}b_1^{pk}(-b_1)^{p-k}\pi ^{-m_1'(pk+p-k)}-\sum _{k=1}^{p-1}{p \choose k}
(b_1^p-\pi^{m_1'(p-1)}b_1)^{k}(a_1)^{p-k}\pi^{-m_1'p^2})$. 
By reducing modulo $\pi$ both equations (1) and (2) we obtain equations (1)' and (2)' for the cover $f_k$ on the level of special fibres. 
Suppose that $f$ has the structure of a 
torsor under a finite and flat $R$-group scheme $G_R$ of rank $p^2$. Then the special fibre $f_k$ of $f$ is a torsor under the special fibre $G_k$ of $G$ 
which is a group scheme of rank $p^2$ and an extension of $\alpha_p$ by $\alpha_p$. In particular both equations (1)' and (2)' that we obtain for $f_k$ must 
define the same $\alpha_p$-torsor $f_{1,k}$ and $f_{2,k}$. This is the case for $f_{1,k}$ since both equations are $\tilde t_1^p=\bar a_1$ and 
$\tilde s_1^p=\bar a_1+\bar b_1^p$ but we will see that this is not the case 
for $f_{2,k}$ if $p>2$. Indeed one can see after some (easy) computations that the equations we obtain for $f_{2,k}$ are: 
$\tilde t_2^p=-\tilde t_1^{p(p-1)+1}$ and $\tilde s_2^p=-\tilde s_1^{p(p-1)+1}+\bar b_1^{p(p-1)+1}+\sum _{k=1}^{p-1}(\bar b_1^{p(k-1)+1}(\tilde s_1-\bar 
b_1)^{p(p-k)}+\bar b_1^{pk}(\tilde s_1-\bar b_1)^{p(p-k-1)+1}))$, where $\tilde s_2$ is the image of $\pi^{\tilde m_1}S_2$ modulo $\pi$,
which gives $\tilde s_2=-\tilde s_1^{p(p-1)+1}+\bar b_1^{p(p-1)+1}+((\tilde s_1-\bar 
b_1)^{p-1}+\bar b_1^{p-1})\sum _{k=1}^{p-1}\bar b_1^{p(k-1)+1}(\tilde s_1-\bar b_1)^{p(p-k-1)+1}$. An easy verification shows that these two equations define 
non isomorphic $\alpha_p$-torsors (although they define isomorphic covers). 
If for example $p=3$ then the two equations are: $\tilde t_2^3=-\tilde t_1^7$ and $\tilde s_2^3=-\tilde s_1^7+\bar b_1^7+
((\tilde s_1-\bar b_1)^2-\bar b_1^2)(\bar b_1(\tilde s_1-\bar b_1)^4+\bar b_1^4(\tilde s_1-\bar b_1))$. Using that $\tilde s_1=\tilde t_1+\bar b_1$ we get 
$\tilde s_2^3= -(\tilde t_1+\bar b_1)^7+\bar b_1^7+(\tilde t_1^2-\bar b_1^2)(\bar b_1\tilde t_1^4+\bar b_1^4\tilde t_1)$ thus 
$\tilde s_2^3=-\tilde t_1^7-34\tilde t_1^3\bar b_1^4-8\tilde t_1\bar b_1^6=-\tilde t_1^7+
2\tilde t_1^3\bar b_1^4+\tilde t_1\bar b_1^6$ and since $2\tilde t_1^3\bar b_1^4+\tilde t_1\bar b_1^6=2\bar a_1\bar b_1^4+\tilde t_1\bar b_1^6$ 
is not a cube we deduce that the two equations do not define the same $\alpha_3$-torsor.

\pop {.5}
\par
\noindent
{\bf \gr 2.4.6. Corollary.}\rm \ {\sl Assume $p> 2$. Let $X$ be a formal affine $R$-scheme of finite type with $X_k$ integral. Then there exists a $p^2$-cyclic cover
$f:Y\to X$ which is \'etale on the generic fibre $X_K$ of $X$ such that the special fibre $Y_k$ of $Y$ is integral and $f$ doesn't have the structure of a torsor under 
a finite and flat $R$-group scheme of rank $p^2$.}

\pop {.5}
\par
\noindent
{\bf \gr Proof.}\rm\ One can easily find examples of cyclic $p^2$-covers $f:Y\to X$ which have the structure of an \'etale torsor above $X_K$ and such that
we are in the case c-1) of 2.4.3. The result follows then by 2.4.5.

\pop {.5}
\par
\noindent
{\bf \gr 2.4.7. Question.}\rm\ With the same notations as in 2.4.3 is it possible to find necessary and sufficient conditions, for example on $\delta_1$ 
and $\delta_2$, such that $f$ has the structure of a torsor under a finite and flat $R$-group scheme of rank $p^2$?

\pop {.5}
\par
Next we define the ``{\it degeneration data}'' arising from the reduction of an \'etale $\Bbb Z/p^2\Bbb Z$-torsor.

\pop {.5}
\par
\noindent
{\bf \gr 2.4.8. Definition.}\rm\ Let $f_K:Y_K\to X_K$ be an \'etale $\Bbb Z/p^2\Bbb Z$-torsor
with $X=\Spec A$ affine as in 2.4.3. Then we define the {\it degeneration type} of the torsor $f_K$ as follows:
$f_K$ has a degeneration of type A, or of type (\'etale, \'etale), if we are in the case a) of 2.4.3, a degeneration of type B, 
or of type (\'etale, radicial), if we are in the case b) of 2.4.3 and a degeneration of type C, or of type (radicial, radicial), if we are in the case c) 
of 2.4.3. Further we define the {\it degeneration data} associated to a degeneration type as follows:

\pop {.5}
\par
\noindent
{\bf a)}\ A degeneration data of type A consists of an element of $H^1_{\et}(X_k,\Bbb Z/p^2\Bbb Z)$.

\pop {.5}
\par
\noindent
{\bf b)}\ A degeneration data of type B consists of an element of $H^1_{\fppf}(X_k,G_k)$ where $G_k$ is the finite commutative group scheme
extension of $\Bbb Z/p\Bbb Z$ by $\alpha_p$ as defined in 2.3.1.

\pop {.5}
\par
\noindent
{\bf c)}\ A degeneration data of type C consists of an element of $H^1_{\fppf}(X_k,H_k)\oplus \Gamma (X_k,\Cal O_{X_k})^{p-1}$.
where $H_k$ is the finite commutative group scheme extension of $\alpha_p$ by $\alpha_p$ as defined in 2.3.1

\pop {.5}
\par
\noindent
{\bf \gr 2.4.9. Proposition.}\rm\ {\sl Assume that $X$ is affine as in 2.4.3 and let $f_K:Y_K\to X_K$ be an \'etale $\Bbb Z/p^2\Bbb Z$-torsor
which has a degeneration of type A (resp. B or C). Then $f_K$ induces canonically a degeneration data of type A (resp. of type B or C)}.

\pop {.5}
\par
\noindent
{\bf \gr Proof.}\rm\ This is a direct consequence of 2.4.3. If $f_K$ has a degeneration of type A then the special fibre $f_k$ of $f$ is 
an \'etale $\Bbb Z/p^2\Bbb Z$-torsor and the assertion follows in this case. Assume that $f_K$ has a degeneration of type B. Then 
the special fibre $f_k$ of $f$ is canonically a $G_k$-torsor and the assertion follows in this case too. Finally assume that the torsor $f_K$ has a 
degeneration of type C.
Then it follows from 2.4.3 that the special fibre $f_k$ of the cover $f$ is defined by the equations: $\tilde t_1^p=\bar a_1$ and 
$\tilde t_2^p=-\bar c_1^p\tilde t_1-2\bar c_2^p\tilde t_1^{p+1}-...-(p-1)\bar c_{p-1}^p\tilde t_1^{p(p-2)+1}-\tilde t_1^{p(p-1)+1}+\bar g$ 
(resp.  $\tilde t_1^p=\bar a_1$ and $\tilde t_2^p=-\bar c_1^p\tilde t_1-2\bar c_2^p\tilde t_1^{p+1}-...-(p-1)\bar c_{p-1}^p\tilde t_1^{p(p-2)+1}+\bar g$, or 
$\tilde t_1^p=\bar a_1$ and $\tilde t_2^p=\bar g$) where $\bar c_1$,...,$\bar c_{p-1}$ (resp. $\bar g$) are functions on $X_k$ 
(eventually equal to $0$) which are uniquely determined (resp. determined up to addition of element of the form $\bar h^p$ where 
$\bar h$ is a function on $X_k$). The pair $(\bar a_1, \bar g)$ defines then canonically an element of 
$H^1_{\fppf}(X_k,H_k)$ and the tuple $(\bar c_1,...,\bar c_{p-1})$ an element of $\Gamma (X_k,\Cal O_{X_k})^{p-1}$. 
Thus we get canonically an element of $H^1_{\fppf}(X_k,H_k)\oplus \Gamma (X_k,\Cal O_{X_k})^{p-1}$ associated to $f_K$ in this case.

\pop {.5}
\par
\noindent
{\bf \gr 2.4.10.}\rm\ It follows from 2.4.8 that an \'etale $\Bbb Z/p^2\Bbb Z$-torsor above the generic fibre $X_K$ of $X$
induces canonically a {\it degeneration data} of type either $A$, $B$ or $C$. Reciprocally we have the following result of {\it lifting} of such a 
degeneration data.

\pop {.5}
\par
\noindent
{\bf \gr 2.4.11. Proposition.}\rm\ {\sl Assume given a degeneration data, say $\Cal D$, of type either $A$, $B$ or $C$ as in 2.4.6. Then there exists a 
$\Bbb Z/p^2\Bbb Z$-torsor $f_K:Y_K\to X_K$ such that the degeneration data associated to $f_K$ via 2.4.8 equals $\Cal D$.}

\pop {.5}
\par
\noindent
{\bf \gr Proof.}\rm\ The proof in the case where the degeneration data is of type $A$ or $B$ is similar to the proof in 2.2.3. Assume that the degeneration data 
is of type $C$ and consists of the pair $(\bar a_1, \bar a_2)$, where $\bar a_1$ and  $\bar a_2$ are functions on $X_k$ which are not $p$-powers, and the
tuple of functions $(\bar c_1,...,\bar c_{p-1})$. Let $a_1$ and $a_2$ (resp. $c_1$,...,$c_{p-1}$) be regular functions on $X$ which lift 
$\bar a_1$ and  $\bar a_2$ (resp. which lifts $\bar c_1$,...,$\bar c_{p-1}$). Let $n=pn'=p^2n''>0$ be an integer. Consider the $\Bbb Z/p^2\Bbb Z$-torsor
given by the equations: $(T_1^p,T_2^p)-(T_1,T_2)=(a_1\pi^{-n'p},f(a_1)\pi^{-pm}+a_2\pi^{-pm+n'(p-1)})$ where $f(a_1)=c_1a_1+...+c_{p-1}a_1^{p-1}$ and 
$m=n'p$. Then it follows easily from the proof of 2.4.3 that the degeneration data associated to $f_K$ via 2.4.7 equals $D$. In this lifting we have
$\delta_1=n'(p-1)$ and $\delta_2=n''(p(p-1)+1)(p-1)$. We have also the other following possibility for the lifting. Namely consider the $\Bbb Z/p^2\Bbb Z$-torsor
given by the equations: $(T_1^p,T_2^p)-(T_1,T_2)=(a_1\pi^{-n'p},f(a_1)\pi^{-pm}+g\pi^{-pm+n'(p-1)})$ where $m$ is a positive integer such that 
$mp-n'(p-1)>n'(p(p-1)+1)$ and $mp-n'(p-1)=pm'$. In this later case we have $\delta_1=n'(p-1)$ and $\delta_2=m'(p-1)$.

\pop {2}
\par
\noindent
{\bf \gr 2.5. Degeneration of $\Bbb Z/p\Bbb Z$-torsors on the boundaries of formal fibres}.\rm \

\pop {.5}
\par
In this section we assume that the {\it residue field $k$ of $R$ is algebraically closed.}
In what follows we explain the degeneration of $\Bbb Z/p\Bbb Z$-torsors on the boundary $\Cal X \simeq 
\Spf R[[T]]\{T^{-1}\}$ of formal fibres of germs of formal $R$-curves. Here $R[[T]]\{T^{-1}\}$ denotes the ring of 
formal power series $\sum_{i\in \Bbb Z} a_iT^{i}$ with $\lim _{i\to -\infty} \vert a_i \vert =0$ where 
$\vert\ \vert$ is an absolute value of $K$ associated to its valuation. Note that $R[[T]]\{T^{-1}\}$ is a 
complete discrete valuation ring with uniformising parameter $\pi$ and residue field $k((t))$ where 
$t\equiv T\mod \pi$. The function $T$ is called a {\it parameter} of the formal fibre $\Cal X$. 
The following result, which describes the degeneration of $\Bbb Z/p\Bbb Z$-torsors above the formal fibre
$\Spf R[[T]]\{T^{-1}\}$, will be used in the next paragraph III in order to prove a formula comparing the 
dimensions of the spaces of vanishing cycles in a Galois cover of degree $p$ between formal germs of $R$-curves.

\pop {.5}
\par
\noindent
{\bf \gr 2.5.1. Proposition.}\rm \ {\sl Let $A:=R[[T]]\{T^{-1}\}$ and let
$f:\Spf B\to \Spf A$ be a non trivial Galois cover of degree $p$.
Assume that the ramification index of the corresponding extension of
discrete valuation rings equals $1$. Then $f$ is a torsor under a
finite and flat $R$-group scheme $G_R$ of rank $p$. Let $\delta$ be
the degree of the different in the above extension. Then the following cases
occur:
\par
{\bf a )} $\delta=0$. In this case $f$ is a torsor under the \'etale group 
$(\Bbb Z/p\Bbb Z)_R$ and for a suitable choice of the parameter $T$ of $A$
the torsor $f$ is given by an equation $X^p-X=T^{m}$ for some negative integer $m$ 
which is prime to $p$. In this case $X^{1/m}$ is a parameter for $B$.
\par
{\bf b )} $0 <\delta=n(p-1)$ for some positive integer $n$. In this case $f$ is a torsor under the group
scheme $\Cal M_{n,R}$. Moreover for a suitable choice of the parameter $T$ the torsor $f$ is
given by an equation $X^p-\pi^{n(p-1)}X=T^m$ with $m\in \Bbb Z$ is prime to $p$. In this case $X^{1/m}$
is a parameter for $B$.

\pop {.5}
\par
\noindent
{\bf \gr Proof.}\rm\ First it follows from 2.2.1 that we are either in case a) or in case b). 
We start first with the case a). In this case $f$ is an \'etale torsor given by an equation
$X^p-X=u=\sum _{i\in \Bbb Z} a_iT^i\in A$. On the level of special fibres the torsor 
$f_k:=\Spec B/\pi B\to \Spec A/\pi A$ is the \'etale torsor given by the equation
$x^p-x=\sum _{i\ge m} \bar a_it^i\in A$ where $\bar a_i$ is the image of $a_i$ modulo $\pi$
and $m\in \Bbb Z$ is some integer. Assume for example that the integer $m=pm'$ is 
divisible by $p$. Then after adding $\bar a_m^{1/p}t^{m'}-\bar a_mt^m$ into the defining equation for $f_k$ we can replace $\bar a_mt^m$ 
by $\bar a_m^{1/p}t^{m'}$. Repeating this process eventually we can finally assume that the integer 
$m$ is prime to $p$ in which case $\sum _{i\ge m} \bar a_it^i=t^m\bar v$ and $u=T^mv$ 
where $v\in A$ is a unit whose image modulo $\pi$ equals $\bar v$.
Further the integer $m$ is then necessarily negative since the residue field extension
$f_k:=\Spec B/\pi B\to \Spec A/\pi A$ must ramify. Finally after extracting an $m$-th root 
of $v$ and replacing $T$ by $Tv^{1/m}$ we arrive to an equation of the form $X^p-X=T^m$.
Secondly assume that we are in the case b). Then $f$ is a torsor under the finite and flat 
group scheme $\Cal M_{n,R}$, for some positive integer $n$, given by an equation 
$X^p-\pi^{n(p-1)}X=u=\sum _{i\in \Bbb Z} a_iT^i\in A$ and $u$
is not a $p$-power modulo $\pi$. In this case on the level of special fibres the torsor 
$f_k:=\Spec B/\pi B\to \Spec A/\pi A$ is the $\alpha_p$-torsor given by the equation
$x^p=\sum _{i\ge m} \bar a_it^i\in A/\pi A$ where $\bar a_i$ is the image of $a_i$ modulo $\pi$
and $m\in \Bbb Z$ is some integer. Assume for example that the integer $m=pm'$ is 
divisible by $p$. Then the term $\bar a_mt^m$ is a $p$-power and 
we can eliminate it from the defining equation for $f_k$ and since $\sum _{i\ge m} \bar a_it^i\in A/\pi A$
is not a $p$-power we can repeat this process eventually and assume after finitely many steps
that $m$ is prime to $p$. In this case $\sum _{i\ge m} \bar a_it^i=t^m\bar v$ and 
$u=T^mv$ where $v\in A$ is a unit whose image modulo $\pi$ equals $\bar v$.
Finally after extracting an $m$-th root of $v$ and replacing $T$ by $Tv^{1/m}$ we arrive to an 
equation of the form $X^p-\pi^{n(p-1)}X=T^m$.

\pop {.5}
\par
\noindent
{\bf \gr 2.5.2. Definition.}\rm \ With the same notations as in 2.5.1 we define the {\it conductor}
of the torsor $f$ to be the integer $-m$. Further we define the {\it degeneration type}
of the torsor $f$ to be $(0,m)$ in the case a) and $(n,m)$ in the case b).

\pop {.5}
\par
\noindent
{\bf \gr 2.5.3. Remark.}\rm \ The above proposition implies in particular that Galois covers $f:\Spf B\to \Spf A$
as in 2.5.1 are classified by their degeneration type as defined in 2.5.2. More precisely given two such Galois covers which have the same
degeneration type then there exists a (non canonical) Galois equivariant isomorphism between both covers.

\pop {2}
\par
\noindent
{\bf \gr 2.6. Degeneration of $\Bbb Z/p^2\Bbb Z$-torsors on the boundaries of formal fibres}.\rm \

\pop {.5}
\par
We use the same notations and hypothesis as in 2.5. Next we explain the degeneration of 
\'etale $\Bbb Z/p^2\Bbb Z$-torsors on the boundaries of formal fibres.

\pop {.5}
\par
\noindent
{\bf \gr 2.6.1. Proposition.}\rm \ {\sl Let $A:=R[[T]]\{T^{-1}\}$ and let
$\tilde f:\Cal Y:=\Spf B\to \Cal X:=\Spf A$ be a non trivial cyclic Galois cover of degree $p^2$ with $\Cal Y_k$ irreducible.
Assume that the ramification index of the corresponding extension of
discrete valuation rings equals $1$. Then $\tilde f$ factorises canonically 
as $\tilde f=f_1\circ f_2$ with $f_2:\Cal Y\to \Cal Y_1:=\Spf B_1$ and $f_1:\Cal Y_1\to X$
where $f_i$ is a Galois cover of degree $p$. Let $\delta$ (resp. $\delta_1$
and $\delta_2$) be the degree of the different in the above extension (resp. 
in $f_1$ and $f_2$). Then the following cases occur:

\par
{\bf a )} $\delta=0$. In this case $\tilde f$ is a torsor under the \'etale $R$-group 
$G=(\Bbb Z/p^2\Bbb Z)_R$ and for a suitable choice of the parameter $T$ of $A$
the torsor $\tilde f$ is given by an equation: 
$$(X_1^p,X_2^p)-(X_1,X_2)=(1/T^{m_1},f(T))$$ 
where $m_1$ is a 
positive integer prime to $p$ and $f(T)\in A$ is such that its image $\bar f(t)$ modulo $\pi$ equals 
$\sum _{i\ge -m_2}c_it^i$ where $m_2\ge 0$ is prime to $p$ and $c_{-m_2}\neq 0$. In this case the torsor $f_1$ (resp. $f_2$) has a reduction of type $(0,-m_1)$ 
(resp. $(0,\tilde m_2:=-pm_2+m_1(p-1))$ if $m_2\ge pm_1$ and $(0,\tilde m_2:=-m_1(p(p-1)+1))$ otherwise).

\par
{\bf b)} $\delta_1=0$ and $\delta_2>0$. In this case $\tilde f$ is a torsor under 
the group scheme $W_{(0,n)}$ for some positive integer $n$ and for a suitable choice of the parameter 
$T$ of $A$ the cover $\tilde f$ is generically given by the equations: 
$$X_1^p-X_1=1/T^{m_1}$$ 
and: 
$$X_2^p-X_2=f(T)/\pi ^{pn}+p^{-1}\sum _{k=1}^{p-1}{p \choose k}X_1^{pk}(-X_1)^{p-k}$$
with $f(T)\in A$ is such that its image $\bar f(t)=\sum _{i\in I}c_it^i$ modulo $\pi$ 
is not a $p$-power.  Let $m_2:=\inf \{i\in I,\gcd(i,p)=1, c_i\neq 0\}$. We call the integer $m_2$ the {\it conductor} of $f(T)$. Then the torsor 
$f_1$ (resp. $f_2$) has a reduction of type $(0,-m_1)$ (resp. $(n,\tilde m_2)$ where $\tilde m_2:=pm_2+m_1(p-1)$).

\par
{\bf c)} $\delta_1<0$. In this case the cover $\tilde f$ is generically given for a suitable choice of the parameter $T$ by the equations: 
$$X_1^p-X_1=T^{m_1}/\pi ^{n_1p}$$ 
for some integer $m_1$ prime to $p$ and: 
$$X_2^p-X_2=f(T)/\pi^{pn_2}+p^{-1}
\sum _{k=1}^{p-1}{p \choose k}X_1^{pk}(-X_1)^{p-k}$$ 
where $f(T)=\sum _{i\in I} a_iT^i\in A$
is not a $p$-power modulo $\pi$. Let $m_2:=\inf \{i\in I,\gcd(i,p)=1$, and $a_i$ is a unit$\}$ be the {\it conductor} 
of $f$. Then the following distinct cases occur:

\par
{\bf c-1)} $pn_2-(p-1)n_1>n_1(p(p-1)+1)$. Then $pn_2-(p-1)n_1=pn_2'$ is divisible by $p$ and the torsor $f_1$ (resp. $f_2$) has a 
reduction of type $(n_1,m_1)$ (resp. $(n_2',m_2p-m_1(p-1)))$.

\par
{\bf c-2)} $pn_2-(p-1)n_1<n_1(p(p-1)+1)$. Then $n_1(p(p-1)+1)=pn_2''$ is divisible by $p$ and the torsor $f_1$ (resp. $f_2$) has a 
reduction of type $(n_1,m_1)$ (resp. $(n_2'',m_1(p(p-1)+1))$).

\par
{\bf c-3)} $pn_2-(p-1)n_1=n_1(p(p-1)+1)$. Then $n_1(p(p-1)+1)=pn_2''$ is divisible by $p$ and the torsor $f_1$ (resp. $f_2$) has a 
reduction of type $(n_1,m_1)$ (resp. $(n_2'',\tilde m_2)$ where $\tilde m_2=\inf (m_1(p(p-1)+1),m_2p-m_1(p-1)))$).

\pop {.5}
\par
\noindent
{\bf \gr Proof.}\rm\ First it follows from 2.4.3 that the above cases are all the possible cases of 
degeneration of the $\Bbb Z/p^2\Bbb Z$-torsor $f_K$. We start with the case a) and assume that $\delta=0$. 
In this case $\tilde f$ is an \'etale torsor given by an equation $(X_1^p,X_2^p)-(X_1,X_2)=(g(T),f(T))$
where $g(T)$ (resp. $f(T)$) are elements of $A$. Moreover by 2.5.1 we can assume that for a suitable 
choice of the parameter $T$ we have $g(T)=1/T^{m_1}$ for some positive integer $m_1$ prime to $p$. One has to compute the 
conductors of the torsors $f_1$ and $f_2$. As for $f_1$ we have $T=X_1^{-p/m_1}(1-X_1^{(1-p)})^{-1/m_1}$ from which 
it follows that $Z:=X_1^{-1/m_1}$ is a parameter for $\Cal Y_1$ and we have $T=Z^p(1-Z^{m_1(p-1)})^{-1/m_1}$. 
The torsor $f_2$ is given by 
the equation $X_2^p-X_2=f(T)+p^{-1}\sum _{k=1}^{p-1}{p \choose k} X_1^{pk}(-X_1)^{p-k}$ and after adding elements of the form $b^p-b$
where $b\in A$ we can assume that the image $\bar f(t)$ of $f(T)$ modulo $\pi$ equals $\sum _{i\ge -m_2}c_it^i$ where 
$m_2\ge 0$ is prime to $p$. Now we have $1/t^{m_2}=1/z^{pm_2}(1-z^{m_1(p-1)})^{m_2/m_1}=1/z^{pm_2}-
(m_2/m_1)1/z^{pm_2-m_1(p-1)}+...$,
where $z$ is the image of $Z$ modulo $\pi$, which can be replaced after some transformation by $1/z^{m_2}-
(m_2/m_1)1/z^{pm_2-m_1(p-1)}+...$. On the other hand the summand in $p^{-1}\sum _{k=1}^{p-1}{p \choose k}X_1^{pk}(-X_1)^{p-k}$ which gives 
the highest contribution to the different is $-X_1^{p(p-1)+1}=1/z^{m_1(p(p-1)+1)}$. From this we deduce that the conductor for the 
torsor $f_2$ equals $-pm_2+m_1(p-1)$ if $m_2\ge pm_1$ and $-m_1(p(p-1)+1)$ otherwise.
\par
Assume next that we are in the case b). In this case $\tilde f$ is generically given by 
an equation $(X_1^p,X_2^p)-(X_1,X_2)=(g(T),f(T)/\pi^{pn})$
where $g(T)$ (resp. $f(T)$) are elements of $A$, $n$ is a positive integer, and the image of $f(T)$ modulo $\pi$ is not a $p$-power.
Moreover by 2.5.1 we can assume that for a suitable choice of the parameter $T$ we have $g(T)=1/T^{m_1}$ for some positive integer 
$m_1$ prime to $p$. Consider the image $\bar f(t)$ of $f(T)$ modulo $\pi$ which equals $\sum _{i\ge m_2}c_it^i$. Then we can assume 
after some transformations which eliminate the $p$-powers that the integer $m_2$ is prime to $p$. The special 
fibre of the torsor $f_2$ is given
by the equation: $\tilde x_2^p=\bar f(t)=\sum _{i\ge m_2}c_it^i$. Now $t^{m_2}=z^{pm_2}(1-z^{m_1(p-1)})^{-m_2/m_1}=z^{pm_2}
+(m_2/m_1)z^{pm_2+m_1(p-1)}+...$ which can be transformed, after eliminating the term 
$z^{pm_2}$ which is a $p$-power, to $(m_2/m_1)z^{pm_2+m_1(p-1)}+...$. From this we deduce that the conductor of $f_2$ equals $pm_2+m_1(p-1)$.
\par
Finally assume that we are in the case c). In this case the cover $\tilde f$ is given by 2.4.3 and 2.5.1, and for 
a suitable choice of the parameter $T$, by the equations: 
$$X_1^p-\pi^{n_1(p-1)}X_1=T^{m_1}$$ 
and 
$$X_2^p-X_2=f(T)/\pi^{pn_2}+p^{-1}\sum _{k=1}^{p-1}{p \choose k}X_1^{pk}(-X_1)^{p-k}\pi^{-n_1(p(p-k)+k)}$$
where $f(T)=\sum _{i\in I} a_iT^i\in A$. Note that the image $\bar f(T)=\sum _{i\in I} \bar a_it^i$ modulo $\pi$ 
is necessarily a $p$-power in $B_1/\pi B_1$ since
the image $t$ of $T$ is. We write $f(T)=f((X_1^p-\pi^{n_1(p-1)}X_1)^{1/m_1})=\sum _{i\in I} a_i(X_1^p-\pi^{n_1(p-1)}X_1)^{i/m_1}$.
\par
The expression $f((X_1^p-\pi^{n_1(p-1)}X_1)^{1/m_1})$ can be transformed after addition of elements of the form $b^p-b$ (where $b\in B_1$) to:
$$\sum _ia_i^{1/p}X_1^{i/m_1}\pi ^{-n_2}-\sum _i (i/{m_1})a_iX_1^{ip/m_1}X_1^{1-p}\pi^{-pn_2+n_1(p-1)}+g(X_1)\pi^{-n'}$$
where $g(X_1)\in B_1$ and $n'<pn_2-n_1(p-1)$ is a positive integer. Assume first that $n_2p-n_1(p-1)>n_1(p(p-1)+1)$ 
then $n_2p-n_1(p-1)=n_2'p$ is necessarily divisible by $p$ and 
$f_2$ is a torsor under the group scheme $\Cal M_{n_2',R}$. Its special fibre $f_{2,k}$ is the $\alpha_p$-torsor given by the equation
$\tilde x_2^p=-\sum _{i\ge m_2} (i/{m_1})\bar a_ix_1^{ip/m_1}x_1^{1-p}$ where $\bar a_i$ (resp. $x_1$) is the image of $a_i$ (resp. $X_1$)
modulo $\pi$ and $m_2:=\inf \{i\in I,\gcd(i,p)=1$, and $a_i$ is a unit$\}$ is the conductor of $f$. 
Since $x_1^{1/m_1}$ is a parameter for $B_1/\pi B_1$ we deduce immediately that the conductor 
for the torsor $f_2$ equals $m_2p-m_1(p-1)$. Next assume that $n_2p-n_1(p-1)<n_1(p(p-1)+1)$ then $n_1(p(p-1)+1)=n_2''p$ is 
necessarily divisible by $p$ and 
$f_2$ is a torsor under the group scheme $\Cal M_{n_2'',R}$. Its special fibre $f_{2,k}$ is the $\alpha_p$-torsor given by the equation
$\tilde x_2^p=-x_1^{p(p-1)+1}$ from which we deduce immediately that the conductor of $f_2$ equals $m_1(p(p-1)+1)$. Finally consider the 
last case where $n_2p-n_1(p-1)=n_1(p(p-1)+1)$ in which case $n_1(p(p-1)+1)=n_2''p$ is necessarily divisible by $p$ and 
$f_2$ is a torsor under the group scheme $\Cal M_{n_2'',R}$. Its special fibre $f_{2,k}$ is the $\alpha_p$-torsor given by the equation
$\tilde x_2^p=-x_1^{p(p-1)+1}-\sum _{i\ge m_2}(i/m_1)\bar a_ix_1^{ip/m_1}x_1^{1-p}$ from which we deduce that the conductor of $f_2$
equals $\inf(m_1(p(p-1)+1),m_2p-m_1(p-1))$ which equals $m_1(p(p-1)+1)$ if $m_2\le pm_1$ and $m_2p-m_1(p-1)$ otherwise.

\pop {.5}
\par
\noindent
{\bf \gr 2.6.2. Definition.}\rm\ With the same notations as in 2.6.1 let $(n_1,m_1)$ (resp. $(n_2,m_2)$) be the degeneration type of
the torsor $f_1$ (resp. $f_2$). Then we define the degeneration type of the cover $\tilde f$ to be $\{(n_1,m_1),(n_2,m_2)\}$. Note that the inequality
$n_2\ge n_1(p(p-1)+1)/p$ holds by 2.4.3.

\pop {.5}
\par
\noindent
{\bf \gr 2.6.3. Remark.}\rm\ Contrary to what happend for $p$-cyclic covers (cf. Remark 2.5.3) it is no more true that $p^2$-cyclic 
covers above formal boundaries are determined by their degeneration type as defined in 2.6.2.

\pop {.5}
\par
\noindent
{\bf \gr 2.6.4. Definition.}\rm\ Let $\{(n_1,m_1),(\tilde n_2,\tilde m_2)\}$ be a pair of pair of integers. We say that this pair is {\it admissible} if
the following holds: either $n_1=\tilde n_2=0$ in which case  $m_1$ and $\tilde m_2$ are negative and $\tilde m_2=\inf(m_1(p(p-1)+1),m_2p-m_1(p-1))$
for some negative integer $m_2$ prime to $p$ or $n_1=0$ and $\tilde n_2\neq 0$ in which case $m_1\le -1$ and $\tilde m_2=pm_2+m_1(p-1)$ for some integer $m_2$ prime to $p$
or finally $n_1\neq 0$ in which case either $\tilde n_2> n_1(p(p-1)+1)/p$ and $\tilde m_2=pm_2-m_1(p-1)$ for some integer $m_2$ prime to $p$ 
or $\tilde n_2=n_1(p(p-1)+1)/p$ in which case $\tilde m_2=m_1(p(p-1)+1)$ or $\tilde m_2=m_2p-m_1(p-1)$ for some integer $m_2$ prime to $p$ with $m_2<pm_1$. 
It follows from 2.6.1 that the degeneration type $\{(n_1,m_1),(\tilde n_2,\tilde m_2)\}$ of the cover $\tilde f$ is an admissible pair.

\pop {1}
\par
\noindent
{\bf \gr III. Computation of vanishing cycles and examples for cyclic $p$-covers.}
\pop {.5}
\par
The main result of this section is Theorem 3.2.3 which gives a {\bf formula} which compares the dimensions of the spaces 
of vanishing cycles in a Galois cover $\Tilde f:\Cal Y\to \Cal X$ with group $\Bbb Z/p\Bbb Z$ between formal
germs of $R$-curves where $R$ is a complete discrete valuation ring of equal characteristic $p>0$ in terms of the degeneration 
type of $\tilde f$ above the boundaries of $\Cal X$ as defined in 2.5.2. This formula enables one in principle to compare the dimensions 
of the spaces of vanishing cycles in a Galois cover $\Tilde f:\Cal Y\to \Cal X$ with group a $p$-group as we will illustrate
for the case of cyclic covers of degree $p^2$ in IV. In all this section we use the following notations:
$R$ is a complete discrete valuation ring of equal characteristic $p>0$. We denote by $K$ 
the fraction field of $R$ by $\pi $ a uniformising parameter and $k$ the residue field. 
We also denote by $v_K$ the valuation of $K$
which is normalised by $v_K(\pi)=1$. We assume also that the residue field $k$ is {\it algebraically closed}.

\pop {.5}
\par
\noindent
{\bf \gr 3.1.}\rm \ By a (formal) $R$-curve we mean a (formal) $R$-scheme of finite type which is normal flat and whose 
fibres have dimension $1$. For an $R$-scheme $X$ we denote by $X_K:=X\times _{\Spec R}\Spec K$ the {\it generic}
fibre of $X$ and $X_k:=X\times _{\Spec R}\Spec k$ its {\it special} fibre.
In what follows by a {\it (formal) germ} $\Cal X$ of an $R$-curve we mean that $\Cal X :=\Spec \Cal O_{X,x}$ is the 
(resp. $\Cal X:=\Spf \hat {\Cal O}_{X,x}$ is the formal completion of the) spectrum of the local ring of an
$R$-curve $X$ at a closed point $x$. We refer to [S-1] 3.1 for the definition of the integers $\delta_x$, $r_x$,
and the genus $g_x$ of the point $x$.

\pop {.5}
\par
\noindent
{\bf \gr 3.2. The compactification process.}\rm \ Let $\Cal X:=\Spf \hat {\Cal O}_{X,x}$ be the formal germ of an $R$-curve at a
closed point $x$ with $\Cal X_k$ reduced. Let $\Tilde f:\Cal Y\to \Cal X$ be a Galois cover with group 
$\Bbb Z/p\Bbb Z$ and $\Cal Y$ local. We assume 
that the special fibre of $\Cal Y_k$ is reduced (this can always be achieved after a finite extension of $R$).
We will construct a compactification of the above cover $\Tilde f$ and as an application we will use this compactification in 
order to compute the arithmetic genus of the closed point of $\Cal Y$. More precisely we will construct a Galois cover 
$f:Y\to X$ of degree $p$ between proper algebraic $R$-curves, a closed point $y\in Y$ and its image $x=f(y)$, such that the 
formal germ of $X$ (resp. of $Y$) at $x$ (resp. at $y$) equals $\Cal X$ (resp. $\Cal Y$) and such that the Galois cover
$f_x:\Spf {\hat {\Cal O}}_{Y,y}\to \Spf {\hat {\Cal O}}_{X,x}$ induced by $f$ between the formal germs at
$y$ and $x$ is isomorphic to the above given cover $\tilde f:\Cal Y\to \Cal X$. The construction of such a compactification 
has been done in [S-1] 3.3.1 in the inequal characteristic case. we first start with the case where the formal germ $\Cal X$ has
only one boundary.

\pop {.5}
\par
\noindent
{\bf \gr 3.2.1. Proposition.}\rm \ {\sl Let $D:=\Spf R<1/T>$ be the formal closed disc centered at $\infty$ (cf. [B-L], 1, for the 
definition of $R<1/T>$). Let $\Cal D:=\Spf R[[T]]\{T^{-1}\}$ and let $\Cal D\to D$ be the canonical morphism. Let $\tilde
f:\Cal Y\to \Cal D$ be a non trivial torsor under a finite and flat $R$-group scheme of rank $p$ such that
the special fibre of $\Cal Y$ is reduced. Then there exists a Galois cover $f:Y\to D$ with group $\Bbb Z/p\Bbb Z$ whose pull back
to $\Cal D$ is isomorphic to the above given torsor $\Tilde f$. More precisely with the same notations introduced in 
2.5 we have the following possibilities:
\par
{\bf a)}\ The torsor $\tilde f$ is \'etale and has a reduction of type $(0,-m)$. In this case consider the Galois cover $f:Y\to D$ given
generically by the equation $Z^p-Z=1/{T^m}$. This cover is an \'etale torsor and its special fibre $f_k:Y_k\to X_k$ is \'etale and 
$Y_k$ is smooth. Moreover the genus of the smooth compactification of $Y_k$ equals $(m-1)(p-1)/2$. 
\par
{\bf b)}\ The cover $\tilde f$ is a torsor under the group scheme $\Cal M_{n,R}$ for some positive integer $n$ and has a reduction
of type $(n,m)$ for some integer $m$ prime to $p$. Consider the following two cases:
\par
{\bf b-1)}\ $m>0$. In this case consider the Galois cover $f:Y\to D$ given generically by the equation 
$Z^p-\pi^{n(p-1)}Z=T^m$. This cover is ramified above $\infty $ with conductor $m$ and its special fibre $f_k:Y_k\to X_k$ is
radicial. Moreover $Y_k$ is smooth and its smooth compactification has genus $0$.
\par
{\bf b-2)}\ $m<0$. In this case consider the Galois cover $f:Y\to D$ given generically by the equation $Z^p-\pi^{n(p-1)}Z=T^m$.
This cover is an \'etale torsor on the generic fibre and its special fibre $f_k:Y_k\to X_k$ is radicial. Moreover $Y_k$ has a unique singular point $y$
which is above $\infty $ and $g_y=(-m-1)(p-1)/2$.}

\pop {.5}
\par
\noindent
{\bf \gr Proof.}\rm\ The proof is similar to the proof of proposition 3.3.1 in [S-1]. For the convenience of the reader
we treat the case b-1). In this case consider the Galois $p$-cover $H\to \Bbb P^1_R$, with $H$ normal, above the projective 
$R$-line with parameter $T$ defined generically by the equation $Z^p-\pi^{n(p-1)}Z=T^m$. This cover is ramified on the generic 
fibre only above the point $\infty$ with conductor $m$ from which we deduce that the genus of the generic fibre $H_K$ of $H$ 
equals $(m-1)(p-1)/2$. On the level of the special fibres the cover $H_k\to\Bbb P^1_k$ is an $\alpha_p$ torsor outside the point
$\infty$ defined by the equation $z^p=t^m$. The genus of the singularity above the point $t=0$ can be then easily calculated and 
equals $(m-1)(p-1)/2$ (cf. [Sa-1] 3.3.1). From this we deduce that $H_k$ is smooth outside $t=0$ since the arithmetic genus of
$H_K$ and $H_k$ are equal.

\pop {.5}
\par
In the next proposition we deal with the general case.
\pop {.5}
\par
\noindent
{\bf \gr 3.2.2. Proposition}. \rm\ {\sl Let $\Cal X:=\Spf {\hat \Cal O_x}$ be the formal germ of an $R$-curve at a closed point $x$ 
and let $\{\Cal X_i\}_{i=1}^n$ be the boundaries of $\Cal X$. Let $\tilde f:\Cal Y\to \Cal X$ be a Galois cover with group 
$\Bbb Z/p\Bbb Z$ and with $\Cal Y$ local. Assume that $\Cal Y_k$ and $\Cal X_k$ are reduced. Then there exists a Galois cover 
$f:Y\to X$ of degree $p$ between proper algebraic $R$-curves $Y$ and $X$, a closed point $y\in Y$ and its image $x=f(y)$,
such that the formal germ of $X$ (resp. of $Y$) at $x$ (resp. at $y$) equals $\Cal X$ (resp. equals $\Cal Y$) and such
that the Galois cover $\Spf {\hat {\Cal O}}_{Y,y}\to \Spf {\hat {\Cal O}}_{X,x}$ induced by $f$ between the formal germs at
$y$ and $x$ is isomorphic to the above given cover $\tilde f:\Cal Y\to \Cal X$. Moreover the formal completion of $X$ along its special fibre
has a covering which consists of $n$ closed formal discs $D_i$ which are patched with $\Cal X$ along the boundaries $\Cal D_i$ 
and the special fibre $X_k$ of $X$ consists of $n$ smooth projective lines which intersect at the point $x$. In particular the
arithmetic genus of $X_K$ equals $g_x$.}

\pop {.5}
\par
\noindent
{\bf \gr Proof.}\rm\ Similar to the proof of proposition 3.3.2 in [S-1].

\pop {1}
\par
The next result is the main one of this section. It provides an explicit {\bf formula} which compares the dimensions of the spaces of
vanishing cycles in a Galois cover of degree $p$ between formal fibres of curves in equal characteristic $p>0$.

\pop {.5}
\par
\noindent
{\bf \gr 3.2.3. Theorem.}\rm \ {\sl Let $\Cal X:=\Spf {\hat {\Cal O}_x}$ be the formal germ of an $R$-curve at a closed point $x$ 
with $\Cal X_k$ reduced. Let $\tilde f:\Cal Y\to \Cal X$ be a Galois cover with group $\Bbb Z/p\Bbb Z$ with $\Cal Y$ local and
$\Cal Y_k$ reduced. Let $\{\wp_i\}_{i\in I}$ be the minimal prime ideals of $\hat {{\Cal O}}_x$ which contain $\pi$ and let
$\Cal X_{i}:=\Spf \hat \Cal O_{\wp _i}$ be the formal completion of the localisation of $\Cal X$ at $\wp _i$. For each $i\in I$ the above
cover $\tilde f$ induces a torsor $\tilde f_i:\Cal Y_i\to \Cal X_{i}$ under a finite and flat $R$-group scheme of rank $p$ above the boundary
$\Cal X_i$ (cf. 2.5.1). Let $(n_i,m_i)$ be the reduction type of $\tilde f_i$ (cf. 2.5.2). Let $y$ be the closed point of $\Cal Y$. Then 
one has the following {\bf ``local Riemann-Hurwitz formula''}:
$$2g_y-2=p(2g_x-2)+d_{\eta}-d_s$$
Where $d_{\eta}$ is the degree of the divisor of ramification in the morphism $\tilde f_K:\Cal Y_K\to \Cal X_K$ induced by $\tilde f$ 
on the generic fibres, where $\Cal X_K:=\Spec ({\hat {\Cal O}_x} \otimes _R K)$ and $\Cal Y_K:=\Spec ({\hat {\Cal O}_{\Cal Y,y}}
\otimes _R K)$, and $d_s:=\sum _{i\in I_{\et}}(-m_i-1)(p-1)+\sum _{i\in I_{\rad}}(-m_i-1)(p-1)$ where $I^{\rad}$ is the subset of 
$I$ consisting of those $i$ for which $n_i\neq 0$ and $I^{\et}$ is the subset of $I$ consisting of those $i$ for which $n_i=0$ and $m_i\neq 0$.

\pop {.5}
\par
\noindent
{\bf \gr Proof.}\rm \ The proof is similar, using 3.2.1,  to the proof of theorem 3.4 in [S-1] with the appropriate modifications.
We repeat briefly the argument for the convenience of the reader. By Proposition 3.2.2 one can compactify the
above morphism $\tilde f$. More precisely there we constructed  a Galois cover $f:Y\to X$ of degree $p$ between proper algebraic 
$R$-curves a closed point $y\in Y$ and its image $x=f(y)$ such that the formal germ of $X$ (resp. of $Y$) at $x$ (resp. at $y$) 
equals $\Cal X$ (resp. equals $\Cal Y$) and such that the Galois cover $\Spf {\hat {\Cal O}}_{Y,y}\to \Spf {\hat {\Cal
O}}_{X,x}$ induced by $f$ between the formal germs at $y$ and $x$ is isomorphic to the given cover $\tilde f:\Cal Y\to \Cal X$. 
The special fibre of $\Cal X$ consists (by construction) of $\card (I)$-distinct smooth projective lines which intersect at the closed 
point $x$. The formal completion of $X$ along its special fibre has a covering which consists of $\card (I)$ formal closed unit discs 
which are patched with the formal fibre $\Cal X$ along the boundaries $\Cal X_i$. The above formula follows then by comparing the 
arithmetic genus of the generic fibre $Y_K$ of $Y$ and the arithmetic genus of its special fibre $Y_k$. Using the precise
informations given in Proposition 3.2.1 one can easily deduce that $g(Y_K)=p g_x+(1-p)+d_{\eta}/2+\sum _{i\in I_{>}} (m_i+1)(p-1)/2$ where
$I_>$ is the subset of $I$ consisting of those $i$ for which the degeneration type above the boundary $\Cal X_i$ is $(n_i,m_i)$
with $n_i>0$ and $m_i>0$. On the other hand one has $g(Y_k)=g_y+\sum _{i\in I_{<}}(-m_i-1)(p-1)/2 +
\sum _{i\in I_{\et}} (-m_i-1)(p-1)$, where $I_{<}$ is the subset of $I$ consisting of those $i$ for which the degeneration
type above the boundary $\Cal X_i$ is $(n_i,m_i)$ with $n_i>0$ and $m_i<0$, and $I_{\et}$ is the subset of $I$ consisting of those $i$ 
for which the degeneration type above the boundary $\Cal X_i$ is $(0,m_i)$. Now since $Y$ is flat we obtain $g(Y_K)=g(Y_k)$ and the above 
formula directly follows.

\pop {.5}
\par
\noindent
{\bf \gr 3.3.\ $p$-Cyclic covers above germs of semi-stable curves.}\rm
\par
In what follows and as a consequence of theorem 3.2.3 we will deduce some results 
in the case of a Galois cover $\Cal Y\to \Cal X$ where $\Cal X$ is the formal germ of a {\it semi-stable} $R$-curve at a closed point. 
These results will play an important role in V in order to exhibit and realise the degeneration data
which describe the semi-stable reduction of Galois covers of degree $p$ in equal characteristic $p$. We start with the case of
a Galois cover of degree $p$ above a germ of a smooth point.
\pop {.5}
\par
\noindent
{\bf \gr 3.3.1. Proposition.}\rm\ {\sl Let $\Cal X:=\Spf R[[T]]$ be the formal germ of an $R$-curve at a {\bf smooth} point $x$ and let
$\Cal X_{\eta}:=\Spf R[[T]]\{T^{-1}\}$ be the boundary of $\Cal X$. Let $f:\Cal Y\to \Cal X$ be a Galois cover of degree $p$
with $\Cal Y$ local. Assume that the special fibre of $\Cal Y$ is reduced. Let $y$ be the unique closed point of $\Cal Y_k$.
Let $d_{\eta}$ be the degree of the divisor of ramification in the morphism $f:\Cal Y_K\to \Cal X_K$. Then
$d_{\eta}=r(p-1)$ is divisible by $p-1$. We distinguish two cases:
\par
1)\ $\Cal Y_k$ is {\bf unibranche} at $y$. Let $(n,m)$ be the degeneration type of $f$ above the boundary $\Cal X_{\eta}$ (cf. 2.5.2).
Then necessarily $r+m-1\ge 0$, and $g_y=(r+m-1)(p-1)/2$.
\par
2)\ $\Cal Y_k$ has $p$-{\bf branches} at $y$. Then the cover $f$ has an \'etale split reduction of type $(0,0)$
on the boundary, i.e. the induced torsor above $\Spf R[[T]]\{T^{-1}\}$ is trivial, in which case $g_y=(r-2)(p-1)/2$.}

\pop {.5}
\par
As an immediate consequence of 3.3.1 one can immediately see whether the point $y$ is smooth
or not. More precisely we have the following:

\pop {.5}
\par
\noindent
{\bf \gr 3.3.2. Corollary.}\rm\ {\sl We use the same notation as in 3.3.1. Then $y$ is a smooth point which is equivalent to $g_y=0$ if and
only if $r=1-m$ which implies that $m\le 1$. In particular if $f$ has a degeneration of type $(n,m)$ on the boundary with $n>0$ and $m>0$ then this happend only if 
$r=0$ and $m=1$.}

\pop {.5}
\par
Next we will give examples of Galois covers of degree $p$ above the formal germ of a smooth point which cover all the
possibilities for the genus and the degeneration type on the boundary. Both in 3.3.3 and 3.3.4 we use the same notations as in 3.3.1. 
We first begin with examples with genus $0$.

\pop {.5}
\par
\noindent
{\bf \gr 3.3.3. Examples.}\rm\ The following are examples given by explicit equations of the different cases, depending on the
possible degeneration type over the boundary, of Galois covers $f:\Cal Y\to \Cal X$ of degree $p$ above $\Cal X=\Spf R[[T]]$ 
and where $\bold {g_y=0}$ (here $y$ denotes the closed point of $\Cal Y$).
\par
1)\ For $m>0$ an integer prime to $p$ consider the cover given generically by the equation $X^p-X=T^{-m}$. Here $r=m+1$ and this cover 
has a reduction of type $(0,-m)$ on the boundary.
\par
2)\ For $\tilde m:=-m$ a negative integer prime to $p$ and a positive integer $n$ consider the cover given generically by the equation 
$X^p-\pi^{n(p-1)}X=T^{\tilde m}$. Here $r=m+1$ and this cover has a reduction of type $(n,\tilde m)$ on the boundary.
\par
3)\ For a positive integer $n$ consider the cover given generically by the equation $X^p-\pi^{n(p-1)}X=T$. Here $r=0$ and this cover has a 
reduction of type $(n,1)$ on the boundary.

\pop {.5}
\par
Next we give examples of Galois covers of degree $p$ above formal germs of smooth points which lead to a singularity with
positive genus.

\pop {.5}
\par
\noindent
{\bf \gr 3.3.4. Examples.}\rm\ The following are examples given by explicit equations of the different cases depending on the
possible reduction type of Galois covers $f:\Cal Y\to \Cal X$ of degree $p$ above $\Cal X=\Spf R[[T]]$ and where $\bold {g_y>0}$.

\par
1 ) Let $m>0$ and $m'>m$ be integers prime to $p$ and consider the cover given generically 
by the equation $X^p-X=\pi/T^{m'}+1/T^m$. This cover has a degeneration of type $(0,-m)$ on the boundary, the point 
$y$ above $x$ is singular, and its genus equals $(m'-m)(p-1)/2$.

\par
2 ) Let $m$, $m'$, and $n$ be positive integers with $m$ and $m'$ prime to $p$ and consider the cover given generically by the equation 
$X^p-X=T^m/\pi^{pn}+\pi/T^{m'}$. This cover has a degeneration of type $(n,m)$ on the boundary, the point $y$ above $x$ is singular, 
and its genus equals $(m'+m)(p-1)/2$.

\par
3 ) Let $m$, $m'$ and $n$ be positive integers such that $m$ and $m'$ are prime to $p$ and $m'>m$. Consider the cover given generically 
by the equation $X^p-X=T^{-m}\pi^{-pn}+\pi/T^{m'}$. This cover has a degeneration of type $(n,-m)$ on the boundary, the point $y$ 
above $x$ is singular, and its genus equals $(m'-m)(p-1)/2$.

\pop {.5}
\par
Next we examine the case of Galois covers of degree $p$ above formal germs at {\bf double points}.

\pop {.5}
\par
\noindent
{\bf \gr 3.3.5. Proposition.}\rm\ {\sl  Let $\Cal X:=\Spf R[[S,T]]/(ST-\pi ^e)$ be the formal germ of an $R$-curve at an ordinary 
double point $x$ of thickness $e$. Let $\Cal X_1:=\Spf R[[S]]\{S^{-1}\}$ and $\Cal X_2:=\Spf R[[T]]\{T^{-1}\}$ be the boundaries of 
$\Cal X$. Let $f:\Cal Y\to \Cal X$ be a Galois cover with group $\Bbb Z/p\Bbb Z$ and with $\Cal Y$ local. Assume that the special fibre of 
$\Cal Y$ is reduced. We assume that $\Cal Y_k$ has {\bf two branches} at the point $y$. Let $d_{\eta}:=r(p-1)$ be the degree of the 
divisor of ramification in the morphism $f:Y_K\to X_K$. Let $(n_i,m_i)$ be the degeneration type on the boundaries of $\Cal X$ 
for $i=1,2$. Then necessarily $r+m_1+m_2\ge 0$ and $g_y=(r+m_1+m_2)(p-1)/2$}.

\pop {.5}
\par
\noindent
{\bf \gr 3.3.6. Proposition.}\rm\ {\sl We use the same notations as in Proposition 3.3.5. We consider the remaining cases:
\par
1)
{\bf $\Cal Y_k$ has $p+1$ branches at} $y$ in which case we can assume that $\Cal Y$ is completely split above $\Cal X_1$.
Let $(n_2,m_2)$ be the degeneration type on the second boundary $\Cal X_2$ of $\Cal X$. Then necessarily $r+m_2-1\ge 0$ and
$g_y=(r+m_2-1)(p-1)/2$.
\par
2) {\bf $\Cal Y_k$ has $2p$ branches at $y$} in which case $\Cal Y$ is completely split above the two boundaries of $\Cal X$ and
$g_y=(r-2)(p-2)/2$.}

\pop {.5}
\par
With the same notations as in proposition 3.3.5 and as an immediate consequence one can recognise whether the point $y$ is a double point 
or not. More precisely we have the following:

\pop {.5}
\par
\noindent
{\bf \gr 3.3.7. Corollary.}\rm\ {\sl We use the same notations as in 3.3.5. Then $y$ is an ordinary double point which is equivalent to 
$g_y=0$ if and only if $x$ is an ordinary double point of thickness divisible by $p$ and $r=m_1+m_2$. Moreover if $r=0$ then $g_y=0$ 
is equivalent to $m_1+m_2=0$}.

\pop {.5}
\par
Next we give examples of Galois covers of degree $p$ above the formal germ of a double point which lead to singularities
with genus $0$, i.e. double points, and such that $r=0$. These examples will be used in V in order to realise the ``degeneration data''
corresponding to Galois covers of degree $p$ in equal characteristic $p>0$.

\pop {.5}
\par
\noindent
{\bf \gr 3.3.8. Examples.}\rm \ The following are examples given by explicit equations of the different cases, depending on the possible
degeneration type on the boundaries, of Galois covers $f:\Cal Y\to \Cal X$ of degree $p$ above $\Cal X=\Spf R[[S,T]]/(ST-\pi ^e)$ 
with $\bold {r=0}$ and where $\bold {g_y=0}$ for a suitable choice of $e$. Note
that $e=pt$ must be divisible by $p$. In all the following examples we have $r=0$.

\par
1 )\ $p$-{\bf Purity}: if $f$ as above has an \'etale reduction type on the boundaries and $r=0$ then $f$ is necessarily \'etale and hence is
completely split since $\Cal X$ is strictly henselian.

\par
2 )\ Consider the cover given generically by the equation $X^p-X=1/T^m=S^m/\pi^{mpt}$ where $m$ is a positive integer prime to $p$ which 
leads to a reduction on the boundaries of type $(0,-m)$ and $(mt,m)$.

\par
3 )\ Let $n$ and $m$ be positive integers such that $m$ is prime to $p$ and $n-tm>0$. Consider the cover given generically by the equation  
$X^p-X=T^m/\pi^{np}=S^{-m}/\pi^{p(n-tm)}$ which leads to a reduction on the boundaries of type $(n,m)$ and $(n-tm,-m)$.

\pop {.5}
\par
In fact one can describe Galois covers of degree $p$ above formal germs of double points (in equal characteristic $p$) which are \'etale 
above the generic fibre and with genus $0$. Namely they are all of the form given in the above examples 3.3.8. In particular these covers are 
uniquely determined up to isomorphism by their degeneration type on the boundaries. More precisely we have the following:

\pop {.5}
\par
\noindent
{\bf \gr 3.3.9. Proposition.}\rm\ {\sl Let $\Cal X$ be the formal germ of an $R$-curve at an ordinary double point $x$. 
Let $f:\Cal Y\to \Cal X$ be a Galois cover of degree $p$ with $\Cal {Y}_k$ reduced and local and with $f_K:\Cal Y_K\to \Cal X_K$ {\bf \'etale}. Let ${\Cal X_i}$ 
for $i=1,2$ be the boundaries of $\Cal X$. Let $f_i:\Cal Y_i\to \Cal X_i$ be the torsors induced by $f$ above $\Cal X_i$ and let $\delta _i$ 
be the corresponding degree of the different (cf. 2.5.1). Let $y$ be the closed point of $\Cal Y$ and assume that $\bold {g_y=0}$. 
Then there exists an isomorphism $\Cal X\simeq \Spf R[[S,T]]/(ST-\pi^{tp})$ such that if say $\Cal X_2$ is the boundary corresponding 
to the prime ideal $(\pi,S)$ one of the following holds:
\par
{\bf a )}\ The cover $f$ is generically given by the equation $X^p-X=1/T^m=S^m/\pi^{mpt}$ where $m$ is a positive integer prime to $p$. This cover 
leads to a reduction on the boundaries of $\Cal X$ of type $(0,-m)$ and $(mt,m)$. Here $t>0$ can be any integer. In this case $\delta _1=0$
and $\delta_2=mt(p-1)$.
\par
{\bf b )}\ The cover $f$ is generically given by an equation $X^p-X=T^m/\pi^{np}=1/\pi^{p(n-tm)}S^m$ where $m>0$ is an integer prime to $p$ 
and $n>0$ such that $n-tm>0$. This cover leads to a reduction on the boundaries of $\Cal X$ of type $(n,m)$ and $(n-tm,-m)$. 
In this case $\delta_1=n(p-1)$ and $\delta_2=(n-tm)(p-1)$.}

\pop {.5}
\par
\noindent
{\bf \gr Proof.}\rm \ The proof is similar to the proof of 4.2.5 in [S-1] in the inequal characteristic case. 
\pop {.5}
\par
\noindent
{\bf \gr 3.3.10. variation of the different.}\rm\ The following result, which is a direct consequence of Proposition 3.3.9,
describes the variation of the degree of the different from one boundary to another in a cover $f:\Cal Y\to \Cal X$ between formal germs at double points.

\pop {.5}
\par
\noindent
{\bf \gr 3.3.11. Proposition.}\ \rm {\sl Let $\Cal X$ be the formal germ of an $R$-curve at an ordinary double point $x$. Let $f:\Cal Y\to \Cal X$
be a Galois cover of degree $p$ with $\Cal {Y}_k$ reduced and local and with $f_K:\Cal Y_K\to \Cal X_K$ {\bf \'etale}. Let $y$ be the closed point 
of $\Cal Y$. Assume that $\bold {g_y=0}$ which implies necessarily that the thickness $e=pt$ of the double point $x$ is divisible by $p$. 
For each integer $0<t'<t$ let $\Cal X_{t'}\to \Cal X$ be the blow-up of $\Cal X$ at the ideal $(\pi ^{pt'},T)$. The special fibre of 
$\Cal X_{t'}$ consists of a projective line $P_{t'}$ which meets two germs of double points $x$ and $x'$. Let $\eta$ be the generic point of  
$P_{t'}$ and let $v_{\eta}$ be the corresponding discrete valuation of the function field of $\Cal X$. Let $f_{t'}:\Cal Y_{t'}\to \Cal X_{t'}$ be 
the pull back of $f$ which is a Galois cover of degree $p$ and let $\delta (t')$ be the degree of the different induced by this cover above  
$v_{\eta}$ (cf. 2.5.1). Also denote by ${\Cal X_i}$ for $i=1,2$ the boundaries of $\Cal X$. Let $f_i:\Cal Y_i\to \Cal X_i$ be the
torsors induced by $f$ above $\Cal X_i$. let $(n_i,m_i)$ be their degeneration type and let $\delta _i$ be the corresponding degree of
the different. Say $\delta _1=\delta (0)$ $\delta _2=\delta (t)$ and $\delta (0)\le \delta (t)$. We have $m:=-m_1=m_2$ say is positive. 
Then the following holds: for $0\le t_1\le t_2\le t$ we have $\delta (t_2)=\delta(t_1)+m(p-1)(t_2-t_1)$
and $\delta (t')$ is an increasing function of $t'$}.

\pop {1}
\par
\noindent
{\bf \gr IV. Computation of vanishing cycles and examples for $p^2$-cyclic covers.}
\pop {.5}
\par
\noindent
{\bf \gr 4.1.}\rm\ In this section we use the same notations and hypothesis as in III and we compute the dimensions of the spaces of 
vanishing cycles arising form a cyclic $p^2$-cover above the formal germ of an $R$-curve. Further we provide examples for such covers 
above formal germs of semi-stable curves.

\pop {1}
\par
\noindent
{\bf \gr 4.1.1. Proposition.}\rm\ {\sl Let $\Cal X:=\Spf {\hat {\Cal O}_x}$ be the formal germ of an $R$-curve at a closed point $x$ 
with $\Cal X_k$ reduced. Let $\tilde f:\Cal Y\to \Cal X$ be a Galois cover with group $\Bbb Z/p^2\Bbb Z$ with $\Cal Y$ local and
$\Cal Y_k$ reduced. Let $\{\wp_i\}_{i\in I}$ be the minimal prime ideals of $\hat {{\Cal O}}_x$ which contain $\pi$ and let
$\Cal X_{i}:=\Spf \hat \Cal O_{\wp _i}$ be the formal completion of the localisation of $\Cal X$ at $\wp _i$. For each $i\in I$ the above
cover $\tilde f$ induces a $p^2$-cyclic cover $\tilde f_i:\Cal Y_i\to \Cal X_{i}$ above the boundary $\Cal X_i$ (cf. 2.5.1). Let 
$\{(n_{i,1},m_{i,1}),(n_{i,2},m_{i,2})\}$ (resp. $\{(n_i,m_i)\}$) be the degeneration type of $\tilde f_i$ (cf. 2.6.2) 
if $\Cal {Y}_{i,k}$ is irreducible (resp. if $\Cal Y_{i,k}$ has $p$-components cf. 2.5.2). Let $y$ be the closed point of $\Cal Y$. Then 
one has the following {\bf ``local Riemann-Hurwitz formula''}:
$$2g_y-2=p^2(2g_x-2)+d_{\eta}-d_s$$
Where $d_{\eta}=pd_{\eta,1}+d_{\eta,2}$ is the degree of the divisor of ramification in the morphism $\tilde f_K:\Cal Y_K\to \Cal X_K$ 
induced by $\tilde f$ where $\Cal X_K:=\Spec ({\hat {\Cal O}_x} \otimes _R K)$ and $\Cal Y_K:=\Spec ({\hat {\Cal O}_{\Cal Y,y}}\otimes _R K)$ 
and $d_s:=pd_{s,1}+d_{s,2}$ where $d_{s,1}=\sum _{i\in I_{\et}}(-m_{i,1}-1)(p-1)+\sum _{i\in I_{\rad}}(-m_{i,1}-1)(p-1)$ 
(resp. $d_{s,2}=\sum _{i\in I_{\et}}(-m_{i,2}-1)(p-1)+\sum _{i\in \tilde I_{\et}}(-m_{i}-1)(p-1)
+\sum _{i\in I_{\rad}}(-m_{i,2}-1)(p-1)+\sum _{i\in \tilde I_{\rad}}(-m_{i}-1)(p-1)$)  where 
$I^{\rad}$ is the subset of $I$ consisting of those $i$ for which $n_{i,1}\neq 0$ (resp.  $n_{i,2}\neq 0$) and $I^{\et}$ is the subset of 
$I$ consisting of those $i$ for which $n_{i,1}=0$ (resp. $n_{i,2}=0$). Here $\tilde I$ denotes the subset
of $I$ consisting of those $i$ for which $\Cal Y_i$ has $p$-components.}

\pop {.5}
\par
\noindent
{\bf \gr Proof.}\rm\ Follows directly from 3.2.3

\pop {.5}
\par
\noindent
{\bf \gr 4.2. $p^2$-Cyclic covers above germs of semi-stable curves.}\rm
\par
In what follows and as a consequence of 4.1.1 we will deduce some results 
in the case of a $p^2$-cyclic cover $f:\Cal Y\to \Cal X$ where $\Cal X$ is the formal germ of a semi-stable $R$-curve at a closed point. 
We start with the case of a Galois cover of degree $p^2$ above a germ of a smooth point.
\pop {.5}
\par
\noindent
{\bf \gr 4.2.1. Proposition.}\rm\ {\sl Let $\Cal X:=\Spf R[[T]]$ be the germ of a formal $R$-curve at a {\bf smooth} point $x$ and let
$\Cal X_{\eta}:=\Spf R[[T]]\{T^{-1}\}$ be the boundary of $\Cal X$. Let $f:\Cal Y\to \Cal X$ be a Galois cover of degree $p^2$
with $\Cal Y$ local. Assume that the special fibre of $\Cal Y$ is reduced. Let $y$ be the unique closed point of $\Cal Y_k$.
Let $d_{\eta}$ be the degree of the divisor of ramification in the morphism $f:\Cal Y_K\to \Cal X_K$. Then
$d_{\eta}=r(p-1)$ is divisible by $p-1$. The following cases occur:
\par
1)\ $\Cal Y_k$ is {\bf unibranche} at $y$. Let $\{(n_1,m_1),(n_2,m_2)\}$ be the degeneration type of $f$ above the boundary $\Cal X_{\eta}$ 
(cf. 2.5.1). Then necessarily $r+pm_1+m_2-p-1\ge 0$ and $g_y=(r+pm_1+m_2-p-1)(p-1)/2$.
\par
2)\ $\Cal Y_k$ has $p$-{\bf branches} at $y$. Let $\{(n,m)\}$ be the degeneration type of $f$ above the boundary $\Cal X_{\eta}$.
Then necessarily $r+pm-p-2\ge 0$ and $g_y=(r+pm-p-2)(p-1)/2$.

\par
3)\ $\Cal Y_k$ has $p^2$-{\bf branches} at $y$. Then the cover $f$ has an \'etale completely split reduction 
on the boundary, i.e. the induced torsor above $\Spf R[[T]]\{T^{-1}\}$ is trivial, in which case necessarily $r-2p-2\ge 0$ and 
$g_y=(r-2p-2)(p-1)/2$.}

\pop {.5}
\par
As an immediate consequence of 4.2.1 one can immediately see whether the point $y$ is smooth or not. More precisely we have the following:

\pop {.5}
\par
\noindent
{\bf \gr 4.2.2. Corollary.}\rm\ {\sl We use the same notation as in 4.2.1. Then $y$ is a smooth point which is equivalent to $g_y=0$ if and
only if we are in case 1) and $r+pm_1+m_2=p+1$. In particular if $f$ has a degeneration on the boundaries of type
$\{(n_1,m_1),(n_2,m_2)\}$ with $m_1$ and $m_2$ positive then this happend only if $r=0$ and $m_1=m_2=1$}.

\pop {.5}
\par
Next we will give examples of cyclic covers of degree $p^2$ above the formal germ of a smooth point which cover all the
possibilities for the genus and the degeneration type on the boundary. Both in 4.2.3 and 4.2.4 we use the same notations as in 4.2.1. 
We first begin with examples with genus $0$.

\pop {.5}
\par
\noindent
{\bf \gr 4.2.3. Examples.}\rm\ The following are examples given by explicit equations of the different cases, depending on the possible 
degeneration type over the boundary, of cyclic covers $\tilde f:\Cal Y\to \Cal X$ of degree $p^2$ above $\Cal X=\Spf R[[T]]$ and where 
$\bold {g_y=0}$ (here $y$ denotes the closed point of $\Cal Y$).
\par
1)\ Let $m_1$ and $m_2$ be two positive integers both prime to $p$ and consider the cover which is generically given by the equations: 
$(X_1^p,X_2^p)-(X_1,X_2)=(1/T^{m_1},f(T))$ where $f(T)=\sum _{i\ge -m_2}c_iT^i\in R[[T]][T^{-1}]$ and $c_{m_2}$ is a unit in $R$.
Assume that $m_2\ge pm_1$. Then in this case $r=m_1+1+pm_2+p$ and this cover has a degeneration on the boundary of type 
$\{(0,-m_1),(0,\tilde m_2:=-pm_2+m_1(p-1))\}$. If $m_2<pm_1$ then $r=m_1+1+p^2m_1+p$ and this cover has a degeneration on the boundary of type 
$\{(0,-m_1),(0,-m_1(p(p-1)+1))\}$.

\par
2)\ Let $m_1$ and $m_2$ be two positive integers both prime to $p$ such that  $m_2\ge pm_1$ and consider the cover given generically 
by the equations: $(X_1^p,X_2^p)-(X_1,X_2)=$
\newline
$(1/T^{m_1},f(T)/\pi^{pn})$ where $n$ is a positive integer and $f(T)=\sum _{i\ge -m_2}c_iT^i\in R[[T]][T^{-1}]$ where 
$c_{-m_2}$ is a unit in $R$. Then in this case $r=m_1+1+pm_2+p$ and this cover has a degeneration on the boundary of type 
$\{(0,-m_1),(n,\tilde m_2:=-pm_2+m_1(p-1))\}$. If $m_2<pm_1$ then $\Cal Y$ is not smooth.

\par
3)\ Consider the cover given generically by the equations: $(X_1^p,X_2^p)-(X_1,X_2)=(T/\pi^{n_1p},f(T)/\pi^{n_2p})$ where $f(T)\in R[[T]][T^{-1}]$ 
has conductor $m_2=1$ (cf. 2.6.1 for the definition of the conductor). Assume that $n_2p\ge n_1(p(p-1)+1)$. In this case $r=0$ and $m_1=m_2=1$ 
and this cover has a degeneration on the boundary of type $\{(n_1,1),(n_2,1)\}$.

\par
4)\ Let $m_1$ and $m_2$ be two positive integers both prime to $p$. Consider the cover given generically by the equations 
$(X_1^p,X_2^p)-(X_1,X_2)=(T^{-m_1}/\pi^{n_1p},f(T)/\pi^{pn_2})$ where $f(T)\in R[[T]][T^{-1}]$ has conductor $-m_2$ (cf. 2.6.1). 
First assume that $pn_2':=n_2p-n_1(p-1)\ge n_1(p(p-1)+1)$ and  $m_2\ge pm_1$. Then in this case $r=m_1+1+pm_2+p$ and this cover has 
a degeneration on the boundary of type $\{(n_1,-m_1),(n_2',\tilde m_2:=-m_2p+m_1(p-1))\}$. Further in this case if $m_2<pm_1$ then $\Cal Y$ is not smooth.
Second assume that $n_2p-n_1(p-1)<pn_2':=n_1(p(p-1)+1)$ and $m_2\le pm_1$ then in this case $r=m_1+1+m_1p^2+p$ and this cover has a degeneration on the boundary of type 
$\{(n_1,-m_1),((n_2',\tilde m_2:=-m_1(p(p-1)+1))\}$. In this later case if $m_2>pm_1$ then $\Cal Y$ is not smooth.

\pop {.5}
\par
\noindent
{\bf \gr 4.2.4. Corollary.}\rm\ {\sl Let $\Cal X:=\Spf R[[T]]$ be the germ of a formal $R$-curve at a {\bf smooth} point $x$ and let
$\Cal X_{\eta}:=\Spf R[[T]]\{T^{-1}\}$ be the boundary of $\Cal X$. Let $\tilde f:\Cal Y\to \Cal X$ be a Galois cover of degree $p^2$
with $\Cal Y$ local and which ramifies above the generic fibre $\Cal X_K$. Let $\{(n_1,m_1),(\tilde n_2,\tilde m_2)\}$ be the degeneration type 
of the cover $\tilde f$ above the formal boundary
$\Cal X_{\eta}$. Assume that $\Cal Y$ is smooth. If $n_1=0$ and $n_2=0$ (resp. $n_1=0$ and $n_2\neq 0$) then necessarily $m_1\le -1$ and 
$\tilde m_2=\inf (m_1(p(p-1)+1),m_2p+m_1(p-1))$
for some negative integer $m_2$ prime to $p$ (resp. $\tilde m_2=m_2p+m_1(p-1)$ for some integer $m_2$ prime to $p$ with $-m_2\ge pm_1$). 
If $n_1\neq 0$ then necessarily $m_1\le -1$ and $\tilde m_2\le -1$ and either $\tilde n_2>n_1(p(p-1)+1)/p$ in which case 
$\tilde m_2=pm_2-m_1(p-1)$ for some integer $m_2$ prime to $p$ with $-m_2\ge -pm_1$ or $\tilde n_2=n_1(p(p-1)+1)/p$
in which case $\tilde m_2=m_1(p(p-1)+1)$ or $\tilde m_2=m_2p-m_1(p-1))$ for some integer $m_2$ prime to $p$ such that $m_2<pm_1$.}

\pop {.5}
\par
\noindent
{\bf \gr Proof.}\rm\ Follows easily from 2.6.1 and the examples in 4.2.3 which cover all the possibilities for $f$ as above.

\pop {.5}
\par
\noindent
{\bf \gr 4.2.5. Definition.}\rm\ An admissible pair $\{(n_1,m_1),(\tilde n_2,\tilde m_2)\}$ as defined in 2.6.4 is said to satisfy the 
condition (*) if it satisfies the numerical contraints in 4.2.4 above.

\pop {.5}
\par
Next we give examples of $p^2$-cyclic covers above formal germs of smooth points which lead to singularities with positive genus.

\pop {.5}
\par
\noindent
{\bf \gr 4.2.6. Examples.}\rm\ The following are examples given by explicit equations of the different cases, depending on the
possible reduction type, of cyclic covers $\tilde f:\Cal Y\to \Cal X$ of degree $p^2$ above $\Cal X=\Spf R[[T]]$ and where $\bold {g_y>0}$.

\par
1 )\ Let $m_1$, $m'_1$, $\tilde m_2$ and $\tilde m_2'$ be positive integers all prime to $p$ such that $m_1<m_1'$ and $\tilde m_2<\tilde m_2'$.
Consider the cover given generically by the equations  $(X_1^p,X_2^p)-(X_1,X_2)=(\pi/T^{m_1'}+1/T^{m_1},f(T))$ where $f(T)=\sum _{i\ge -\tilde m_2'}c_iT^i\in 
R[[T]][T^{-1}]$ with $c_{-\tilde m_2'}\in \pi R$ and $\bar f(t)=\sum _{i\ge -\tilde m_2}\bar c_it^i\in  k[[T]]$ with $\bar c_{-\tilde m_2}\neq 0$
(i.e. $\tilde m_2$ is the conductor of $f_2$). Assume that $\tilde m'_2>pm_1'$. In this case 
$r=m_1'+pm_2'+p+1$  and this cover has a degeneration on the boundary of type $\{(0,-m_1),(0,m_2:=-p\tilde m_2+m_1(p-1))\}$. Moreover
$g_y=(p\tilde m_2'-p\tilde m_2+m_1'-m_1)(p-1)/2$.

\par
2)\  Let $m_1$, $m'_1$, $\tilde m_2$, and $\tilde m_2'$ be positive integers all prime to $p$ such that $m_1<m_1'$ and $\tilde m_2<\tilde m_2'$.
Consider the cover given generically by the equations $(X_1^p,X_2^p)-(X_1,X_2)=(\pi/T^{m_1'}+1/T^{m_1},f(T)/\pi^{pn})$ 
where $n$ is a positive integer and 
$f(T)$ is as in 1). Assume that $\tilde m'_2>pm_1'$. In this case $r=m_1'+p\tilde m_2'+p+1$ and
this cover has a degeneration on the boundary of type $\{(0,-m_1),(n,m_2:=-p\tilde m_2+m_1(p-1))\}$. Moreover 
$g_y=(p\tilde m_2'-p\tilde m_2+m_1'-m_1)(p-1)/2$.

\par
3)\  Let $m_1$, $m'_1$, $\tilde m_2$, and $\tilde m_2'$ be positive integers all prime to $p$ such that $m_1<m_1'$ and $\tilde m_2<\tilde m_2'$.
Consider the cover given by the equations $(X_1^p,X_2^p)-(X_1,X_2)=(\pi/T^{m_1'}+T^{-m_1}\pi ^{n_1p},f(T)/\pi^{pn_2})$ where $n_1$ and $n_2$ are 
positive integers and  $f(T)$ is as in 1). 
Assume that $\tilde m'_2>pm_1'$ in which case we have $r=m_1'+p\tilde m_2'+p+1$. We distinguish the following cases:
\par
3.a)\ Assume that $pn_2':=n_2p-n_1(p-1)\ge n_1(p(p-1)+1)$. Then this cover has a degeneration on the boundary of type 
$\{(n_1,-m_1),(n_2',m_2:=-\tilde m_2p+m_1(p-1))\}$ and $g_y=(p\tilde m_2'-p\tilde m_2+m_1'-m_1)(p-1)/2$.

\par
3.b)\ Assume that $n_2p-n_1(p-1)<pn_2':=n_1(p(p-1)+1)$ then this cover has a degeneration on the boundary of type 
$\{(n_1,-m_1),(n_2',m_2:=-m_1(p(p-1)+1))\}$ and $g_y=(p\tilde m_2'-p^2m_1+m_1'-m_1)(p-1)/2$.

\par
4)\  Let $m_1$, $m'_1$, $\tilde m_2$, and $\tilde m_2'$ be positive integers all prime to $p$ such that $m_1<m_1'$ and $\tilde m_2<\tilde m_2'$.
Consider the cover given generically by the equations $(X_1^p,X_2^p)-(X_1,X_2)=(\pi/T^{m_1'}+T^{-m_1}\pi ^{n_1p},f(T)/\pi^{pn_2})$ where 
$n_1$ and $n_2$ are positive integer and  $f(T)$ is as in 1). 
Assume that $\tilde m'_2<pm_1'$ in which case we have $r=m_1'+m_1'p^2+p+1$. We distinguish the following cases:
\par
4.a)\ Assume that $pn_2':=n_2p-n_1(p-1)\ge n_1(p(p-1)+1)$. Then this cover has a degeneration on the boundary of type 
$\{(n_1,-m_1),(n_2',m_2:=-\tilde m_2p+m_1(p-1))\}$ and $g_y=(p^2m_1'-\tilde m_2p+m_1'-m_1)(p-1)/2$.

\par
4.b)\ Assume that $n_2p-n_1(p-1)<pn_2':=n_1(p(p-1)+1)$ then this cover has a degeneration on the boundary of type 
$\{(n_1,-m_1),(n_2',m_2:=-m_1(p(p-1)+1))\}$ and $g_y=(p^2m_1'-m_1p^2+m_1'-m_1)(p-1)/2$.

\pop {.5}
\par
Next we examine the case of cyclic covers of degree $p^2$ above formal germs at {\bf double points}.

\pop {.5}
\par
\noindent
{\bf \gr 4.2.7. Proposition.}\rm\ {\sl  Let $\Cal X:=\Spf R[[S,T]]/(ST-\pi ^e)$ be the formal germ of an $R$-curve at an ordinary 
double point $x$ of thickness $e$ and let $\Cal X_1:=\Spf R[[S]]\{S^{-1}\}$ and $\Cal X_2:=\Spf R[[T]]\{T^{-1}\}$ be the boundaries of 
$\Cal X$. Let $f:\Cal Y\to \Cal X$ be a Galois cover with group $\Bbb Z/p^2\Bbb Z$ and with $\Cal Y$ local. Assume that the special fibre of 
$\Cal Y$ is reduced. We assume that $\Cal Y_k$ has {\bf two branches} at the point $y$. Let $d_{\eta}:=r(p-1)$ be the degree of the 
divisor of ramification in the morphism $f:Y_K\to X_K$. Let $\{(n_{i,1},m_{i,1}),(n_{i,2},m_{i,2})\}$ be the degeneration type on the 
boundaries of $\Cal X$ for $i=1,2$. Then necessarily $r+pm_{1,1}+pm_{2,1}+m_{1,2}+m_{2,2}\ge 0$ and 
$g_y=(r+pm_{1,1}+pm_{2,1}+m_{1,2}+m_{2,2})(p-1)/2$}.

\pop {.5}
\par
Next we examine the remaining cases.

\pop {.5}
\par
\noindent
{\bf \gr 4.2.8. Proposition.}\rm\ {\sl We use the same notation as in Proposition 4.2.7. We consider the remaining cases:
\par
1)\ {\bf $\Cal Y_k$ has $p+1$ branches at} $y$ in which case we can assume that the degeneration type of $f$ above 
the boundaries $\Cal X_1$ is $\{(n_{1},m_{1})\}$. Then necessarily $r+pm_{2,1}+p+m_{1}+m_{2,2}\ge 0$ and 
$g_y=(r+pm_{2,1}+p+m_{1}+m_{2,2})(p-1)/2$.

\par
2)\ {\bf $\Cal Y_k$ has $2p$ branches at $y$} in which case we can assume that $n_{1,1}=n_{2,1}=m_{1,1}=m_{2,1}=0$. 
Then necessarily $r+m_{2,2}+m_{2,1}-2p\ge 0$ and $g_y=(r+m_{2,2}+m_{2,1}-2p)(p-1)/2$.

\par
3)\ {\bf $\Cal Y_k$ has $p^2+1$ branches at $y$} and we can assume that $f$ is completely split above the boundary $\Cal X_1$. Then 
necessarily $r+pm_{2,1}+m_{2,2}-p-1\ge 0$ and $g_y=(r+pm_{2,1}+m_{2,2}-p-1)(p-1)/2$.
 
\par
4)\ {\bf $\Cal Y_k$ has $p^2+p$ branches at $y$} and we can assume that $f$ is completely split above the boundary $\Cal X_1$.
Then necessarily $r+m_{2,2}-2p-1=0$ and $g_y=(r+m_{2,2}-2p-1)(p-1)/2$.
\par
5)\ {\bf $\Cal Y_k$ has $2p^2$ branches at $y$} and $f$ is completely split on both boundaries of $\Cal X$. In this case 
$g_y=(r-2p-2)(p-1)/2$.}

\pop {.5}
\par
With the same notation as in proposition 4.2.7 and as an immediate consequence one can recognise whether the point $y$ is a double point 
or not. More precisely we have the following:

\pop {.5}
\par
\noindent
{\bf \gr 4.2.9. Corollary.}\rm\ {\sl We use the same notations as in 4.2.5. Then $y$ is an ordinary double point which is equivalent to 
$g_y=0$ if and only if $x$ is an ordinary double point of thickness divisible by $p^2$ and $r+pm_{1,1}+pm_{2,1}+m_{1,2}+m_{2,2}=0$. Moreover if 
$r=0$ then $g_y=0$ is equivalent to $m_{1,1}+m_{2,1}=0$ and $m_{1,2}+m_{2,2}=0$}.

\pop {.5}
\par
Next we give examples of $p^2$-cyclic covers above the formal germ of a double point which lead to singularities
with genus $0$, i.e. double points, and such that $r=0$. These examples will be used in VI in order to realise the ``degeneration data''
corresponding to cyclic covers of degree $p^2$ in equal characteristic $p>0$.

\pop {.5}
\par
\noindent
{\bf \gr 4.2.10. Examples.}\rm \ The following are examples given by explicit equations of the different cases, depending on the possible
degeneration type on the boundaries, of cyclic covers $\tilde f:\Cal Y\to \Cal X$ of degree $p^2$ above $\Cal X=\Spf R[[S,T]]/(ST-\pi ^e)$ 
with $\bold {r=0}$ and where $\bold {g_y=0}$ for a suitable choice of $e$. Note
that $e=p^2t$ must be divisible by $p^2$. In all the following examples we have $r=0$.

\par
1 )\ $p$-{\bf Purity}: if $\tilde f$ as above has an \'etale reduction type on the boundaries and $r=0$ then $f$ is necessarily \'etale and hence is
completely split since $\Cal X$ is strictly henselian.

\par
2 )\ Consider the cover given generically by the equation $(X_1^p,X_2^p)-(X_1,X_2)=(1/T^{m_1},f(T))$ where $m_1$ is a positive integer prime 
to $p$ and $f(T)\in R[[T]][T^{-1}]$ is such that its image $\bar f(t)$ modulo $\pi$ equals $\sum _{i\ge -m_2}c_it^i$ where $m_2\ge 0$ is prime to $p$.
This cover leads to a reduction on the boundaries of type $\{(0,-m_1),(0,\tilde m_2:=-pm_2+m_1(p-1))\}$ and $\{(tm_1p,m_1),(pm_2t-(p-1)m_1t
,pm_2-m_1(p-1))\}$ if $m_2> m_1p$, and of type $\{(0,-m_1),(0,\tilde m_2:=-m_1(p(p-1)+1))\}$ and $\{(tm_1p,m_1),
(tm_1(p(p-1)+1),m_1(p(p-1)+1))\}$ if $m_2\le m_1p$.

\par
3 )\ Consider the cover given generically by the equation  
$X_1^p-X_1=1/T^{m_1}$ and $X_2^p-X_2=f(T)/\pi ^{pn}+p^{-1}\sum _{k=1}^{p-1}{p \choose k}X_1^{pk}(-X_1)^{p-k}$
with $f(T)\in  R[[T]][T^{-1}]$ is such that its image $\bar f(t)=\sum _{i\ge m'}c_it^i$ modulo $\pi$ is not a $p$-power, and $n$ is a positive integer such that
$pn-m_2p^2t>0$. Let $m_2$ be the conductor of $f(T)$. If $n-m_2pt>m_1p^2t$ then $\Cal Y$ is smooth and 
this cover leads to a reduction on the boundaries of type $\{(0,-m_1),(n,pm_2+m_1(p-1))\}$ and $\{(m_1pt,m_1),(n-m_2pt-(p-1)m_1t,-pm_2-m_1(p-1))\}$.
Further in this case if $n-m_2pt\le m_1p^2t$ then $\Cal Y$ is not smooth.

\par
4)\ Consider the cover given generically by the equation
$X_1^p-\pi^{n_1(p-1)}X_1=T^{m_1}$ for some integer $m_1$ prime to $p$ and: 
$$X_2^p-X_2=f(T)/\pi^{pn_2}+p^{-1}\sum _{k=1}^{p-1}{p \choose k}X_1^{pk}(-X_1)^{p-k}/\pi^{n_1(p(p-k)+k)}$$ 
where $f(T)=\sum _{i\in I} a_iT^i\in  R[[T]][T^{-1}]$
is not a $p$-power and where $n_1$ (resp. $n_2$) is a positive integer such that $n_1-ptm_1>0$ (resp.  $n_2-ptm_2>0$). 
Here $m_2$ denotes the conductor of $f_2$. We distinguish the following cases:

\par
4-1)\  $pn_2-(p-1)n_1>n_1(p(p-1)+1)$. If $pn_2-m_2p^2t-(p-1)(n_1-m_1pt)>(n_1-m_1pt)(p(p-1)+1)$ which necessarily implies that $m_1p>m_2$ then $\Cal Y$ 
is smooth and this cover leads to a reduction on the boundaries of type
$\{(n_1,m_1),((n_2p-n_1(p-1))/p,m_2p-m_1(p-1))\}$ and $\{(n_1-m_1pt,-m_1),((n_2p-m_2p^2t-(p-1)(n_1-m_1pt))/p,-m_2p+m_1(p-1))\}$. Further in this case if
$pn_2-m_2p^2t-(p-1)(n_1-m_1pt)\le (n_1-m_1pt)(p(p-1)+1)$ then $\Cal Y$ is not smooth.

\par
4-2)\ $pn_2-(p-1)n_1<n_1(p(p-1)+1)$. If $pn_2-m_2p^2t-(p-1)(n_1-m_1pt)<(n_1-m_1pt)(p(p-1)+1)$ (resp. $pn_2-m_2p^2t-(p-1)(n_1-m_1pt)=(n_1-m_1pt)(p(p-1)+1)$)
which necessarily implies that $m_1p<m_2$ (resp. $m_1p>m_2$) then $\Cal Y$ is smooth and this cover leads to a reduction on the boundaries of type
$\{(n_1,m_1),(n_1(p(p-1)+1)/p,m_1(p(p-1)+1))\}$ and $\{(n_1-m_1pt,-m_1),((n_1-m_1pt)(p(p-1)+1)/p,-m_1(p(p-1)+1))\}$. Further in this case if
$pn_2-m_2p^2t-(p-1)(n_1-m_1pt)>(n_1-m_1pt)(p(p-1)+1)$ then $\Cal Y$ is not smooth.

\par
4-3)\ $pn_2-(p-1)n_1=n_1(p(p-1)+1)$. If $pn_2-m_2p^2t-(p-1)(n_1-m_1pt)>(n_1-m_1pt)(p(p-1)+1)$ which is equivalent to $m_2<m_1p$ 
(resp. $pn_2-m_2p^2t-(p-1)(n_1-m_1pt)<(n_1-m_1pt)(p(p-1)+1)$ which is equivalent to $m_2>m_1p$) then  
$\Cal Y$ is smooth and this cover leads to a reduction on the boundaries of type
$\{(n_1,m_1),(n_1(p(p-1)+1)/p,m_2p-m_1(p-1))\}$ and $\{(n_1-m_1pt,-m_1),((n_2p-m_2p^2t-(p-1)(n_1-m_1pt))/p,-m_2p+m_1(p-1))\}$
(resp. $\{(n_1,m_1),(n_1(p(p-1)+1)/p,m_1(p(p-1)+1))\}$ and $\{(n_1-m_1pt,-m_1),((n_1-m_1pt)(p(p-1)+1)/p,-m_1(p(p-1)+1))\}$).

\pop {.5}
\par
In fact one can describe cyclic covers of degree $p^2$ above formal germs of double points (in equal characteristic $p$) which are \'etale 
above the generic fibre and with genus $0$. Namely they are all of the form given in the above examples 4.2.10. More precisely we have the following:

\pop {.5}
\par
\noindent
{\bf \gr 4.2.11. Proposition.}\rm\ {\sl Let $\Cal X$ be the formal germ of an $R$-curve at an ordinary double point $x$. 
Let $\tilde f:\Cal Y\to \Cal X$ be a cyclic cover of degree $p^2$ with $\Cal {Y}_k$ reduced and local and with $f_K:\Cal Y_K\to \Cal X_K$ {\bf \'etale}. 
Let ${\Cal X_i}$ for $i=1,2$ be the boundaries of $\Cal X$. Let $y$ be the closed point of $\Cal Y$ and assume that $\bold {g_y=0}$. 
Then there exists an isomorphism $\Cal X\simeq \Spf R[[S,T]]/(ST-\pi^{tp^2})$ such that the following holds:

\par
{\bf a)}\ The cover $\tilde f$ is generically given by the equation:
$$(X_1^p,X_2^p)-(X_1,X_2)=(1/T^{m_1},f(T))$$ 
where $m_1$ is a positive integer prime to $p$ and $f(T)\in  R[[T]][T^{-1}]$ is such that its image $\bar f(t)$ modulo $\pi$ equals 
$\sum _{i\ge -m_2}c_it^i$, with $c_{-m_2}\neq 0$, 
where $m_2\ge 0$ is prime to $p$, which leads to a reduction on the boundaries of type $\{(0,-m_1),(0,\tilde m_2:=-pm_2+m_1(p-1))\}$ 
and $\{(tm_1p,m_1),(pm_2t-(p-1)m_1t,pm_2-m_1(p-1))\}$ 
if $m_2> m_1p$ and of type $\{(0,-m_1),(0,\tilde m_2:=-m_1(p(p-1)+1))\}$ and $\{(tm_1p,m_1),(tm_1(p(p-1)+1),m_1(p(p-1)+1))\}$ if $m_2\le m_1p$.

\par
{\bf b )}\ The cover $\tilde f$ is given generically by the equations: 
$$X_1^p-X_1=1/T^{m_1}$$ 
and: 
$$X_2^p-X_2=f(T)/\pi ^{pn}+p^{-1}\sum _{k=1}^{p-1}{p \choose k}X_1^{pk}(-X_1)^{p-k}$$
with $f(T)\in  R[[T]][T^{-1}]$ is such that its image $\bar f(t)=\sum _{i\ge m'}c_it^i$ modulo $\pi$ is not a $p$-power and $n$ is a 
positive integer such that
$n-m_2pt>m_1p^2t$. Here $m_2$ denotes the conductor of $f(T)$. this cover leads to a reduction on the boundaries of type 
$\{(0,-m_1),(n,pm_2+m_1(p-1))\}$ and $\{(m_1pt,m_1),(n-m_2pt-(p-1)m_1t,-pm_2-m_1(p-1))\}$.

\par
{\bf c)}\  The cover $\tilde f$ is given generically by the equations:
$$X_1^p-\pi^{n_1(p-1)}X_1=T^{m_1}$$ 
for some integer $m_1$ prime to $p$ and: 
$$X_2^p-X_2=f(T)/\pi^{pn_2}+p^{-1}\sum _{k=1}^{p-1}{p \choose k}X_1^{pk}(-X_1)^{p-k}/\pi^{n_1(p(p-k)+k)}$$ 
where $f(T)=\sum _{i\in I} a_iT^i\in  R[[T]][T^{-1}]$
is not a $p$-power and where $n_1$ (resp. $n_2$) is a positive integer such that $n_1-ptm_1\ge 0$ (resp.  $n_2-ptm_2\ge 0$). 
Here $m_2$ denotes the conductor of $f_2$ and the following cases occur:

\par
{\bf c-1)}\  $pn_2-(p-1)n_1>n_1(p(p-1)+1)$ and $pn_2-m_2p^2t-(p-1)(n_1-m_1pt)>(n_1-m_1pt)(p(p-1)+1)$ (which necessarily implies that $m_1p>m_2$), in which case 
this cover leads to a reduction on the boundaries of type $\{(n_1,m_1),((n_2p-n_1(p-1))/p,m_2p-m_1(p-1))\}$ and 
$\{(n_1-m_1pt,-m_1),((n_2p-m_2p^2t-(p-1)(n_1-m_1pt))/p,-m_2p+m_1(p-1))\}$.

\par
{\bf c-2)}\ $pn_2-(p-1)n_1<n_1(p(p-1)+1)$ and $pn_2-m_2p^2t-(p-1)(n_1-m_1pt)<(n_1-m_1pt)(p(p-1)+1)$ (resp. $pn_2-m_2p^2t-(p-1)(n_1-m_1pt)=(n_1-m_1pt)(p(p-1)+1)$)
which necessarily implies that $m_1p<m_2$ (resp. $m_1p>m_2$). This cover leads to a reduction on the boundaries of type
$\{(n_1,m_1),(n_1(p(p-1)+1)/p,m_1(p(p-1)+1))\}$ and $\{(n_1-m_1pt,-m_1),((n_1-m_1pt)(p(p-1)+1)/p,-m_1(p(p-1)+1))\}$.
\par
{\bf c-3)}\ $pn_2-(p-1)n_1=n_1(p(p-1)+1)$ and $pn_2-m_2p^2t-(p-1)(n_1-m_1pt)>(n_1-m_1pt)(p(p-1)+1)$, which is equivalent to $m_2<m_1p$, 
(resp. $pn_2-m_2p^2t-(p-1)(n_1-m_1pt)<(n_1-m_1pt)(p(p-1)+1)$ which is equivalent to $m_2>m_1p$). This cover leads to a reduction 
on the boundaries of type $\{(n_1,m_1),(n_1(p(p-1)+1)/p,m_2p-m_1(p-1))\}$ and $\{(n_1-m_1pt,-m_1),((n_2p-m_2p^2t-(p-1)(n_1-m_1pt))/p,-m_2p+m_1(p-1))\}$
(resp. $\{(n_1,m_1),(n_1(p(p-1)+1)/p,m_1(p(p-1)+1))\}$ and $\{(n_1-m_1pt,-m_1),((n_1-m_1pt)(p(p-1)+1)/p,-m_1(p(p-1)+1))\}$).

\pop {1}
\par
\noindent
{\bf \gr V. Semi-stable reduction of cyclic $p$-covers above formal germs of curves in equal characteristic $p>0$.}

\rm
\pop {.5}
\par
In all this paragraph we use the following notations: $R$ is a complete discrete valuation ring of equal characteristic 
$p$ with residue field $k$ which we assume to be {\bf algebraically closed} and fraction field $K:=\Fr R$. 
We denote by $\pi$ a uniformising parameter of $R$.

\pop {.5}
\par
\noindent
{\bf \gr 5.1.}\ Let ${\Cal X}:=\Spec \hat {\Cal O}_{X,x}$ be the 
formal germ of an $R$-curve $X$ at a closed point $x$ and let
$f:\Cal Y\to \Cal X$ be a Galois cover with 
group $G$ such that  $\Cal Y$ is normal and 
local. In this paper we are mainly concerned with the case where $G$ is cyclic of order $p$ or $p^2$. It follows then easily from the
theorem of semi-stable reduction for curves (cf. [De-Mu]), as well as the 
compactification process in 3.2, that after 
eventually a finite extension $R'$ of $R$ with
fractions field $K'$ the formal germ $\Cal Y$ has a semi-stable 
reduction. More precisely there exists a birational and proper morphism 
$\Tilde f :\Tilde {\Cal Y}\to \Cal Y'$ where  $\Cal Y'$ is the 
normalisation of $\Cal Y\times _R R'$ such that 
$\Tilde {\Cal Y}_{K'}\simeq \Cal Y'_{K'}$ and 
the following conditions hold:
\pop {.5}
\par
\noindent
(i) The special fibre $\Tilde {\Cal Y}_{k}:=\Tilde {\Cal Y}\times _{\Spec R'}\Spec k$ of $\Tilde {\Cal Y}$ is reduced.
\par
\noindent
(ii) $\Tilde {\Cal Y}_k$ has only ordinary double points as singularities.
\pop {.5}
\par
Moreover there exists such a semi-stable model $\Tilde f :\Tilde {\Cal Y}
\to \Cal Y'$ which is {\bf minimal} for the above properties. In 
particular the action of $G$ on $\Cal Y'$ extends to an action
on $\Tilde {\Cal Y}$. Let $\Tilde {\Cal X}$ be the quotient of 
$\Tilde {\Cal Y}$ by $G$ which is a semi-stable model of $\Cal X$. One has
the following commutative diagram:

$$
\CD
\Tilde {\Cal Y}     @> \Tilde f >> \Cal Y '\\
   @Vg VV                  @Vf' VV   \\
\Tilde {\Cal X} @>\Tilde g>>  {\Cal X'}
\endCD
$$
\par
One can moreover choose the semi-stable models $\Tilde {\Cal Y}$ and $\Tilde {\Cal X}$ such that the
set of points $B_{K'}:=\{x_{i,K'}\}_{1\le i\le s}$ consisting of the branch 
locus in the morphism $f_K':\Cal Y'_{K'}\to \Cal X'_{K'}$ specialise in {\bf smooth distincts} 
points of $\Cal X'_{k}$ and one can choose such $\Tilde {\Cal X}$ and $\Tilde {\Cal Y}$ which are minimal for these 
properties. The fibre $\Tilde {g}^{-1}(x)$ of the closed point $x$ in $\Tilde {\Cal X}$ is a tree $\Gamma$ of projective lines. 
This tree is canonically endowed with some ``{\it degeneration data}'' that we will exhibit below and in the next section, in the cases
where $G\simeq \Bbb Z/p\Bbb Z$ and  $G\simeq \Bbb Z/p^2\Bbb Z$, and
which take into account the geometry of the special fibre $\Tilde {\Cal Y}_k$ of $\Tilde {\Cal Y}$. This will follow
mainly from the results in II, III and IV.

\pop {.5}
\par
\noindent
{\bf \gr 5.2.}\rm \ We will use the same notations as in 5.1. We assume that $G\simeq \Bbb Z/p\Bbb Z$. We consider the case where $\Cal X\simeq \Spf
A$ is the formal germ of a semi-stable $R$-curve at a {\bf smooth} point $x$ i.e. $A\simeq R[[T]]$. Let $R'$ be a finite extension of 
$R$ as in 5.1 and let $\pi '$ be a uniformiser of $R'$. Below we exhibit the {\bf degeneration data} associated to the semi-stable 
reduction of $\Cal Y$.

\pop {.5}
\par
\noindent
{\bf \gr Deg.1.}\rm \ Let $\wp:=(\pi')$ be the ideal of $A':=A\otimes _R R'$ generated by $\pi'$  and let $\hat A'_{\wp}$ be the 
completion of the localisation of $A'$ at $\wp$. Let $\Cal X'_{\eta}:=\Spf \hat A'_{\wp}$ be the formal boundary of $\Cal X'$ and let 
$\Cal X'_{\eta}\to \Cal X'$ be the canonical morphism. Consider the following cartesian diagram:
$$
\CD
{\Cal  Y'}_{\eta} @>f_{\eta}>>  {\Cal X}'_{\eta}  \\                         
   @VVV                @VVV   \\
\Cal Y'    @> f'>> \Cal X' \\
\endCD
$$
\par
Then $f_{\eta}:\Cal Y'_{\eta}\to \Cal X'_{\eta}$ is a torsor under a commutative finite and flat 
$R'$-group scheme $\Cal M_{n,R'}$ of rank $p$ (cf. 2.5.1) for some integer $n\ge 0$. 
Let $(n,m)$ be the degeneration type of the torsor $f_{\eta}$ (cf. 2.5.2) which is canonically associated to $f$. 
The arithmetic genus $g_y$ of the point $y$ equals $(r+m-1)(p-1)/2$ (cf. 3.3.1) where $d_{\eta}:=r(p-1)$ is 
the degree of the divisor of ramification in the morphism $f'_{K'}:\Cal Y'_{K'}\to \Cal X'_{K'}$.

\pop {.5}
\par
\noindent
{\bf \gr Deg.2.}\rm \ The fibre $\Tilde {g}^{-1}(x)$ of the closed point $x$ of $\Cal X'$ in $\Tilde {\Cal X}$
is a tree $\Gamma $ of projective lines. Let $\Vert (\Gamma ):=\{X_i\}_{i\in I}$  be the set of 
irreducible components of $\Tilde {g}^{-1}(x)$  which are the vertices of the tree $\Gamma $.
The tree $\Gamma$ is canonically endowed with an origin vertex $X_{i_0}$ which is the unique irreducible component of 
$\Tilde {g}^{-1}(x)$ which meets the point $x$. We fix an orientation of the tree $\Gamma$ starting from $X_{i_0}$ 
in the direction of the ends.

\pop {.5}
\par
\noindent
{\bf \gr Deg.3.}\rm \ For each $i\in I$ let $\{x_{i,j}\}_{j\in S_i}$ be the set of points of $X_i$ in which 
specialise some point of $B_{K'}$ ($S_i$ may be empty). Also let  $\{z_{i,j}\}_{j\in D_i}$ be the set of  
double points of $\Tilde {\Cal X}_{k}$ supported by $X_i$. 
In particular $x_{i_0,j_0}:=x$ is a double point of $\Tilde {\Cal X}_{k}$. We denote by $B_{k}$ the set of all points 
$\cup _{i\in I}\{x_{i,j}\}_{j\in S_i}$, which is the set of specialisation of the branch locus $B_{K'}$, and by $D_{k}$ 
the set of double points of $\Tilde {\Cal X}_{k}$.
\pop {.5}
\par
\noindent
{\bf \gr Deg.4.}\rm \ Let $\Cal U:=\tilde {\Cal X}-\{B_{k}\cup D_{k}\}$. Let $\{\Cal U_i\}_{i\in I}$ be the set of connected
components of $\Cal U$. The restriction  $f_i:\Cal V_i\to \Cal U_i$ of $f'$ to $\Cal U_i$ is a torsor under a commutative finite and flat $R'$-group 
scheme $\Cal M_{n_i,R'}$ of rank $p$ for some integer $n_i\ge 0$ (cf. 2.2.1) and $f_{i,k}:\Cal V_{i,k}\to \Cal U_{i,k}:={\Cal U}_i\times _{R'}{k}$ 
is a torsor under the $k$-group scheme $\Cal M_{n_i,R'}\times _{R'}{k}$ which is either \'etale isomorphic to $(\Bbb Z/p\Bbb Z)_k$
or radicial isomorphic to $(\alpha_p)_k$. Further when we move in the graph $\Gamma$ (following the above fixed orientation) 
from a fixed vertex $X_i$ in the direction of a vertex $X_{i'}$
such that $n_{i'}=0$ then the corresponding integers $n_i$ decrease strictly as follows from 3.3.10.

\pop {.5}
\par
\noindent
{\bf \gr Deg.5.}\rm \ Each smooth point $x_{i,j}\in B_k$ is endowed via $f$ with a degeneration data on the boundary of the formal fibre at $x_{i,j}$, 
as in Deg.1 above, and which satisfy certain {\it compatibility  conditions}. More precisely for each point 
$x_{i,j}$ we have the reduction type $(n_{i,j}:=n_i,m_{i,j})$ on the boundary of the formal fibre at this point induced by $g$ and such that 
$r_{i,j}=-m_{i,j}+1$ where $r_{i,j}(p-1)$ is the contribution to $d_{\eta}$ of the point which specialise into $x_{i,j}$. In particular $m_{i,j}\le -1$
since $r_{i,j}\neq 0$.

\pop {.5}
\par  
\noindent
{\bf \gr Deg.6.}\rm \ Each double point $z_{i,j}=z_{i',j'}\in X_i\cap X_{i'}$ of $\Tilde X$ with origin vertex $X_i$ and terminal vertex
$X_{i'}$ is endowed with degeneration data $(n_{i,j}:=n_i,m_{i,j})$ and $(n_{i',j'}:=n_{i'},m_{i',j'})$ induced by $g$ on the two 
boundaries of the formal fibre at this point  and we have $m_{i,j}+m_{i',j'}=0$ (cf. 3.3.7). Let $e_{i,j}$ be the thikness of 
the double point $z_{i,j}$. Then $e_{i,j}=pt_{i,j}$ is necessarily divisible by $p$ and we have 
$n_i-n_{i'}=t_{i,j}m_{i,j}$ as follows from 3.3.10. 

\pop {.5}
\par
\noindent
{\bf \gr Deg.7.}\rm \ It follows after easy calculation that: 
$$g_y=\sum _{i\in I_{\et}}(-2+\sum _{j\in S_i}(m_{i,j}+1)+\sum _{j\in D_i}(m_{i,j}+1))(p-1)/2$$   
where $I_{\et}$ is the subset of $I$ consisting of those $i$ for which the torsor $f_i$ is \'etale (i.e. $n_i=0$).

\pop {.5}
\par
\noindent
{\bf \gr 5.2.1. Example.}\rm \ In the following we give an example where one can exhibit the degeneration 
data associated to a Galois cover $f:\Cal Y\to \Cal X$ of degree $p$ where $\Cal X\simeq \Spf R[[T]]$ is the 
formal germ of a smooth point. More precisely for $m>0$ an integer such that both $m$ and $m+1$ are prime to $p$ consider the cover given generically 
by the equation $X^p-X=(T^{-m}+\pi T^{-m-1})$. Here $r=m+2$ and this cover has a reduction of type 
$(0,-m)$ on the boundary. In particular the arithmetic genus $g_y$ of the closed point $y$ of $\Cal Y$ 
equals $(p-1)/2$. The degeneration data associated to the above cover consists necessarily of
a tree with only one vertex and no edges i.e. a unique projective line $X_1$ with a marked point $x_1$ and an \'etale torsor 
$f_1:V_1\to U_1:=X_1-\{x_1\}$ above $U_1$ with conductor $2$ at the point $x_1$. 
\pop {.5}
\par
The above considerations lead naturally to the following abstract geometric and combinatorial definition of degeneration data.
\pop {.5}
\par
\noindent
{\bf \gr 5.2.2. Definition.}\rm \ A {\bf simple degeneration data} Deg(x) of type $(r,(n,m))$ and rank $p$
consists of the following data:
\pop {.5}
\par
\noindent
{\bf \gr Deg.1.}\rm \ $G_{k}$ is a commutative finite and flat $k$-group scheme of rank $p$ which is either \'etale
if $n=0$, or radicial of type $\alpha_p$ otherwise, $r\ge 0$ is an integer, and $m$ is an integer prime to $p$ such that $r+m-1\ge 0$.

\pop {.5}
\par
\noindent
{\bf \gr Deg.2.}\rm\  $\Gamma:=X_{k}$ is an oriented tree of projective lines over $k$ with vertices $\Vert (\Gamma):=\{X_i\}_{i\in I}$ 
which is endowed with an origin vertex $X_{i_0}$ and a marked point $x:=x_{i_0,j_0}$ 
on $X_{i_0}$. We denote by $\{z_{i,j}\}_{j\in D_i}$ the set of double points or (non oriented) 
edges of $\Gamma$ which are supported by $X_i$. We further assume that the orientation of $\Gamma$ is in the direction going form  
$X_{i_0}$ towards its ends.

\pop {.5}
\par
\noindent
{\bf \gr Deg.3.}\rm \ For each vertex $X_i$ of $\Gamma$ is given a set (may be empty) of smooth 
marked points $\{x_{i,j}\}_{j\in S_i}$.

\pop {.5}
\par
\noindent
{\bf \gr Deg.4.}\rm \ For each $i\in I$ is given a torsor $f_{i}:V_{i}\to U_{i}:= X_{i}-\{\{x_{i,j}\}_{j\in S_i}
\cup \{z_{i,j}\}_{j\in D_i}\}$ under a commutative finite and flat $k$-group scheme $G_{k,i}$ of rank $p$, which
is either \'etale or radicial of type $\alpha_p$, with $V_i$ {\bf smooth}. Moreover for each $i\in I$ is given an
integer $n_i$ which equals $0$ if  $f_i$ is \'etale and is positive otherwise.

\pop {.5}
\par
\noindent
{\bf \gr Deg.5.}\rm \ For each $i\in I$ and $j\in S_i$ are given integers $m_{i,j}$ where $m_{i,j}$ is the conductor of the 
torsor $f_i$ at the point $x_{i,j}$ (cf. [S] 1.5) with $m_{i_0,j_0}=-m$. We assume further that $m_{i,j}\le -1$ if $n_i> 0$. 

\pop {.5}
\par
\noindent
{\bf \gr Deg.6.}\rm \ For each double point $z_{i,j}=z_{i',j'}\in X_i\cap X_{i'}$ is given an integer $m_{i,j}$ (resp. $m_{i',j'}$) prime to $p$ 
where $m_{i,j}$ (resp. $m_{i',j'}$) is the conductor of the torsor $f_{i}$ (resp. $f_{i'}$) at the 
point $z_{i,j}$ (resp. $z_{i',j'}$) (cf. [S] 1.3 and 1.5). These data must satisfy $m_{i,j}+m_{i',j'}=0$.

\pop {.5}
\par
\noindent
{\bf \gr Deg.7.}\rm \ For each double point $z_{i,j}=z_{i',j'}\in X_i\cap X_{i'}$ of $\Gamma$ with origin vertex $X_i$ is given an integer 
$e_{i,j}=pt_{i,j}$ divisible by $p$ such that with the same notations as above we have $n_i-n_{i'}=m_{i,j}t_{i,j}$.
Moreover associated to $x$ is an integer $e=pt$ such that $n-n_{i_0}=mt$.

\pop {.5}
\par
\noindent
{\bf \gr Deg.8.}\rm \  Let $I_{\et}$ be the subset of $I$ consisting of those $i$ for which
$G_{k,i}$ is \'etale. Then the following equality should hold: 
$(r+m-1)(p-1)/2=\sum _{i\in I_{\et}}(-2+\sum _{j\in S_i}(m_{i,j}+1)+\sum _{j\in D_i}(m_{i,j}+1))(p-1)/2$. The integer 
$g:=(r+m-1)(p-1)/2$ is called the {\bf genus} of the degeneration data Deg(x).

\pop {.5}
\par
There is a natural notion of isomorphism of simple degeneration data of a given type and rank $p$. 
We will denote by $\bold {Deg_p}$ the set of isomorphism classes of simple degeneration data of rank $p$. 
The discussion in 5.2 can be reinterpreted as follows:

\pop {.5}
\par
\noindent
{\bf \gr 5.2.3. Proposition.}\rm \ {\sl Let $\Cal X$ be the germ of a formal $R$-curve at a smooth point $x$ and
 let $f:\Cal Y\to \Cal X$ be a cyclic $p$-cover with $\Cal Y$ normal and local. Then one can associate to $f$ canonically 
a simple degeneration data $\Deg(x)\in \bold {Deg_p}$ which describes the semi-stable reduction of $\Cal Y$. In other words there exists 
a canonical ``specialisation'' map $\Sp:H^1_{\et}(\Spec L,\Bbb Z/p\Bbb Z)\to \bold {Deg_p}$ where $L$ is the function field of the geometric fibre of $\Cal X$.}

\pop {.5}
\par
Reciprocally we have the following result of realisation for degeneration data for such covers:

\pop {.5}
\par
\noindent
{\bf \gr 5.2.4. Theorem.}\rm \ {\sl The above specialisation map $\Sp:H^1_{\et}(\Spec L,\Bbb Z/p\Bbb Z)\to \bold {Deg_p}$
defined in 5.2.3 is surjective.}

\pop {.5}
\par
\noindent
{\bf \gr Proof.}\rm \ Consider a degeneration data $\Deg(x)\in \bold {Deg_p}$. We have to show that $\Deg (x)$ is associated to
some cyclic $p$-cover above the formal germ of a smooth $R$-curve after eventually enlarging $R$. We assume that the degeneration data is of type
$(r,(n,m))$. We treat only the case $n=0$, the case where $n>0$ is treated in a similar way. The proof is done by induction on the length of the tree $\Gamma$ of 
$\Deg (x)$. Assume first that the tree $\Gamma$ has minimal length and consists of 
one irreducible component $X=\Bbb P^1_k$ with one marked (double) point $x$ and $r$ smooth distinct marked points $\{x_i\}_{i=1}^r$. Let $U:=X-\{x,x_i\}_i$ and let 
$\bar f:V\to U$ be the torsor given by the data $\Deg.4$ which is necessarily an $\alpha_p$-torsor by 3.3.8 1) (i.e. the integer $n_i:=n'$ associated to the
vertex $X$ in $\Deg.4$ is non zero). First for each $i\in \{1,...,r\}$ consider the formal germ 
$\Cal X_i:=\Spf R[[T_i]]$ and the cyclic $p$-cover $f_i:\Cal Y_i\to \Cal X_i$ given by the equation $Y_i^p-\pi ^{n'}Y_i=T_i^{m_i}$ where $m_i$ is the ``conductor''
associated to the point $x_i$ in $\Deg.5$ and $n'$ is the positive integer associated to $\bar f$ in $\Deg.4$. Let $\Cal X$ be a formal projective $R$-line with special fibre 
$X$. Let $X':=X-\{x\}$ and let $\Cal X'$ (resp $\Cal U$) be the formal fibre of $X'$ (resp. of $U$) in $\Cal X$. The torsor $\bar f$ is given by an equation $t^p=\bar u$
where $\bar u$ is a regular function on $U$. Let $u$ be a regular function on $\Cal U$ which lifts $\bar u$. Then the cover $f:\Cal V\to \Cal U$ given by the equation
$Y^p-\pi^{n'(p-1)}Y=u$ is a torsor under the $R$-group scheme $\Cal M_n$ which lifts the $\alpha_p$-torsor $\bar f$. By construction the torsor $f$ has a reduction on the 
formal boundary at each point $x_i$ of type $(n',m_i)$ which coincides with the degeneration type of the cover $f_i$ above the boundary of $\Cal X_i$. The technique of formal
patching (cf. [S-1], 1) allows then one to construct a $p$-cyclic cover $f':\Cal Y'\to \Cal X'$ which restricted to $\Cal U$ is isomorphic to $f'$ and restricted to 
$\Cal X_i$, for each $i\in \{i,...,r\}$, is isomorphic to $f_i$ (cf. loc. cit.). Let $\Cal X_1\to \Cal X$ be the blow up of $\Cal X$ at the point $x$
and let $X_1$ be the exeptional fibre in $\Cal X_1$ which meets $X$ at the double point $x$. Let $e=pt$ be the integer associated to the marked double point $x$
via $\Deg.7$. After enlarging $R$ we can assume that the double point $x$ of $\Cal X_1$ has thikness $e$. We have $-n'=mt$ by assumption. 
Let $\Cal X_1'$ be the formal fibre of $X_1':=X_1-\{x\}$ in $\Cal X_1$. Let $f_1':\Cal Y_1'\to \Cal X_1'$ be the \'etale $\Bbb Z/p\Bbb Z$-torsor given
by the equation ${Y'}^p-Y'=h^m$ where $h$ is a ``parameter'' on $\Cal X_1'$. Further let $\Cal X_{1,x}\simeq \Spf R[[S,T]]/(ST-\pi^{pt})$ be the formal germ of 
$\Cal X_1$ at the double point $x$. Consider the cover $f_{1,x}:\Cal Y_{1,y} \to \Cal X_{1,x}$ given by the equation $Y^p-Y=S^m=\pi^{ptm}T^{-m}$. 
Then $\Cal Y_{1,y}$ is the formal germ of a double point of thikness $t$ (cf. 3.3.9. a). Moreover the cover $f_1'$ (resp. $f'$) has the same degeneration type 
(by construction) on the boundary corresponding to the double point $x$ as the degeneration type of the cover $f_{1,x}$ above the formal boundary with parameter 
$T$ (resp. above the formal boundary with parameter $S$). A second application of the formal patching techniques allows one to construct a $p$-cyclic cover
$f_1:\Cal Y_1\to \Cal X_1$ which restricted to $\Cal X'$ (resp. $\Cal X_1'$ and $\Cal X_{1,x}$) is isomorphic to $f'$ (resp. to $f_1'$ and $f_x$).
Let $\tilde {\Cal X}$ be the $R$-curve obtained by contracting the irreducible component $X$ in $\Cal X_1$. We denote the image of the double point $x$ in
$\tilde {\Cal X}$ simply by $x$.
The cover $f_1:\Cal Y_1\to \Cal X_1$ induces canonically a $p$-cyclic cover $\tilde f:\tilde {\Cal Y}\to \tilde {\Cal X}$ above $\tilde {\Cal X}$. 
Let $\Cal X_x\simeq R[[T]]$ be the formal
germ of $\tilde {\Cal X}$ at the smooth point $x$. Then $\tilde f$ induces canonically a $p$-cyclic cover $f_{x}:\Cal Y_{y} \to \Cal X_{x}$ where $y$ is the closed point of 
$\Cal Y$ above $x$. Now it is easy to see that the degeneration data associated to $f_x$ via 5.2.3 is isomorphic to the degeneration data $\Deg (x)$ 
we started with. Finally the proof in the general case is very similar and is left to the reader. The only modification in the proof above 
is that one has to considers $p$-cyclic covers
$f_i:\Cal Y_i\to \Cal X_i$ above the formal germs $\Cal X_i:=\Spf R[[T_i]]$ which one obtains by induction hypothesis as realisation of the degeneration data, 
induced by $\Deg (x)$, on the subtrees $\Gamma_i$ of $\Gamma$ which starts from the edge $x_i$ in the direction of the ends and which clearly has length smaller
than the length of $\Gamma$.

\pop {.5}
\par
\noindent
{\bf \gr 5.2.5. Remarks.}

\pop {.5}
\par
\noindent
{\bf \gr 1.}\rm\ One can also define in the same way as in 5.2.2 and using the results of II, III, and IV, the set of isomorphism classes of ``double'' degeneration data 
associated to the minimal semi-stable model of $p$-cyclic covers $f:\Cal Y\to \Cal X$ above the formal germ of a double point. Moreover one can prove, in a similar way
as in 5.2.4, a result of realisation for such a degeneration data.

\pop {.5}
\par
\noindent
{\bf \gr 2.}\rm\ In [M] Maugeais proved (in equal characteristic $p>0$) in theorem 5.4 a global result of lifting for finite ``admissible'' covers of 
degree $p$ between semi-stable curves. The methods used 
in the proof of 5.2.4 are essentially the same as he uses but more direct in the sens that the lemmas 4.2, 4.4, and corollary 4.3 he uses are avoided and we use instead
our results 3.3.3 and 3.3.8 which are a direct consequence of the computation of vanishing cycles.

\pop {.5}
\par
\noindent
{\bf \gr 3.}\rm\ It is easy to construct examples of $p$-cyclic covers 
$f:\Cal Y\to \Cal X$ above the forml germ of a smooth point $\Cal X$
where the special fibre $\Cal Y_k$ is singular and unibranche at the closed point $y$ of $\Cal Y$ 
and such that the configuration of the special fibre of a semi-stable model $\Tilde {\Cal Y}$ of $\Cal Y$ is 
not a tree-like. More precisely consider the simple degeneration data $\Deg (x)$
of type $(n,m)$, with $n>0$, $m>0$, and $n=mt$, which consists of a graph $\Gamma$ with two
vertices $X_1$ and $X_2$ linked by a unique edge $\tilde x$ with given 
marked points $x=x_1$ on $X_1$ and $x_2$ on $X_2$. Further $X_1$ is the original component of $\Gamma$. 
As part of the data are given \'etale torsors of rank $p$:
$f_1:V_1\to U_1:=X_1-\{x_1\}$ with conductor $m$ at $x=x_1$ and  
$f_2:V_2\to U_2:=X_2-\{x_2\}$ with conductor $m'$ at $x_2$. Also is given the thikness $e=pt$ at the 
``double'' point $x$ with $n=tm$. Then it follows from 5.2.4 that there exists after eventually a finite extension 
of $R$ a Galois cover $f:\Cal Y\to \Cal X$ of degree $p$ above the formal 
germ $\Cal X\simeq \Spf R[[T]]$ at the smooth $R$-point $x$ such that the 
simple degeneration data associated to the above cover $f$ is the above 
given one. Moreover by construction the singularity of the closed point $y$ 
of $\Cal Y$ is unibranche and the configuration of the 
semi-stable reduction of $\Cal Y$ consists of two projective curves which 
meet at $p$-double points (the above cover will be \'etale in reduction 
above the double point $\tilde y$). In particular one has $p-1$ cycles in this 
configuration.

\pop {1}
\par
\noindent
{\bf \gr VI. Semi-stable reduction of cyclic $p^2$-covers above formal germs of curves in equal characteristic $p>0$.}

\pop {.5}
\par
\noindent
{\bf \gr 6.1.}\rm \ In all this paragraph we use the same notations as in V.

\pop {.5}
\par
\noindent
{\bf \gr 6.2.}\rm \ Let ${\Cal X}:=\Spf \hat {\Cal O}_{X,x}$ be the 
formal germ of an $R$-curve $X$ at a closed point $x$ and let $f:\Cal Y\to \Cal X$ be a Galois cover with 
group $G$ such that  $\Cal Y$ is normal and local. We assume that $G\simeq \Bbb Z/p^2\Bbb Z$ is cyclic of order $p^2$ and we use the same notations as 
in 5.1 for the minimal semi-stable model $\tilde \Cal Y$ of $\Cal Y$. We consider the case where $\Cal X\simeq \Spf R[[T]]$ is the formal germ of a 
semi-stable $R$-curve at a {\bf smooth} point $x$. Let $R'$ be a finite extension of 
$R$ as in 5.1 and let $\pi '$ be a uniformiser of $R'$. Below we exhibit the {\bf degeneration data} associated to the semi-stable 
reduction of $\Cal Y$ and which are consequences of the results in II, III, and IV. In all what follows, and in order to make the statements of our results
simpler, we make the following assumption (**):
{\it the special fibre of $\Cal Y$ is irreducible and for each irreducible component $X_i$ of the special fibre $\tilde {\Cal X}$ 
the fibre $g^{-1}(X_i)$ of $X_i$ in $\tilde {\Cal Y}$ is irreducible}.

\pop {.5}
\par
\noindent
{\bf \gr Deg.1.}\rm \ Let $\wp:=(\pi')$ be the ideal of $A':=A\otimes _R R'$ generated by $\pi'$  and let $\hat A'_{\wp}$ be the 
completion of the localisation of $A'$ at $\wp$. Let $\Cal X'_{\eta}:=\Spf \hat A'_{\wp}$ be the boundary of $\Cal X'$ and let 
$\Cal X'_{\eta}\to \Cal X'$ be the canonical morphism. Consider the following cartesian diagram:
$$
\CD
{\Cal  Y'}_{\eta} @>f_{\eta}>>  {\Cal X}'_{\eta}  \\                         
   @VVV                @VVV   \\
\Cal Y'    @> f'>> \Cal X' \\
\endCD
$$
\par
Then $f_{\eta}:\Cal Y'_{\eta}\to \Cal X'_{\eta}$ is a finite cyclic cover of degree $p^2$.
Let $\{(n_1,m_1),(n_2,m_2)\}$ be the degeneration type of the cover $f$ (cf. 2.6.2) which is canonically associated to $f$ and which 
is an admissible pair by 2.6.4. The arithmetic genus $g_y$ of the point $y$ equals $(r+pm_1+m_2-p-1)(p-1)/2$ (cf. 4.2.1) where $d_{\eta}:=r(p-1)$ is 
the degree of the divisor of ramification in the morphism $f'_{K'}:\Cal Y'_{K'}\to \Cal X'_{K'}$.

\pop {.5}
\par
\noindent
{\bf \gr Deg.2.}\rm \ The fibre $\Tilde {g}^{-1}(x)$ of the closed point $x$ of $\Cal X'$ in $\Tilde {\Cal X}$
is a tree $\Gamma $ of projective lines. Let $\Vert (\Gamma ):=\{X_i\}_{i\in I}$  be the set of 
irreducible components of $\Tilde {g}^{-1}(x)$  which are the vertices of the tree $\Gamma $.
The tree $\Gamma$ is canonically endowed with an origin vertex $X_{i_0}$ which is the unique irreducible component of 
$\Tilde {g}^{-1}(x)$ which meets the point $x$. We fix an orientation of the tree $\Gamma$ starting from $X_{i_0}$ 
towards its ends.

\pop {.5}
\par
\noindent
{\bf \gr Deg.3.}\rm \ For each $i\in I$ let $\{x_{i,j}\}_{j\in S_i}$ be the set of points of $X_i$ in which 
specialise some point of $B_{K'}$ ($S_i$ may be empty). Also let  $\{z_{i,j}\}_{j\in D_i}$ be the set of  
double points of $\Tilde {\Cal X}_{k}$ supported by $X_i$. 
In particular $x_{i_0,j_0}:=x$ is a double point of $\Tilde {\Cal X}_{k}$. We denote by $B_{k}$ the set of all points 
$\cup _{i\in I}\{x_{i,j}\}_{j\in S_i}$, which is the set of specialisation of the branch locus $B_{K'}$, and by $D_{k}$ 
the set of double points of $\Tilde {\Cal X}_{k}$.
\pop {.5}
\par
\noindent
{\bf \gr Deg.4.}\rm \ Let $\Cal U:=\tilde {\Cal X}-\{B_{k}\cup D_{k}\}$. Let $\{\Cal U_i\}_{i\in I}$ be the set of connected
components of $\Cal U$. The restriction  $f_i:\Cal V_i\to \Cal U_i$ of $f'$ to $\Cal U_i$ is a finite cyclic cover of degree $p^2$ which factorises 
canonically as $\Cal V_i @>f_{i,2}>> \Cal V_i' @>f_{i,1}>> \Cal U_i$ where $f_{i,2}$ (resp.  $f_{i,1}$) is a torsor under the $R$-group scheme
$\Cal M_{n_{i_2}}$ (resp. a torsor under the $R$-group scheme $\Cal M_{n_{i_1}}$). Moreover to the cover $f_i$ is canonically associated a 
degeneration data $\Cal D_i$ of type $A$, $B$, or $C$ and which describes the special fibre $f_{i,k}:\Cal V_{i,k}\to \Cal U_{i,k}:=
{\Cal U}_i\times _{R'}{k}$ of $f_i$ (cf. 2.4.8). The pair $(n_{i_1},n_{i_2})$ satisfies $n_{i_2}\ge n_{i_1}(p(p-1)+1)$ by 2.4.3.

\pop {.5}
\par
\noindent
{\bf \gr Deg.5.}\rm \ Each smooth point $x_{i,j}\in B_k$ is endowed via $f$ with a degeneration data on the boundary of the formal fibre at $x_{i,j}$, 
as in Deg.1 above, and which satisfy certain {\it compatibility  conditions}. More precisely for each point 
$x_{i,j}$ we have the reduction type $\{(n_{i_1,j_1}:=n_{i_1},m_{i_1,j_1}),(n_{i_2,j_2}:=n_{i_2},m_{i_2,j_2})\}$ on the boundary of the formal fibre at 
this point induced by $g$, which is an admissible pair satisfying the condition (*) in 4.2.5, and such that 
$r_{i,j}=-pm_{i_1,j_1}-m_{i_2,j_2}+p+1$ where $r_{i,j}(p-1)$ is the contribution to $d_{\eta}$ of the point which specialise into $x_{i,j}$.

\pop {.5}
\par  
\noindent
{\bf \gr Deg.6.}\rm \ Each double point $z_{i,j}=z_{i',j'}\in X_i\cap X_{i'}$ of $\Gamma$ with origin vertex $X_i$ and terminal vertex
$X_{i'}$ is endowed with degeneration data $\{(n_{i_1,j_1}:=n_{i_1},m_{i_1,j_1}),(n_{i_2,j_2}:=n_{i_2},m_{i_2,j_2})\}$ and 
$\{(n_{i'_1,j'_1}:=n_{i'_1},m_{i'_1,j'_1}),(n_{i'_2,j'_2}:=n_{i'_2},m_{i'_2,j'_2})\}$ induced by $g$ on the two 
boundaries of the formal fibre at this point. These pairs are admissible. Moreover we have $m_{i_1,j_1}+m_{i'_1,j'_1}=0$ 
and $m_{i_2,j_2}+m_{i'_2,j'_2}=0$ (cf. 4.2.9). Let $e_{i,j}$ be the thikness of 
the double point $z_{i,j}$. Then $e_{i,j}=p^2t_{i,j}$ is necessarily divisible by $p^2$ and we have 
$n_{i_1}-n_{i'_1}=pt_{i,j}m_{i_1,j_1}$ and $n_{i_2}-n_{i'_2}=t_{i,j}m_{i_2,j_2}$ as follows from 3.3.10.

\pop {.5}
\par
\noindent
{\bf \gr Deg.7.}\rm \ It follows after easy calculation that: 
\newline
$g_y=\sum _{i\in \tilde I}(-2+\sum _{j\in S_i}(d_{i,j})+\sum _{j\in D_i}(d_{i,j}))(p-1)/2$   
where $\tilde I$ denote the subset of $I$ consisting of those $i$ for which the cover $f_i$ has a reduction of type (\'etale, \'etale) or (\'etale, radicial) (cf. 2.4.8).
Here $d_{i,j}=(m_{i_1}+1)+p(pm_{i_1}+1)$ if $pm_{i_1}>m_{i_2}$ and $d_{i,j}=(m_{i_1}+1)+p(m_{i_2}+1)$ otherwise in the case of 
(\'etale, \'etale) reduction type (resp. $d_{i,j}=(m_{i_1}+1)$ in the case of (\'etale, radicial) reduction type).

\pop {.5}
\par
The above considerations lead naturally to the following abstract geometric and combinatorial definition of degeneration data.
\pop {.5}
\par
\noindent
{\bf \gr 6.2.1. Definition.}\rm \ A {\bf simple irreducible degeneration data} Deg(x) of type $(r,\{(n_1,m_1),(n_2,m_2)\})$ and rank $p^2$
consists of the following data:
\pop {.5}
\par
\noindent
{\bf \gr Deg.1.}\rm \  The pair $(n_1,n_2)$ satisfies $n_2\ge n_1(p(p-1)+1)/p$, $r\ge 0$ is an integer, and $m_1$ (resp. $m_2$) is an integer 
prime to $p$ such that $r+pm_1+m_2-p-1\ge 0$. We further assume that the pair $\{(n_1,m_1),(n_2,m_2)\}$ is admissible (cf. 2.6.4).

\pop {.5}
\par
\noindent
{\bf \gr Deg.2.}\rm\  $\Gamma:=X_{k}$ is an oriented tree of projective lines over $k$ with vertices $\Vert (\Gamma):=\{X_i\}_{i\in I}$ 
which is endowed with an origin vertex $X_{i_0}$ and a marked point $x:=x_{i_0,j_0}$ 
on $X_{i_0}$. We denote by $\{z_{i,j}\}_{j\in D_i}$ the set of double points or (non oriented) 
edges of $\Gamma$ which are supported by $X_i$. Further for each vertex $X_i$ of $\Gamma$ is given a set (may be empty) of smooth marked points $\{x_{i,j}\}_{j\in S_i}$. 
We further assume that the orientation of $\Gamma$ is in the direction going form $X_{i_0}$ towards its ends.

\pop {.5}
\par
\noindent
{\bf \gr Deg.3.}\rm \ For each $i\in I$ is given a degeneration data $\Cal D_i$ of type $A$, $B$, or $C$ as in 2.4.8 associated to
$U_{i}:= X_{i}-\{\{x_{i,j}\}_{j\in S_i}\cup \{z_{i,j}\}_{j\in D_i}\}$
which is of type (\'etale, \'etale), (\'etale, radicial), or (radicial, radicial) respectively.  
Moreover for each $i\in I$ is given a pair of integers $(n_{i_1},n_{i_2})$ which equals $(0,0)$ and $(0,n_{i_2})$ in the cases $A$ and $B$ 
respectively and which satisfies the inequality $n_{i_2}\ge n_{i_1}(p(p-1)+1)/p$.

\pop {.5}
\par
\noindent
{\bf \gr Deg.3'.}\rm \ To the degeneration data $\Cal D_i$ above and the pair of integers  $(n_{i_1},n_{i_2})$ is associated a finite 
cover $\bar f_i:V_i\to U_i$ of degeree $p^2$ which factorises as $V_i@>\bar f_{i,2}>>V_{i,1}@>\bar f_{i,1}>>U_i$ where $\bar f_{i,2}$ and 
$\bar f_{i,1}$ have the structure of a torsor under a finite and flat group scheme of 
rank $p$. In the case where  $\Cal D_i$ is of type A or B then  $\bar f_i$ is simply the corresponding torsor. If $\Cal D_i$ is of type $C$ given by a pair
$(\bar u_1,\bar u_2)$ of regular functions on $U_i$, defined up to addition of $(\bar v_1^p,\bar v_2^p)$, and a tuple $(\bar c_1,...,\bar c_{p-1})$ of regular 
functions on $U_i$. Then the cover $\bar f_i$ is defined by the equations: $t_1^p=\bar u_1$ and $t_2^p=\bar u_2-\bar c_1^pt_1-2\bar c_2^pt_1^{p+1}-...-(p-1)
\bar c_{p-1}t_1^{p(p-2)+1}$ if $n_{i_2}> n_{i_1}(p(p-1)+1)/p$ (resp. $t_2^p=\bar u_2-\bar c_1^pt_1-2\bar c_2^pt_1^{p+1}-...-(p-1)
\bar c_{p-1}t_1^{p(p-2)+1}-t_1^{p(p-1)+1}$ if $n_{i_2}=n_{i_1}(p(p-1)+1)/p$).

\pop {.5}
\par
\noindent
{\bf \gr Deg.4.}\rm \ For each $i\in I$ and $j\in S_i$ are given integers $(m_{i_1,j_1},m_{i_2,j_2})$ such that $m_{i_1,j_1}$ (resp. $m_{i_2,j_2}$)
is the conductor of the torsor $\bar f_{i,1}$ (resp. the torsor $\bar f_{i,2}$) above in $\Deg 3'$ at the point $x_{i,j}$ (resp. the point of $V_{i,1}$ above
$x_{i,j}$). we assume that
the pair $\{(n_{i_1},m_{i_1,j_1}),(n_{i_2},m_{i_2,j_2})\}$ satisfies the condition (*) in 4.2.5. Moreover we assume that
$(m_{{i_0}_1,{j_0}_1},m_{{i_0}_2,{j_0}_2})$
\newline
$=(-m_1,-m_2)$. 

\pop {.5}
\par
\noindent
{\bf \gr Deg.5.}\rm \ For each double point $z_{i,j}=z_{i',j'}\in X_i\cap X_{i'}$ is given a pair of integers
\newline 
$(m_{i_1,j_1},m_{i_2,j_2})$ (resp. $(m_{i'_1,j'_1},m_{i'_2,j'_2})$) prime to $p$ such that $m_{i_1,j_1}$ and $m_{i_2,j_2}$
are the conductors of the torsor $\bar f_{i,1}$ (resp. the torsor $\bar f_{i,2}$) above at the point $z_{i,j}$
(resp. $m_{i'_1,j'_1}$ and $m_{i'_2,j'_2}$ are the conductors of the torsor $\bar f_{i',1}$ (resp. the torsor $\bar f_{i',2}$) above at the point $z_{i,j}$). 
We assume further that the pair $\{(n_{i_1},m_{i_1,j_1}),(n_{i_2},m_{i_2,j_2})\}$ (resp.$\{(n_{i'_1},m_{i'_1,j'_1}),(n_{i'_2},m_{i'_2,j'_2})\}$) is admissible.
These data must satisfy $m_{i_1,j_1}+m_{i'_1,j'_1}=m_{i_2,j_2}+m_{i'_2,j'_2}=0$.

\pop {.5}
\par
\noindent
{\bf \gr Deg.6.}\rm \ For each double point $z_{i,j}=z_{i',j'}\in X_i\cap X_{i'}$ of $\Gamma$ with origin vertex $X_i$ is given an integer 
$e_{i,j}=p^2t_{i,j}$ divisible by $p^2$ such that with the same notations as above we have $n_{i_1}-n_{i'_1}=pt_{i,j}m_{i_1,j_1}$ 
and $n_{i_2}-n_{i'_2}=t_{i,j}m_{i_2,j_2}$.

\pop {.5}
\par
\noindent
{\bf \gr Deg.7.}\rm \  Let $\tilde I$ denote the subset of $I$ consisting of those $i$ for which the degeneration data $\Cal D_i$ $f_i$ is of type 
(\'etale, \'etale) or (\'etale, radicial). Then the following equality should hold: 
$(r+pm_1+m_2-p-1)(p-1)/2 =\sum _{i\in I_{\et}}(-2+\sum _{j\in S_i}(d_{i,j})+\sum _{j\in D_i}(d_{i,j}))(p-1)/2$ where
$d_{i,j}=(m_{i_1}+1)+p(pm_{i_1}+1)$ if $pm_{i_1}>m_{i_2}$ and $d_{i,j}=(m_{i_1}+1)+p(m_{i_2}+1)$ otherwise in the case of 
(\'etale, \'etale) type (resp. $d_{i,j}=(m_{i_1}+1)$ in the case of (\'etale, radicial) type).

\pop {.5}
\par
We further make the following important assumption as part of the degeneration data:

\pop {.5}
\par
\noindent
{\bf \gr Deg.8.}\rm \  For each $i\in I$  let $U_{i}$ be as in $\Deg.3$ and let $\Cal U_i$ be a formal affine 
$R$-scheme with special fibre $U_i$. We assume that for each $i\in I$ there exists a lifting of the degeneration data $\Cal D_i$ as in 2.4.9 to a $p^2$-cyclic cover 
$f_i:\Cal V_i\to \Cal U_i$, which is an \'etale torsor above the generic fibre of $\Cal U_i$, such that the following ``compatibility'' condition holds: 
for each double point $z:=z_{i,j}=z_{i',j'}\in X_i\cap X_{i'}$ the restriction of $f_i$ (resp. restriction of $f_{i'}$) to the boundary 
$\Cal X_{z,i}\simeq \Spf R[[T]]\{T^{-1}\}$ corresponding to the double point $z$ (resp. the boundary 
$\Cal X_{z,i'}\simeq \Spf R[[S]]\{S^{-1}\}$ corresponding to the point $z$) are Galois isomorphic.

\pop {.5}
\par
There is a natural notion of isomorphism of simple irreducible degeneration data of a given type and rank $p^2$. 
We will denote by $\bold {Deg_{p^2}}$ the set of isomorphism classes of simple irreducible degeneration data of rank $p^2$. 
The discussion in 6.2 can be reinterpreted as follows:

\pop {.5}
\par
\noindent
{\bf \gr 6.2.2. Proposition.}\rm \ {\sl Let $\Cal X$ be the germ of a formal $R$-curve at a smooth point $x$ and
let $f:\Cal Y\to \Cal X$ be a cyclic $p^2$-cover with $\Cal Y$ normal and local. We assume that the condition $(**)$ in 6.2 holds. 
Then one can associate to $f$ canonically a simple irreducible degeneration data $\Deg(x)\in \bold {Deg_{p^2}}$. In other words there exists a 
canonical ``specialisation'' map $\Sp:H^1_{\et}(\Spec L,\Bbb Z/p^2\Bbb Z)^{\ir}\to \bold {Deg_{p^2}}$ where $L$ is the function field of the geometric fibre 
of $\Cal X$ and $H^1_{\et}(\Spec L,\Bbb Z/p^2\Bbb Z)^{\ir}$ is the subset of $H^1_{\et}(\Spec L,\Bbb Z/p^2\Bbb Z)$ consisting of those covers satisfying 
the condition $(**)$.}

\pop {.5}
\par
Reciprocally we have the following result of realisation for degeneration data for such covers:

\pop {.5}
\par
\noindent
{\bf \gr 6.2.3. Theorem.}\rm \ {\sl The above specialisation map $\Sp:H^1_{\et}(\Spec L,\Bbb Z/p^2\Bbb Z)^{\ir}\to \bold {Deg_{p^2}}$
defined in 6.2.2 is surjective.}

\pop {.5}
\par
\noindent
{\bf \gr Proof.}\rm\ Consider a simple irreducible degeneration data $\Deg(x)\in \bold {Deg_{p^2}}$. We have to show that $\Deg (x)$ is associated to
some cyclic $p^2$-cover above the formal germ of a smooth $R$-curve after eventually enlarging $R$. We assume that the degeneration data is of type
$(r,\{(n_1',m_1),(n_2',m_2)\})$. We treat only the case where $n_1'=n_2'=0$. The remaining cases are treated in a similar way. 
As in the proof of 5.2.4 we proceed by induction on the length of the tree $\Gamma$ of $\Deg (x)$. Assume first that the tree $\Gamma$ has minimal length and consists of 
one irreducible component $X=\Bbb P^1_k$ with one marked (double) point $x=x_{i_0}$ and $r$ smooth marked points $\{x_i\}_{i=1}^r$. 
Let $\Cal X$ be a formal projective $R$-line with special fibre 
$X$. Let $X':=X-\{x\}$ and $U:=X-\{x,x_i\}_i$. Further let $\Cal X'$ (resp $\Cal U$) be the formal fibre of $X'$ (resp. of $U$) in $\Cal X$.
Let $\Cal D$ be the degeneration data given by $\Deg.3$. This data is necessarily of type $C$ (i.e (radicial, radicial)), since $n'_1=n'_2=0$,
hence is given by an element of $H^1_{\fppf}(U,H_k)\oplus \Gamma (U,\Cal O_{X_k})^{p-1}$ where $H_k$ is the finite commutative group scheme extension of 
$\alpha_p$ by $\alpha_p$ defined in 2.3.1. The degeneration data $\Cal D$ is thus given by a pair of functions $(\bar u_1,\bar u_2)$ on $U$, defined up 
to addition of $(\bar v_1^p,\bar v_2^p)$, and a tuple $(\bar c_1,...,\bar c_{p-1})$ of regular functions on $U$. Also let $(n_1,n_2)$ be the pair of integers
associated to the vertex $X$ via $\Deg.3$ and which satisfy $n_2\ge n_1(p(p-1)+1)/p$. We assume for simplicity that $n_1=p\tilde n_1$ and that $n_2=\tilde n_1(p(p-1)+1)$.
By 2.4.9 there exists a $p^2$-cyclic cover $f':\Cal V\to \Cal U$ which gives rise to the degeneration data $\Cal D$ by 2.4.8. For example we can choose the 
cover $f'$ to be generically given by the equations: $(T_1^p,T_2^p)-(T_1,T_2)=(u_1\pi^{-\tilde n_1p},f(u_1)\pi^{-p^2\tilde n_1}+u_2\pi^{-p^2\tilde n_1+\tilde n(p-1)})$ 
where $u_1$ (resp. $u_2$) is a regular function on $\Cal U$ which lifts $\bar u_1$ (resp. $\bar u_2$) and $f(u_1)=c_1^pu_1+c_2^pu_1^2+...+c_{p-1}u_1^{p-1}$
where $c_i$ is a regular function on $\Cal U$ which lifts $\bar c_i$.

\par
For each $i\in \{1,...,r\}$ let $\Cal X'_i\simeq R[[T_i]]\{T_i^{-1}\}$ be the formal boundary of $\Cal U$ at the point $x_i$. The cover $f'$ induces canonically
a cover $f'_i:\Cal Y'_i\to \Cal X'_i$ which by construction has a reduction of type $\{(n_1,m_{i,1}),(n_2,m_{i,2})\}$:= the admissible pair
satisfying (*) which is associated to the point $x_i$ via $\Deg. 4$. Note that we are necessarily in the case
c-3) of 2.6.1 and the cover $f_i$ is generically given by the equations:
$$X_1^p-X_1=T_i^{m_{i,1}}/\pi ^{\tilde n_1p}$$ 
and: 
$$X_2^p-X_2=f(T_i)/\pi^{p^2\tilde n_1}+p^{-1}\sum _{k=1}^{p-1}{p \choose k}X_1^{pk}(-X_1)^{p-k}$$ 
where $f(T_i)=\sum _{j\in I} a_jT_i^j\in R[[T_i]]\{T_i^{-1}\}$ is not a $p$-power modulo $\pi$. Let $m_2$ be the conductor of $f$ then 
$m_{i,2}=\inf (m_{i,1}(p(p-1)+1),m_2p-m_{i,1}(p-1))$ (cf. loc. cit.). Further we can assume after some transformation which eliminate the $p$-powers that
$f(T_i)=\sum _{j\ge m_2} a_jT_i^j$. Let $\Cal X_i\simeq R[[T_i]]$ be the formal open disc with parameter $T_i$ and consider the $p^2$-cyclic cover
$f_i :\Cal Y_i\to \Cal X_i$ which is generically given by the equations:
$$X_1^p-X_1=T_i^{m_{i,1}}/\pi ^{\tilde n_1 p}$$ 
and: 
$$X_2^p-X_2=f(T_i)/\pi^{p^2\tilde n_1}+p^{-1}\sum _{k=1}^{p-1}{p \choose k}X_1^{pk}(-X_1)^{p-k}$$

Then $\Cal Y_i$ is smooth (cf. 4.2.3 4)). Now by construction the covers $f_i$ and $f'$ coincide when restricted to the formal fibre $\Cal X'_i$. The technique of formal
patching (cf. [S-1], 1) allows then one to construct a $p^2$-cyclic cover $\tilde f:\Cal Y'\to \Cal X'$ which restricted to $\Cal U$ is isomorphic to $f'$ 
and restricted to $\Cal X_i$, for each $i\in \{i,...,r\}$, is isomorphic to $f_i$ (cf. loc. cit.).

Let $\Cal X_{i_0}\simeq \Spf R[[T]]\{T^{-1}\}$ be the formal boundary of $\Cal U$ at the point $x=x_{i_0}$. Then the cover $\tilde f$ induces canonically a  
$p^2$-cyclic cover $f_{i_0}:\Cal Y_{i_0}\to \Cal X_{i,0}$ which by construction has a reduction of type $\{(n_1,-m_1=m_{i_{0,1},j_{0,1}}),(n_2,-m_2:=
m_{i_{0,2},j_{0,2}})\}$:=the pair associated to the point $x_{i_0}$ via $\Deg. 4$. Note that $m_1$ and $m_2$ must be negative since we assumed that $n_1=n_2=0$.
Moreover because $-n_1=m_1pt$ and $-n_2=m_2t$ by $\Deg.6$ and because we assumed $n_2=n_1(p(p-1)+1)$ we see that necessarily $m_2=m_1(p(p-1)+1)$.
The cover $f_{i_0}$ is then generically given by 2.6.1 c) by equations: 
$$X_1^p-X_1=T^{-m_{1}}/\pi ^{n_1p}$$ 
and: 
$$X_2^p-X_2=f(T)/\pi^{pn_2''}+p^{-1}\sum _{k=1}^{p-1}{p \choose k}X_1^{pk}(-X_1)^{p-k}$$ 

where $f(T)=\sum _{j\in J}a_jT^j\in R[[T]]\{T^{-1}\}$ has conductor $m_2'$ and the following holds: either $pn_2''-n_1(p-1)<n_1(p(p-1)+1)$ or
$pn_2''-n_1(p-1)=n_1(p(p-1)+1)$ and $m_2'\ge m_1p$. Moreover we can assume after some transformation which eliminate the $p$-powers that
$f(T)=\sum _{j\ge m_2'} a_jT^j$. Let $\Cal X_1\to \Cal X$ be the blow up of $\Cal X$ at the point $x$
and let $X_1$ be the exeptional fibre in $\Cal X_1$ which meets $X$ at the double point $x=x_{i_0}$. Let $e=p^2t$ be the integer associated to the marked double 
point $x$ via $\Deg.6$. After enlarging $R$ we can assume that the double point $x$ of $\Cal X_1$ has thikness $e$.
Let $\Cal X_1'\simeq \Spf R<S^{-1}>$  be the formal fibre of $X_1':=X_1-\{x\}$ in $\Cal X_1$. Let $f_1':\Cal Y_1'\to \Cal X_1'$ be the 
\'etale $\Bbb Z/p^2\Bbb Z$-torsor given by the equation: 
$$X_1^p-X_1=S^{m_1}$$
and
$$X_2^p-X_2=h'(S)+p^{-1}\sum _{k=1}^{p-1}{p \choose k}X_1^{pk}(-X_1)^{p-k}$$ 

where $h'(S)=f(\pi^{p^2t}S^{-1})/\pi ^{pn_2''}\in R<S^{-1}>$. 

Further let $\Cal X_{1,x}\simeq \Spf R[[S,T]]/(ST-\pi^{p^2t})$ be the formal germ of 
$\Cal X_1$ at the double point $x$. Consider the cover $f_{1,x}:\Cal Y_{1,y} \to \Cal X_{1,x}$ given by the equation:
$$X_1^p-X_1=T^{-m_{1}}/\pi ^{n_1p}$$ 
and: 
$$X_2^p-X_2=f(T)/\pi^{pn_2''}+p^{-1}\sum _{k=1}^{p-1}{p \choose k}X_1^{pk}(-X_1)^{p-k}$$ 
hen $\Cal Y_{1,y}$ is the formal germ of a double point of thikness $t$ (cf. 4.2.8). Moreover the cover $f_{1,x}$ restricted to the boundary
$\Spf R[[S]]\{S^{-1}\}$ is given by the equations:
$$X_1^p-X_1=S^{-m_1}$$
and
$$X_2^p-X_2=h'(S)+p^{-1}\sum _{k=1}^{p-1}{p \choose k}X_1^{pk}(-X_1)^{p-k}$$ 

In particular the cover $f_1'$ (resp. $f'$) has the same degeneration type, 
by construction, on the boundary corresponding to the double point $x$ as the degeneration type of the cover $f_{1,x}$ above the formal boundary with parameter 
$T$ (resp. above the formal boundary with parameter $S$). A second application of the formal patching techniques allows one to construct a $p^2$-cyclic cover
$f_1:\Cal Y_1\to \Cal X_1$ which restricted to $\Cal X'$ (resp. $\Cal X_1'$ and $\Cal X_{1,x}$) is isomorphic to $f'$ (resp. to $f_1'$ and $f_x$).
Let $\tilde {\Cal X}$ be the $R$-curve obtained by contracting the irreducible component $X$ in $\Cal X_1$. We denote the image of the double point $x$ in
$\tilde {\Cal X}$ simply by $x$. The cover $f_1:\Cal Y_1\to \Cal X_1$ induces canonically a $p^2$-cyclic cover $\tilde f:\tilde {\Cal Y}\to \tilde {\Cal X}$ 
above $\tilde {\Cal X}$. Let $\Cal X_x\simeq R[[T]]$ be the formal
germ of $\tilde {\Cal X}$ at the smooth point $x$. Then $\tilde f$ induces canonically a $p^2$-cyclic cover $f_{x}:\Cal Y_{y} \to \Cal X_{x}$ where 
$y$ is the closed point of $\Cal Y$ above $x$. Now it is easy to see that the degeneration data associated to $f_x$ via 6.2.2 is isomorphic to the degeneration 
data $\Deg (x)$ we started with. Now the proof in the general case is very similar and is left to the reader. The only modification is that one has to consider
$p^2$-cyclic covers $f_i:\Cal Y_i\to \Cal X_i$ above the formal germs $\Cal X_i:=\Spf R[[T_i]]$ which one obtains by induction hypothesis as realisation of the 
degeneration data, induced by $\Deg (x)$, on the subtrees $\Gamma_i$ of $\Gamma$ which starts from the edge $x_i$ going to the ends and which have clearly 
length strictly less than the length of $\Gamma$. For this one still has to use the data in
$\Deg.8$ to insure that there exists a lifting $f'$ as above of the degeneration data $\Cal D$ above $U$ which restricted to the formal boundary corresponding
to the point $x_i$ is Galois isomorphic to the restriction of $f_i$ to the formal boundary $\Spf R[[T_i]]\{T_i^{-1}\}$.

\pop {.5}
\par
\noindent
{\bf \gr 6.2.4. Remarks.}

\pop {.5}
\par
\noindent
{\bf \gr 1.}\rm\ The data $\Deg.8$ in 6.2.1 seems necessary for the proof of 6.2.3. Indeed given two $p^2$-cyclic covers above the two boundaries
of the formal germ $\Cal X$ of an $R$-curve at a double point with ``compatible'' degeneration type as in 4.2.11 it is not always possible to find a 
$p^2$-cyclic cover $f:\Cal Y\to \Cal X$ where $\Cal Y$ is the forml germ of a double point and $f$ is \'etale above $\Cal X_K$ and which restricts
to the above given covers above the boundaries. Note that this is possible for $p$-cyclic covers.

\pop {.5}
\par
\noindent
{\bf \gr 2.}\rm\ Using the same techniques as in 6.2.3 it is possible to prove global results of lifting for finite covers between proper and semi-stable
$k$-curves of degree $p^2$ to a $p^2$-cyclic cover between $R'$-curve over some finite extension $R'$ of $R$.

\pop {2}
\par
\noindent
{\bf \gr References}
\rm

\pop {.5}
\par
\noindent
[B-L] S. Bosch and W. L\"utkebohmert, {\sl Formal and rigid geometry}, Math. Ann, 295, 291-317, (1993).

\pop {.5}
\par
\noindent
[De-Mu] P. Deligne and D. Munford, {\sl The irreducibility of the space of 
curves with given genus}, Publ. Math. IHES. 36, 75-109, (1969).

\pop {.5}
\par
\noindent
[D-G] M. Demazure and P. Gabriel, Groupes Alg\'ebriques, Tome 1, Masson and CIE \'Editeur, Paris, North-Holland Publishing Company, 
Amsterdam, (1970).

\pop {.5}
\par
\noindent
[E] H. Epp, {\sl Eliminating wild ramification}, Invent. Math. 19, 235-249, (1973).

\pop {.5}
\par
\noindent
[G-Ma] B. Green and M. Matignon, {\sl Order $p$-automorphisms of the open disc of a $p$-adic field}, 
J. Amer. Math. Soc. 12, 269-303, (1999).

\pop {.5}
\par
\noindent
[H] Y. Henrio, {\sl arbres de Hurwitz et automorphismes d'ordre $p$ de disques et couronnes $p$-adique formels}, 
PHD Thesis, university of Bordeaux, (1999).

\pop {.5}
\par
\noindent
[L] C. Lehr, {\sl Reduction of wildly ramified covers of curves}, PHD Thesis, University of Pennsylvania, (2001).

\pop {.5}
\par
\noindent
[Ma] M. Matignon, {\sl Vers un algorithme pour la r\'eduction stable des rev\^etements $p$-cycliques de la 
droite projective sur un corps $p$-adique.}, Math. Ann. 325, 323-354, (2003).

\pop {.5}
\par
\noindent
[Ma-Y] M. Matignon and T. Youssefi, {\sl Prolongements de morphismes entre fibres formelles et cycles evanecsents}, preprint, (1992).

\pop {.5}
\par
\noindent
[M] S. Maugeais, {\sl Rel\`evement des rev\^etements $p$-cycliques des courbes rationnelles semi-stable},
Math. Ann., 327, 365-393, (2003).

\pop {.5}
\par
\noindent
[R] M. Raynaud, {\sl $p$-groupes et r\'eduction semi-stable des courbes}, The Grothendieck 
\newline
Festschrift, Vol III, 
Progress in Mathematics 88 (1990), Birkh\"auser, 179-197.

\pop {.5}
\par
\noindent
[R-1] M. Raynaud, {\sl Sp\'ecialisation des rev\^etements en caract\'ertistique $p>0$}, Ann. Scient. \'Ec. Norm. Sup. 32, 
87-126, (1999).

\pop {.5}
\par
\noindent
[S] M. Sa\"\i di, {\sl Torsors under finite and flat group schemes of rank $p$ with Galois action}, Math. Z.,
245, 695-710, (2003).

\pop {.5}
\par
\noindent
[S-1] M. Sa\"\i di, {\sl Wild ramification and a vanishing cycles formula}, J. of Algebra, 273, 1, 108-128, (2004).

\pop {.5}
\par
\noindent
[S-2] M. Sa\"\i di, {\sl Rev\^etements mod\'er\'es et groupe fondamental de graphe de groupes}, Compositio Math., 107, 319-338,
(1997).

\pop {.5}
\par
\noindent
[S-3]  M. Sa\"\i di, {\sl Galois covers of degrre $p$: Galois action and semi-stable reduction}, 
\newline
math.AG/0106249.
\pop {.5}
\par
\noindent
[SGA-1] A. Grothedieck, Rev\^etements \'etales et groupe fondamental, Lect. Notes in Math. 224, Springer-Verlag, (1971). 

\pop {.5}
\par
\noindent
[We] S. Wewers, {\sl Three point covers with bad reduction}, J. Amer. Math. Soc., 16, 991-1032, (2003).

\pop {.5}
\par
\noindent
[W] E. Witt, {\sl Zyclische K\"orper und algebren der characteristic $p$ vom 
grad $p^n$. Struktur diskret bewerteter perfekter K\"orper mit vollkommenem 
Restklassenk\"orper der Characteristic $p$}, J. F\"ur die reine und angewandte
Mathematik (Crelle), 176, 126-140, (1937).

\enddocument
\end